\documentclass[10pt]{article}
\usepackage{amsfonts,amsmath,mathrsfs,amssymb}
\usepackage{tikz}
\usepackage{pxfonts}
\usepackage{verbatim}
\usepackage{array}
\usepackage{epsfig}
\usepackage{tabularx}
\usepackage[hang]{subfigure}

\usepackage[round]{natbib}

\newtheorem{theorem}{Theorem}[section]

\newtheorem{defin}{Definition}[section]

\newtheorem{remark}{Remark}[section]
\newtheorem{cor}[theorem]{Corollary}

\def\R{ \mathbb{R}}
\let\tilde\widetilde

\let\bar\overline

\def\eop{\hfill $\square$}

\definecolor{revision}{rgb}{1,0,0}

\def\complaint#1{}
\def\withcomplaints{
\newcounter{mycomplaints}
\def\complaint##1{\refstepcounter{mycomplaints}%
\ifhmode%
\unskip%
{\dimen1=\baselineskip \divide\dimen1 by 2 %
\raise\dimen1\llap{\tiny -\themycomplaints-}}\fi%
\marginpar{\tiny [\themycomplaints]: \textcolor{revision} {##1}}}%
}

\withcomplaints 

\newcommand{\Z}{\mathbb{Z}}

\newcommand{\bt}{{\bf t}}
\newcommand{\p}{{p}}

\newcommand{\bu}{{u}}

\newcommand{\wt}{\widetilde}

\newcommand{\pfw}{(\langle \wt G, L \rangle , \wt \p)} 
\newcommand{\pofw}{(\langle G, n \rangle,\p)} 
\newcommand{\pog}{\langle G, n \rangle} 
\newcommand{\spog}{\langle G, g \rangle} 

\newcommand{\grey}{\textcolor{gray}}
\definecolor{pink}{rgb}{1, .509, .8}
\definecolor{sky}{rgb}{.4, .914, 1}

\begin{document}

\title {Finite motions from periodic frameworks with added symmetry}

\author{Elissa Ross\footnote{supported in part under a grant from NSERC Canada} \\
 Department of Mathematics and Statistics\\ York University\\ 4700 Keele Street\\ Toronto, ON M3J1P3, Canada
\\Bernd Schulze\footnote{supported by the DFG Research Unit 565 `Polyhedral
Surfaces'}\\ Institute of Mathematics, MA 6-2, TU Berlin,\\ Strasse des 17.
Juni 136, D-10623 Berlin, Germany \\
and Walter Whiteley\footnote{Supported by a grant from NSERC (Canada).}\\
 Department of Mathematics and Statistics\\ York University\\ 4700 Keele Street\\ Toronto, ON M3J1P3, Canada
}

\maketitle

\begin{abstract}
Recent work from authors across disciplines has made substantial contributions to counting rules (Maxwell type theorems) which predict when an infinite  periodic structure would be rigid or flexible while preserving the periodic pattern,  as an engineering type framework, or equivalently, as an idealized molecular framework.    Other work has shown that for finite frameworks, introducing symmetry modifies the previous general counts, and under some circumstances this symmetrized Maxwell type count can predict added finite flexibility in the structure.

In this paper we combine these approaches to present new Maxwell type counts for the columns and rows of a modified orbit matrix for structures that have both a periodic structure and additional symmetry within the periodic cells.  In a number of cases, this count for the combined group of symmetry operations demonstrates there  is added finite flexibility in what would have been rigid when realized without the symmetry.   Given that many crystal structures have these added symmetries, and that their flexibility may be key to their physical and chemical properties, we present a summary of the results as a way to generate further developments of both a practical and theoretic interest.
\end{abstract}

\section{Introduction}
\label{sec:Introduction}

The theory of periodic frameworks has undergone rapid and extensive development in the last four years \citep{BorceaStreinuII, BorceaStreinuI, Theran, Ross}.   We now have necessary conditions (call them Maxwell counts) for such frameworks to be rigid, either with a fixed lattice of translations or with a flexible lattice of translations.   Underlying much of the recent work are finite `lattice rigidity matrices' for the equivalence classes of vertices and edges under the infinite group of translations $\Z^{d}$ in $d$-space.   With the corresponding count of periodicity-preserving trivial motions under these constraints (typically $d$ translations), the number of rows, $e$, and columns, $dv+l$ (where $l$ is the number of lattice parameters) of these `orbit matrices' lead to necessary Maxwell type counts for a framework to be infinitesimally rigid \citep{BorceaStreinuI,
Theran, Ross}:  $e\geq dv + l -d$.

The theory of finite symmetric frameworks has also experienced some breakout results, building on a decade or more of initial Maxwell type necessary conditions for frameworks of various symmetry groups \citep{FGsymmax,FG4, cfgsw}.  In some key cases, these symmetry conditions predict finite motions for frameworks realized generically within the symmetry constraints, but whose graphs would be generically rigid without  symmetry \citep{KG1,bricard}.   Recently, key results of this work have been expressed in terms of `orbit rigidity matrices' for the equivalence classes of vertices and edges under the group of symmetry operations ${\cal S}$ \citep{BS4, BS6, BSWWorbit}.  With modified counts for the symmetry-preserving trivial motions $t_{\cal S}$, and with $e_0$ and $v_0$ denoting the number of edge orbits and vertex orbits under the group action of $\mathcal{S}$, respectively, these matrices lead to Maxwell type necessary counts for frameworks to be infinitesimally rigid: $e_{0} \geq d v_{0} - t_{\cal S}$.

Given that many crystal structures combine both periodic structure and symmetry within the unit cells, it is natural to investigate the interactions of these two types of group operations. So we will consider frameworks with `combined symmetry groups' of the form $\Z^{d}\rtimes \cal S $, where $\Z^{d}$ is the group of translations of the framework, $\cal S $ is the group of additional symmetries of the framework, and $\rtimes$ denotes the semi-direct product of $\cal S $ acting on $\Z^{d}$. Note that every symmetry operation in such a group can be written as a unique product of an element of $\Z^{d}$ and an element of $\cal{S}$. However, since $\cal S $ is typically not normal in $\Z^{d}\rtimes \cal S $, the groups $\Z^{d}\rtimes \cal S $ are in general not direct products. Details on the semi-direct product can be found in any abstract algebra text, such as \citet{DummitFoote}. In Section \ref{sec:PlanePeriodicSymmetry} we will introduce combined `orbit matrices' for the groups $\Z^{d}\rtimes \cal S $.  Combined with the counts of the trivial motions which preserve both symmetry and periodicity, this will provide extended Maxwell type necessary counts for infinitesimal rigidity.
In this setting we:
\begin{enumerate}
\item  count the rows of the combined orbit matrix: one row per orbit of edges $r=e_{0}$;
\item  count the columns of the combined orbit matrix: one vector column per orbit of vertices plus columns for symmetry-preserving lattice deformations: $c = dv_{0}+\ell_{\cal S}$;
\item  the dimension of the space of trivial motions (translations) left by symmetries:  $t_{\cal S}$.
\end{enumerate}
The minimum dimension of the space of non-trivial symmetry-preserving infinitesimal periodic motions of the periodic structure is:
$$  m= c - t_{\cal S}\ \  -r \quad \text{or}  \quad  m= dv_{0} +\ell_{\cal S}  - t_{\cal S}\ \ -e_{0}.$$
This is compared with the corresponding count on the graph without symmetry, where with orbits of size $k_{{\cal S}}$ and no fixed edges or vertices, for the fully flexible lattice, we would anticipate:  $$m= d(k_{{\cal S}}v_{0}) + {d+1 \choose 2} -d \ \ - (k_{{\cal S}}e_{0}).$$
 In addition, if we choose the positions of the vertices generically within the symmetry (i.e., make one generic choice for each orbit of vertices) then the predicted infinitesimal motions will be finite flexes \citep{asiroth, BS6, BSWWorbit}.

The results are a surprise -- adding symmetry can sometimes cause additional flexibility beyond what the original graph without symmetry would exhibit in the periodic lattice.   These more flexible examples include symmetries such as inversive symmetry, or half-turn symmetry with a mirror, found in a number of crystals, such as zeolites.  Recent studies have confirmed that flexibility is a feature of natural zeolites \citep{Thorpe} and contributes to their physical and chemical properties. In turn, this suggests that predicted flexibility in a computer designed theoretical \ `zeolite' would be a criterion for selecting which theoretical compounds should be synthesized for further testing.

When adding symmetry to a periodic lattice structure, we must consider the flexibility that this symmetry allows in the lattice structure.  Inversive symmetry will be a key example, as it fits all possible lattice deformations (it occurs in `triclinic lattices'), and the addition of this symmetry to the framework generates flexes from frameworks that previously were minimally rigid, while preserving the full range of possible flexes of the lattice itself.

In contrast, only certain types of lattices leave open the addition of a half-turn symmetry in $3$-space.  A half-turn parallel to a side of the lattice requires that side to be perpendicular to the remaining  parallelogram face.  This leaves only four of the six possible flexes of the lattice (monoclinic lattices), but it does predict additional flexes.    Similarly, mirrors of symmetry can fit parallel to faces of the lattice, and restrict the shapes to monoclinic lattices, with the variable angle now parallel to the mirror.

We can also have a larger symmetry group, with several generators.  For example, monoclinic prismatic crystals, such as some forms of zeolite, have the symmetry group $\mathcal{C}_{2h}$ which has both a half-turn symmetry and a mirror.  These restrict the possible lattice shapes to lattices with a parallelogram base (perpendicular to the axis, parallel to the mirror) and the vertical prism at right angles to the base.   Other forms  of zeolite have the added symmetry of $D_{2h}$, forcing the base to be a rectangle.  Each symmetry group for the crystal structure and the associated crystal system requires some  specific terms in the analysis, although patterns emerge, and we will present tables with rows for the combinations.

In the larger theory of rigidity of frameworks, infinitesimal flexes of `generic frameworks' transfer to finite flexes, for appropriate versions of generic.  This holds for generic frameworks without symmetry, for frameworks generic within the symmetry class, and for periodic frameworks.   That property extends to these combined symmetry periodic frameworks, so we are talking about flexibility on a finite scale, at generic realizations for representatives of the orbits of the expanded group.

The outline of the paper is as follows.   In \S2, \S3, and \S4, we present the basic definitions and associated rigidity matrices for: (i) finite frameworks (\S2); (ii) finite frameworks with symmetry and the orbit matrices (\S3); and (iii) periodic frameworks with the associated lattice matrices (\S4).

In \S5, we introduce our method of inserting symmetry into the analysis of a periodic framework in dimensions $2$ and $3$.   This includes a short summary of the wallpaper pattern types (dimension $2$) and  crystal systems (dimension $3$) which arise as the possible lattices and restricted lattice variations for various symmetry groups in the corresponding dimension.  \S6 then presents some key examples in the plane for the groups $\Z^{2}\rtimes {\cal C}_{2}$ and $\Z^{2}\rtimes {\cal C}_{s}$ with various lattice flexibilities, as well as summary tables over all plane groups of the form
$\Z^{2}\rtimes {\cal S}$.  \S7 presents key examples  in $3$-space, with corresponding tables.

In \S8 we briefly describe a range of extensions which are accessible using these methods, including:  extensions to include fixed vertices and edges; extensions to additional plane and space groups; extensions to higher dimensions; and the companion static analysis of the frameworks.

The analysis is not complete.  While we have covered groups of the form $\Z^{d}\rtimes {\cal S}$,  some plane and space groups are not covered, namely those with 6-fold symmetry, or with glide reflections.    As \S8 illustrates, there is lots of room for additional exploration.  We hope that this introduction, and the follow-up papers with detailed proofs for what is claimed here, will offer interested researchers tools to explore the  range of examples of interest in their context.

\section{Preliminaries on the rigidity of (finite) frameworks}
\label{sec:Preliminaries}

\indent A \emph{framework}  in $\mathbb{R}^{d}$ is a pair $(G,p)$, where $G$ is a finite simple (no loops or multiple edges) graph with vertex set $V(G)$ and edge set $E(G)$, and  $p:V(G)\to \mathbb{R}^d$ is a map such that $p(i) \neq p(j)$ for all $\{i,j\}\in E(G)$.  We also say that $(G,p)$ is a $d$-dimensional \emph{realization} of the \emph{underlying graph} $G$ \citep{W1}. For $i\in V(G)$, we say that $p_i:=p(i)$ is the \emph{joint} of $(G,p)$ corresponding to $i$, and for $\{i,j\}\in E(G)$, we say that $p_{\{i,j\}}:=p(\{i,j\})$ is the \emph{bar} of $(G,p)$ corresponding to $\{i,j\}$. Moreover, we let $v:=|V(G)|$ and $e:=|E(G)|$.
It is often useful to identify $p$ with a vector in $\mathbb{R}^{dv}$ by using the order on $V(G)$. In this case we also refer to $p$ as a \emph{configuration} of $v$ points in $\mathbb{R}^{d}$.

A framework $(G,p)$ in $\mathbb{R}^d$ is \emph{flexible} if there exists a continuous path, called a \emph{finite flex} or \emph{mechanism}, $p(t):[0,1]\to \mathbb{R}^{dv}$  such that \begin{itemize}
\item[(i)] $p(0)=p$;
\item[(ii)]  $\|p(t)_i-p(t)_j\|=\|p_i-p_j\|$ for  all $0\leq t\leq 1$ and all $\{i,j\}\in E(G)$;
\item[(iii)]  $\|p(t)_k-p(t)_l\|\neq\|p_k-p_l\|$ for  all $0< t\leq 1$ and some pair $\{k,l\}$ of vertices of $G$. \end{itemize}
 Otherwise $(G,p)$ is said to be \emph{rigid}. For some alternate equivalent definitions of a rigid and flexible framework see \citet{asiroth}, for example.

An \emph{infinitesimal motion} of a framework $(G,p)$ in $\mathbb{R}^d$ is a function $u: V(G)\to \mathbb{R}^{d}$ such that
\begin{equation}
\label{infinmotioneq}
(p_i-p_j)\cdot (u_i-u_j)=0 \quad\textrm{ for all } \{i,j\} \in E(G)\textrm{,}\end{equation}
where $u_i$ denotes the vector $u(i)$ for each $i$.

An infinitesimal motion $u$ of $(G,p)$ is an \emph{infinitesimal rigid motion} (or \emph{trivial infinitesimal motion}) if there exists a skew-symmetric matrix $S$ (a rotation) and a vector $t$ (a translation) such that $u_i=Sp_i+t$ for all $i\in V(G)$. Otherwise $u$ is an \emph{infinitesimal flex} (or \emph{non-trivial infinitesimal motion}) of $(G,p)$.

$(G,p)$ is \emph{infinitesimally rigid} if every infinitesimal motion of $(G,p)$ is an infinitesimal rigid motion. Otherwise $(G,p)$ is said to be \emph{infinitesimally flexible} \citep{W1}.

While an infinitesimally rigid framework is always rigid, the converse does not hold in general. \citet{asiroth} however, showed that for `generic' configurations, infinitesimal rigidity and rigidity are in fact equivalent.

The \emph{rigidity matrix} of $(G,p)$ (which in structural engineering is also known as the \emph{compatibility matrix} of $(G,p)$ \citep{KG2, cfgsw}) is the $e \times dv$ matrix
\begin{displaymath} \mathbf{R}(G,p)=\bordermatrix{& & & & i & & & & j & & & \cr & & & &  & & \vdots & &  & & &
\cr \{i,j\} & 0 & \ldots &  0 & (p_{i}-p_{j}) & 0 & \ldots & 0 & (p_{j}-p_{i}) &  0 &  \ldots&  0 \cr & & & &  & & \vdots & &  & & &
}
\textrm{,}\end{displaymath}
that is, for each edge $\{i,j\}\in E(G)$, $\mathbf{R}(G,p)$ has the row with
$(p_{i}-p_{j})_{1},\ldots,(p_{i}-p_{j})_{d}$ in the columns $d(i-1)+1,\ldots,di$, $(p_{j}-p_{i})_{1},\ldots,(p_{j}-p_{i})_{d}$ in
the columns $d(j-1)+1,\ldots,dj$, and $0$ elsewhere \citep{W1}. See also Example 3.1.

Note that if we identify an infinitesimal motion $u$ of $(G,p)$ with a column vector in $\mathbb{R}^{dv}$ (by using the order on $V(G)$), then the equations in (\ref{infinmotioneq}) can be written as $\mathbf{R}(G,p)u=0$. So, the kernel of the rigidity matrix $\mathbf{R}(G,p)$ is the space of all infinitesimal motions of $(G,p)$. It is well known that the infinitesimal rigid motions arising from $d$ translations and $\binom{d}{2}$ rotations of $\mathbb{R}^{d}$ form a basis of the space of infinitesimal rigid motions of $(G,p)$, provided that the points $p_{1},\ldots ,p_{v}$ span an affine subspace of $\mathbb{R}^{d}$ of dimension at least $d-1$ \citep{W1}. Thus, for such a framework $(G,p)$, we have $\textrm{nullity }\big(\mathbf{R}(G,p)\big)\geq d+\binom{d}{2}=\binom{d+1}{2}$ and $(G,p)$ is infinitesimally rigid if and only if $\textrm{nullity } \big(\mathbf{R}(G,p)\big)=\binom{d+1}{2}$ or equivalently, $\textrm{rank }\big(\mathbf{R}(G,p)\big)=d v - \binom{d+1}{2}$.

In particular, it follows that we can sometimes detect infinitesimal flexes in frameworks - and, by the result of \citet{asiroth}, even predict \emph{finite} flexes in  generic frameworks - by simply counting vertices and edges:

\begin{theorem} [Maxwell's rule] \label{maxtheorem} Let $(G,p)$ be a $d$-dimensional framework whose joints span an affine subspace of $\mathbb{R}^{d}$ of dimension at least $d-1$. If
\begin{equation}\label{maxrule} e<dv-\binom{d+1}{2}\textrm{,}
\end{equation}
then $(G,p)$ has an infinitesimal flex.

If the joints of $(G,p)$ are in generic position, then there even exists a finite flex of $(G,p)$.
\end{theorem}

\section{Symmetry in frameworks}
\label{sec:OrbitMatrices}

\subsection{Symmetric frameworks and motions}
\label{subsec:SymmetricOrbitMatrices1}

Given a finite simple graph $G$ with vertex set $V(G)=\{1,\ldots,n\}$, and a map $p:V(G)\to \mathbb{R}^d$, a \emph{symmetry operation} of the framework $(G,p)$ in $\mathbb{R}^{d}$ is an isometry $s$ of $\mathbb{R}^{d}$ such that for some $\alpha_s\in \textrm{Aut}(G)$, we have
\begin{equation} \label{eq:symop} s(p_i)=p_{\alpha_s(i)}\quad \textrm{for all } i\in V(G)\textrm{, }\nonumber\end{equation} where $\textrm{Aut}(G)$ denotes the automorphism group of the graph $G$ \citep{Hall, BS2, BS1}.  The set of all symmetry operations of a framework $(G,p)$ forms a group under composition, called the \emph{point group} of $(G,p)$ \citep{bishop, Hall, BS2, BS1}. Since translating a framework does not change its rigidity properties, we may assume wlog that the point group of a framework is always a \emph{symmetry group}, i.e., a subgroup of the orthogonal group $O(\mathbb{R}^{d})$ \citep{BS2, BS1}.

Throughout this paper, we will highlight the Schoenflies notation for the symmetry operations and symmetry groups, as this is one of the standard notations in the literature (see \citet{bishop, cfgsw, FGsymmax, FG4, Hall, KG1, KG2, BS2, BS4, BS1}, for example). In the later tables for crystallographic groups, we will show three notations in parallel, to ensure clearer communication with multiple audiences. 

In the Schoenflies notation, the groups we will focus on in our examples and tables are denoted by $\mathcal{C}_s$, $\mathcal{C}_n$, $\mathcal{C}_{nv}$, $\mathcal{C}_{nh}$, $\mathcal{C}_i$, $\mathcal{D}_{n}$, and $\mathcal{D}_{nh}$. For dimension $2$ and $3$, $\mathcal{C}_s$ is a symmetry group consisting of the identity $Id$ and a single reflection $s$, and $\mathcal{C}_n$ is a cyclic group generated by an $n$-fold rotation $C_n$. The only other possible type of symmetry group in dimension $2$ is the group $\mathcal{C}_{nv}$ which is a dihedral group generated by a pair $\{C_n, s\}$. In dimension $3$, $\mathcal{C}_{nv}$ denotes any symmetry group that is generated by a rotation $C_n$ and a reflection $s$ whose corresponding mirror contains the rotational axis of $C_n$, whereas a symmetry group $\mathcal{C}_{nh}$ is generated by a rotation $C_n$ and the reflection $s$ whose corresponding mirror is perpendicular to the $C_n$-axis. The group $\mathcal{C}_i$ consists of the identity $Id$ and an inversion $i$ in $3$-space. Finally, $\mathcal{D}_{n}$ is generated by an $n$-fold rotation $C_n$ and a $2$-fold rotation $C_2$ whose rotational axes are perpendicular to each other, and $\mathcal{D}_{nh}$ is generated by the generators  $C_n$ and $C_2$ of a group $\mathcal{D}_{n}$ and by a reflection $s$ whose mirror is perpendicular to the $C_n$-axis.

Given a symmetry group $\mathcal{S}$ in dimension $d$ and a graph $G$, we let $\mathscr{R}_{(G,\mathcal{S})}$ denote the set of all $d$-dimensional realizations of $G$ whose point group is either equal to $\mathcal{S}$ or contains $\mathcal{S}$ as a subgroup \citep{BS2, BS1}. In other words, the set $\mathscr{R}_{(G,\mathcal{S})}$ consists of all frameworks $(G,p)$ for which there exists a map $\phi:\mathcal{S}\to \textrm{Aut}(G)$ so that
\begin{equation}\label{class} s\big(p_i\big)=p_{\Phi(s)(i)}\textrm{ for all } i\in V(G)\textrm{ and all } s\in \mathcal{S}\textrm{.}\end{equation}
If a framework $(G,p)\in \mathscr{R}_{(G,\mathcal{S})}$ satisfies the equations in (\ref{class}) for the map $\Phi:\mathcal{S}\to \textrm{Aut}(G)$, we say that $(G,p)$ is \emph{of type $\Phi$}.
It is shown in \citet{BS4, BS1} that if the map $p$ of a framework  $(G,p)\in \mathscr{R}_{(G,\mathcal{S})}$ is injective, then $(G,p)$ is of a unique type $\Phi$ and $\Phi$ is necessarily also a homomorphism. For simplicity, we therefore assume that the map $p$ of any framework $(G,p)$ considered in this paper is injective (i.e., $p_i\neq p_j$ if $i\neq j$). In particular, this allows us (with a slight abuse of notation) to use the terms $p_{s(i)}$ and  $p_{\Phi(s)(i)}$ interchangeably, where $i\in V(G)$ and $s\in \mathcal{S}$. In general, if the type $\Phi$ is clear from the context, we often simply write $s(i)$ instead of $\Phi(s)(i)$.

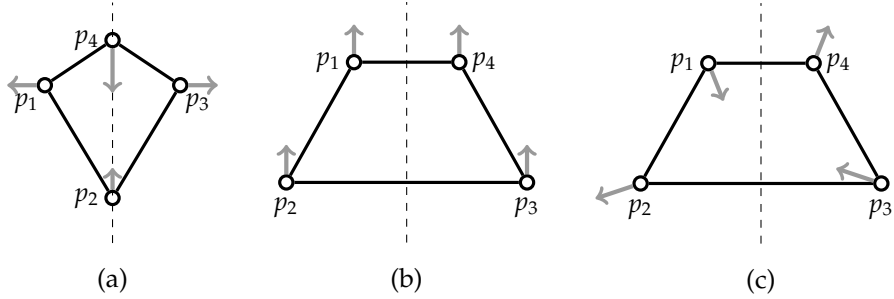
\begin{figure}[htp]
\begin{center}
\begin{tikzpicture}[very thick,scale=1]
\tikzstyle{every node}=[circle, draw=black, fill=white, inner sep=0pt, minimum width=5pt];
        \path (1.5,-1) node (p2)[label = left: $p_{2}$] {} ;
       \path (0.6,0.5) node (p1)  [label = below left: $p_{1}$]{} ;
   \path (2.4,0.5) node (p3) [label = below right: $p_{3}$] {} ;
   \path (1.5,1.1) node (p4) [label = left: $p_{4}$] {} ;
   \draw (p1) -- (p4);
      \draw (p3) -- (p4);
     \draw (p2) -- (p3);
      \draw (p2) -- (p1);
             \draw [dashed, thin] (1.5,-1.6) -- (1.5,1.6);
              \draw [ultra thick, ->, black!40!white] (p1) -- (0.1,0.5);
      \draw [ultra thick, ->, black!40!white] (p3) -- (2.9,0.5);
      \draw [ultra thick, ->, black!40!white] (p4) -- (1.5,0.4);
      \draw [ultra thick, ->, black!40!white] (p2) -- (1.5,-0.6);
      \node [draw=white, fill=white] (a) at (1.5,-2.1) {(a)};
    \end{tikzpicture}
    \hspace{0.5cm}
            \begin{tikzpicture}[very thick,scale=1]
\tikzstyle{every node}=[circle, draw=black, fill=white, inner sep=0pt, minimum width=5pt];
    \path (-0.7,0.8) node (p1) [label = left: $p_{1}$] {} ;
    \path (0.7,0.8) node (p4) [label = right: $p_{4}$]{} ;
    \path (-1.6,-0.8) node (p2) [label = below: $p_{2}$] {} ;
     \path (1.6,-0.8) node (p3) [label = below: $p_{3}$] {} ;
      \draw (p1) -- (p4);
    \draw (p1) -- (p2);
    \draw (p3) -- (p4);
    \draw (p2) -- (p3);
     \draw [dashed, thin] (0,-1.6) -- (0,1.6);
     \draw [ultra thick, ->, black!40!white] (p1) -- (-0.7,1.3);
      \draw [ultra thick, ->, black!40!white] (p4) -- (0.7,1.3);
      \draw [ultra thick, ->, black!40!white] (p2) -- (-1.6,-0.3);
      \draw [ultra thick, ->, black!40!white] (p3) -- (1.6,-0.3);
      \node [draw=white, fill=white] (b) at (0,-2.1) {(b)};
        \end{tikzpicture}
        \hspace{0.5cm}
        \begin{tikzpicture}[very thick,scale=1]
\tikzstyle{every node}=[circle, draw=black, fill=white, inner sep=0pt, minimum width=5pt];
    \path (-0.7,0.8) node (p1) [label = left: $p_{1}$]  {} ;
    \path (0.7,0.8) node (p4)[label = right: $p_{4}$] {} ;
    \path (-1.6,-0.8) node [label = below: $p_{2}$](p2)  {} ;
     \path (1.6,-0.8) node [label = below: $p_{3}$](p3)  {} ;
      \draw (p1) -- (p4);
    \draw (p1) -- (p2);
    \draw (p3) -- (p4);
    \draw (p2) -- (p3);
     \draw [dashed, thin] (0,-1.6) -- (0,1.6);
     \draw [ultra thick, ->, black!40!white] (p1) -- (-0.5,0.3);
      \draw [ultra thick, ->, black!40!white] (p4) -- (0.9,1.3);
      \draw [ultra thick, ->, black!40!white] (p2) -- (-2.2,-1);
      \draw [ultra thick, ->, black!40!white] (p3) -- (1,-0.6);
      \node [draw=white, fill=white] (b) at (0,-2.1) {(c)};
        \end{tikzpicture}
\end{center}
\vspace{-0.3cm}
\caption{Infinitesimal motions of frameworks in the plane: (a) a $\mathcal{C}_s$-symmetric infinitesimal flex; (b) a $\mathcal{C}_s$-symmetric infinitesimal rigid motion; (c) an infinitesimal flex which is not $\mathcal{C}_s$-symmetric.}
\label{fulsym}
\end{figure}

An infinitesimal motion $u$ of a framework $(G,p)\in \mathscr{R}_{(G,\mathcal{S})}$  is \emph{$\mathcal{S}$-symmetric} if \begin{equation}\label{fulsymmot} s\big(u_i\big)=u_{s(i)}\textrm{ for all } i\in V(G)\textrm{ and all } s\in \mathcal{S}\textrm{,}\end{equation} i.e., if $u$ is unchanged under all symmetry operations in $\mathcal{S}$ (see also Figure \ref{fulsym}(a) and (b)).

Note that if for $(G,p)\in \mathscr{R}_{(G,\mathcal{S})}$, we choose a set of representatives $\{1,\ldots, v_0\}$ for the orbits  $\mathcal{S}(i)=\{s(i)|\,s\in \mathcal{S}\}$ of vertices of $G$ under the group action of $\mathcal{S}$, then the positions of \emph{all} joints of  $(G,p)$ are uniquely determined by the positions of the joints $p_1,\ldots,p_{v_0}$ and the symmetry constraints imposed by $\mathcal{S}$. Similarly, an $\mathcal{S}$-symmetric infinitesimal motion $u$ of $(G,p)$ is uniquely determined by the velocity vectors $u_1,\ldots, u_{v_0}$ for the representative vertices.

The following extension of the theorem of \citet{asiroth} shows that an analysis of the `$\mathcal{S}$-symmetric' infinitesimal rigidity properties of a symmetric framework can be used to also detect \emph{finite} flexes in the framework, provided that its joints are positioned generically within the symmetry.

\begin{theorem}
\label{thm:flexes} \citep{BS6} (see also \citet{BS4})
Let $\mathcal{S}$ be a symmetry group in dimension $d$, and let $(G,p)\in \mathscr{R}_{(G,\mathcal{S})}$ be a framework whose joints span all of $\mathbb{R}^d$, in an affine sense. If $(G,p)$ is generic modulo the symmetry group $\mathcal{S}$, i.e., the vertices of a set of representatives for the vertex orbits under the action of $\mathcal{S}$ are placed in `generic' positions (see \citet{BS6, BS1, BSWWorbit} for details), and $(G,p)$ also possesses an $\mathcal{S}$-symmetric infinitesimal flex, then $(G,p)$ also has a finite flex which preserves all the symmetries in $\mathcal{S}$ throughout the path.
\end{theorem}

\subsection{Orbit rigidity matrices for symmetric frameworks}
\label{subsec:SymmetricOrbitMatrices}

To determine whether a given framework $(G,p)\in \mathscr{R}_{(G,\mathcal{S})}$ possesses an $\mathcal{S}$-symmetric infinitesimal flex,
we can use the techniques from group representation theory described in \citet{FGsymmax, KG2, BS6, BS2}.
In the recent paper \citet{BSWWorbit} the `orbit matrix' was introduced as a simplifying alternative to detect symmetric infinitesimal motions in symmetric frameworks and to predict finite flexes for configurations which are generic within the symmetry.  In fact, it is shown in \citet{BSWWorbit} that the orbit matrix is equivalent to the submatrix block $\tilde{\mathbf{R}}_1(G,p)$ studied in \citet{BS6}, but the construction is transparent, and the entries in the matrix can be explicitly derived, without using techniques from group representation theory.

For a $d$-dimensional framework $(G,p)\in \mathscr{R}_{(G,\mathcal{S})}$ which has no joint that is `fixed' by a non-trivial symmetry operation in $\mathcal{S}$ (i.e., $(G,p)$ has no joint $p_i$ with $s(p_i)=p_i$ for some $s\in \mathcal{S}$, $s\neq id$),
the construction of the orbit matrix becomes particularly easy (see Definition \ref{orbitmatrixdef}), because, in this case, the orbit matrix has a set of $d$ columns for each orbit of vertices under the group action of $\mathcal{S}$.

\begin{defin}\label{orbitmatrixdef} \citep{BSWWorbit} Let $\mathcal{S}$ be a symmetry group in dimension $d$ and let $(G,p)\in \mathscr{R}_{(G,\mathcal{S})}$ be a framework which has no joint that is `fixed' by a non-trivial symmetry operation in $\mathcal{S}$. Further, let $\mathscr{O}_{V(G)}=\{1,\ldots, v_0\}$ be a set of representatives for the orbits  $\mathcal{S}(i)=\{s(i)|\,s\in \mathcal{S}\}$ of vertices of $G$. For each edge orbit $\mathcal{S}(e)=\{s(e)|\,s\in \mathcal{S}\}$ of $G$, the \emph{orbit matrix} $\mathbf{O}(G,p,\mathcal{S})$ of $(G,p)$ has the following corresponding ($dv_0$-dimensional) row vector:
\begin{description}
\item[Case 1:] If the two end-vertices of the edge $e$ lie in distinct vertex orbits, then there exists an edge in $\mathcal{S}(e)$ that is of the form $\{a,s(b)\}$ for some $s\in \mathcal{S}$, where $a,b\in\mathscr{O}_{V(G)}$.
The row we write in $\mathbf{O}(G,p,\mathcal{S})$ is:
        \begin{displaymath}\renewcommand{\arraystretch}{0.8}
        \bordermatrix{  & &  a &  &  b &  \cr & 0 \ldots 0 & \big(p_a-s(p_b)\big) & 0  \ldots  0 & \big(p_b-s^{-1}(p_a)\big) & 0  \ldots  0}\textrm{.}
    \end{displaymath}
\item[Case 2:] If the two end-vertices of the edge $e$ lie in the same vertex orbit, then there exists an edge in $\mathcal{S}(e)$ that is of the form $\{a,s(a)\}$ for some $s\in \mathcal{S}$, where $a\in\mathscr{O}_{V(G)}$. The row we write in $\mathbf{O}(G,p,\mathcal{S})$ is:
\begin{displaymath}
\bordermatrix{ & &  a & \cr & 0  \ldots  0 & \big(2p_a-s(p_a)-s^{-1}(p_a)\big) & 0  \ldots  0}\textrm{.}
    \end{displaymath}
\end{description}
\end{defin}

\medskip

\begin{figure}[htp]
\begin{center}
\begin{tikzpicture}[very thick,scale=1]
\tikzstyle{every node}=[circle, draw=black, fill=white, inner sep=0pt, minimum width=5pt];
        \path (0,0) node (p1) [label = left: $p_{1}$] {} ;
       \path (0,-1.1) node (p2) [label = left: $p_{2}$] {} ;
   \path (2.5,-1.8) node (p3) [label = right: $p_{3}$] {} ;
   \path (2.5,-0.7) node (p4) [label = right: $p_{4}$] {} ;
   \draw (p1) -- (p4);
      \draw (p3) -- (p4);
     \draw (p2) -- (p3);
      \draw (p2) -- (p1);
\filldraw[fill=black, draw=black]
    (1.25,-0.9) circle (0.01cm);
\node [rectangle, draw=white, fill=white] (a) at (1.25,-0.7) {\small center};
\node [rectangle, draw=white, fill=white] (b) at (1.25,-2.4) {(a)};
         \end{tikzpicture}
         \hspace{2cm}
        \begin{tikzpicture}[auto, node distance=2cm, very thick]
\tikzstyle{vertex1}=[circle, draw, fill=white, inner sep=1pt, minimum
width=3pt];
\tikzstyle{vertex2}=[circle, draw, fill=yellow, inner sep=1pt, minimum
width=3pt];
\tikzstyle{vertex3}=[circle, draw, fill=blue, inner sep=1pt, minimum
width=3pt];
\tikzstyle{vertex4}=[circle, draw, fill=green, inner sep=1pt, minimum
width=3pt];
\tikzstyle{gain} = [fill=white, inner sep =1pt,  anchor=center];
\tikzstyle{vertexbb}=[rectangle, draw=white, fill=white, inner sep=0pt, minimum width=5pt];
\node[vertex1] (2) at (1, -1.9) {$1$};
\node[vertex1] (3) at (2, -0.5) {$2$};
\path
(2) edge [bend left] (3);
\pgfsetarrowsend{stealth}[ shorten >=2pt]
\path
(2) edge [bend right] node[gain] {\small $C_2$} (3);
\pgfsetarrowsend{}
\node [vertexbb] (b) at (1.5,-2.7) {(b)};
\end{tikzpicture}
\end{center}
\vspace{-0.3cm}
\caption{The framework $(G,p)\in\mathscr{R}_{(G,\mathcal{C}_2)}$ (a) and its corresponding symmetric orbit graph (b).}
\label{quadc2pic}
\end{figure}
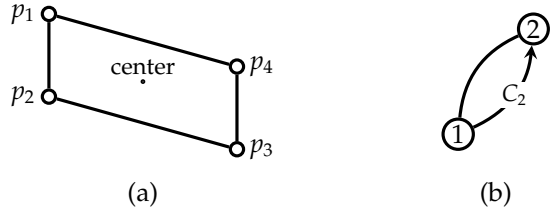

\noindent {\bf Example 3.2.1}
To illustrate the above definition, we consider the $2$-dimensional framework $(G,p)$ with point group $\mathcal{C}_2=\{id, C_2\}$ depicted in Figure \ref{quadc2pic} as an example. If we denote $p_1=(a,b)$, $p_2=(c,d)$, $p_3=(-a,-b)$, and $p_4=(-c,-d)$, then the rigidity matrix of $(G,p)$ is
\begin{displaymath}\bordermatrix{
                &1&2& {3}=C_2(1) &{4}=C_2(2) \cr
                \{1,2\}&(a-c,b-d) &  (c-a,d-b)  &0\ 0 & 0 \ 0\cr
                 \{1,C_2(2)\}& (a+c,b+d) & 0 \ 0  & 0 \ 0 & (-a-c,-b-d)\cr
             C_2(\{1,2\})& 0 \ 0  & 0 \ 0 & (c-a,d-b) &  (a-c,b-d) \cr
             C_2(\{1,C_2(2)\})& 0 \ 0  &  (a+c,b+d)       & (-a-c,-b-d) & 0 \ 0}
\end{displaymath}
The orbit matrix $\mathbf{O}(G,p,\mathcal{C}_2)$ of $(G,p)$ will only have two rows, one for each representative of the edge orbits under the  action of $\mathcal{C}_2$. (Note that if we are only interested in infinitesimal motions and self-stresses of $(G,p)$ that are $\mathcal{C}_2$-symmetric, then it indeed suffices to focus on the first two rows of the rigidity matrix of $(G,p)$. The other two rows are clearly redundant in this symmetric context!). Further, $\mathbf{O}(G,p,\mathcal{C}_2)$ will have only four columns, because $G$ has only two vertex orbits under the action of $\mathcal{C}_2$, represented by the vertices $1$ and $2$, for example, and each of the joints $p_1$ and $p_2$ has two degrees of freedom in the plane. Since both edge orbits satisfy Case 2 in Definition \ref{orbitmatrixdef}, $\mathbf{O}(G,p,\mathcal{C}_2)$ has the following form:
\begin{displaymath}\bordermatrix{
                &1&2 \cr
                \{1,2\}&(p_1-p_2) &  (p_2-p_1)\cr
                 \{1,C_2(2)\}& \big(p_1-C_2(p_2)\big) & \big(p_2-C_2^{-1}(p_1)\big) \cr
             }=\bordermatrix{
                &1&2 \cr
                &(a-c,b-d) &  (c-a,d-b)\cr
                 & (a+c,b+d) & (c+a,d+b) \cr
             }
\end{displaymath}
\eop

We can use the `symmetric orbit graph' to describe the underlying combinatorial structure for the orbit matrix of a symmetric framework (see also Figure \ref{quadc2pic} (b)):

\begin{defin} \label{defsymorbgr} The \emph{symmetric orbit graph} $G_{\mathcal{S}}$ of a framework $(G,p)\in \mathscr{R}_{(G,\mathcal{S})}$ is a labeled multigraph (it may contain loops and multiple edges) whose vertex set $\{1,\ldots, v_0\}$ is a set of representatives of the vertex orbits of $G$ under the action of $\mathcal{S}$, and whose edge set is defined as follows. For each edge orbit of $G$ under the action of $\mathcal{S}$, there exists one edge in $G_{\mathcal{S}}$: for an edge orbit satisfying Case 1 of Definition \ref{orbitmatrixdef}, $G_{\mathcal{S}}$ has a directed edge connecting the vertices $a$ and $b$. If the edge is directed from $a$ to $b$, it is labeled with $s$, and if the edge is directed from $b$ to $a$, it is labeled with $s^{-1}$. For simplicity we omit the label and the direction of the edge if $s=id$. Similarly, for an edge orbit satisfying Case 2 of Definition \ref{orbitmatrixdef}, $G_{\mathcal{S}}$ has a loop at the vertex $a$ which is labeled with $s$.
\end{defin}


The key result for the orbit matrix is the following:

\begin{theorem}\label{orbitmatrixthm} \citep{BSWWorbit}
Let $S$ be a symmetry group and let $(G,p)$ be a framework in $\mathscr{R}_{(G,\mathcal{S})}$. Then the solutions to $\mathbf{O}(G,p,\mathcal{S})u = 0$ are isomorphic to the space of $\mathcal{S}$-symmetric infinitesimal motions of the original framework $(G,p)$.
\end{theorem}

As an immediate consequence of Theorems  \ref{thm:flexes} and \ref{orbitmatrixthm}, we have the following Maxwell type counting rule for detecting finite `symmetry-preserving' flexes in symmetric frameworks:

\begin{theorem}\label{symMaxthm} \citep{BS4, BS6, BSWWorbit}
Let $\mathcal{S}$ be a symmetry group in dimension $d$ and let $(G,p)$ be a framework in $\mathscr{R}_{(G,\mathcal{S})}$ which has no joint that is `fixed' by a non-trivial symmetry operation in $S$. Further, let $e_0$ and $v_0$ denote the number of edge orbits and vertex orbits under the action of $\mathcal{S}$, respectively, and  let $triv_{\mathcal{S}}$ denote the dimension of the space of $\mathcal{S}$-symmetric infinitesimal rigid motions of $(G,p)$. If \begin{equation}\label{eqflex} e_0<dv_0 - triv_{\mathcal{S}} \textrm{,}\end{equation} then $(G,p)$ has an $\mathcal{S}$-symmetric infinitesimal flex. If the joints of $(G,p)$ also span all of $\mathbb{R}^d$ (in an affine sense) and are in generic position modulo $\mathcal{S}$, then there even exists a finite flex of $(G,p)$ which preserves the symmetries in $\mathcal{S}$ throughout the path.
\end{theorem}

The dimension $triv_{\mathcal{S}}$ of the space of $\mathcal{S}$-symmetric infinitesimal rigid motions of $(G,p)\in \mathscr{R}_{(G,\mathcal{S})}$ can easily be computed using the techniques described in \citet{BS4, BS2}. In particular, in dimension 2 and 3,  $triv_{\mathcal{S}}$ can be deduced immediately from the character tables given in \citet{cfgsw}.
Thus, in order to check condition (\ref{eqflex}), it is only left to determine the size of the orbit matrix $\mathbf{O}(G,p,\mathcal{S})$, which in turn requires only a simple count of the vertex orbits and edge orbits of the graph $G$ under the action of $\mathcal{S}$ (see also Example 3.2.1 and Table~\ref{table:pointgroups}).

Note that for a symmetry group $\mathcal{S}$ in \emph{any} dimension $d$, the dimension $t_{\mathcal{S}}$ of the space of $S$-symmetric infinitesimal \emph{translations} can also be obtained in a very intuitive way, without using the techniques in \citet{BS4, BS2}:
if for a symmetry operation $s\in \mathcal{S}$, we let $F_s$ denote the \emph{symmetry element} corresponding to $s$ (i.e., $F_s=\{a\in\mathbb{R}^{d} \,| s(a)=a\}$), then $t_{\mathcal{S}}$ is simply the dimension of the \emph{symmetry element of the group $\mathcal{S}$}, i.e., the dimension of the linear subspace $\bigcap_{s\in \mathcal{S}}F_s$ of $\mathbb{R}^d$, because the initial velocity vectors of an $\mathcal{S}$-symmetric infinitesimal translation must all be contained in the space $\bigcap_{s\in \mathcal{S}}F_s$.

For example, for the `reflectional' symmetry group $\mathcal{C}_s=\{id,\sigma\}$ in dimension $d$, it is easy to see that the space of $\mathcal{C}_s$-symmetric infinitesimal translations is of dimension $(d-1)$, since it consists of those translations whose velocity vectors are elements of the $(d-1)$-dimensional mirror-plane $F_{\sigma}=\mathbb{R}^d\cap F_{\sigma}=F_{id}\cap F_{\sigma}=\bigcap_{s\in \mathcal{C}_s}F_s$ corresponding to $\sigma$ (see also Figure \ref{fulsym} (b)).

However, finding the dimension of the space of  $\mathcal{S}$-symmetric infinitesimal \emph{rotations} heuristically, without the techniques in \citet{BS4, BS2}, becomes increasingly hard, if not impossible, in dimensions $>3$.

\medskip

\begin{figure}
    \begin{center}
        \subfigure[] { \begin{tikzpicture}[very thick,scale=1]
\tikzstyle{every node}=[circle, draw=black, fill=white, inner sep=0pt, minimum width=5pt];
    \draw [thick, blue] (0.1,0.2) -- (2.1,-1.1);
      \draw [thick, blue](2.3,-2.3) -- (2.1,-1.1);
     \draw [thick, blue](0.3,-1.1) -- (2.3,-2.3);
      \draw [thick, blue](0.3,-1.1) -- (0.1,0.2);
            \draw [thick, ->, red](0.1,0.18) -- (0.52,-0.26);
      \draw [thick, ->, red](2.3,-2.3) -- (1.88,-1.84);
      \draw [thick, ->, red](0.3,-1.1) -- (-0.5,-1.65);
      \draw [thick, ->, red](2.1,-1.1) -- (2.9,-0.55);
      \filldraw[fill=black, draw=black]
    (1.2,-1.05) circle (0.004cm);
              \end{tikzpicture}}
         \hspace{1.7cm}
        \subfigure[] { \includegraphics [width=.30\textwidth]{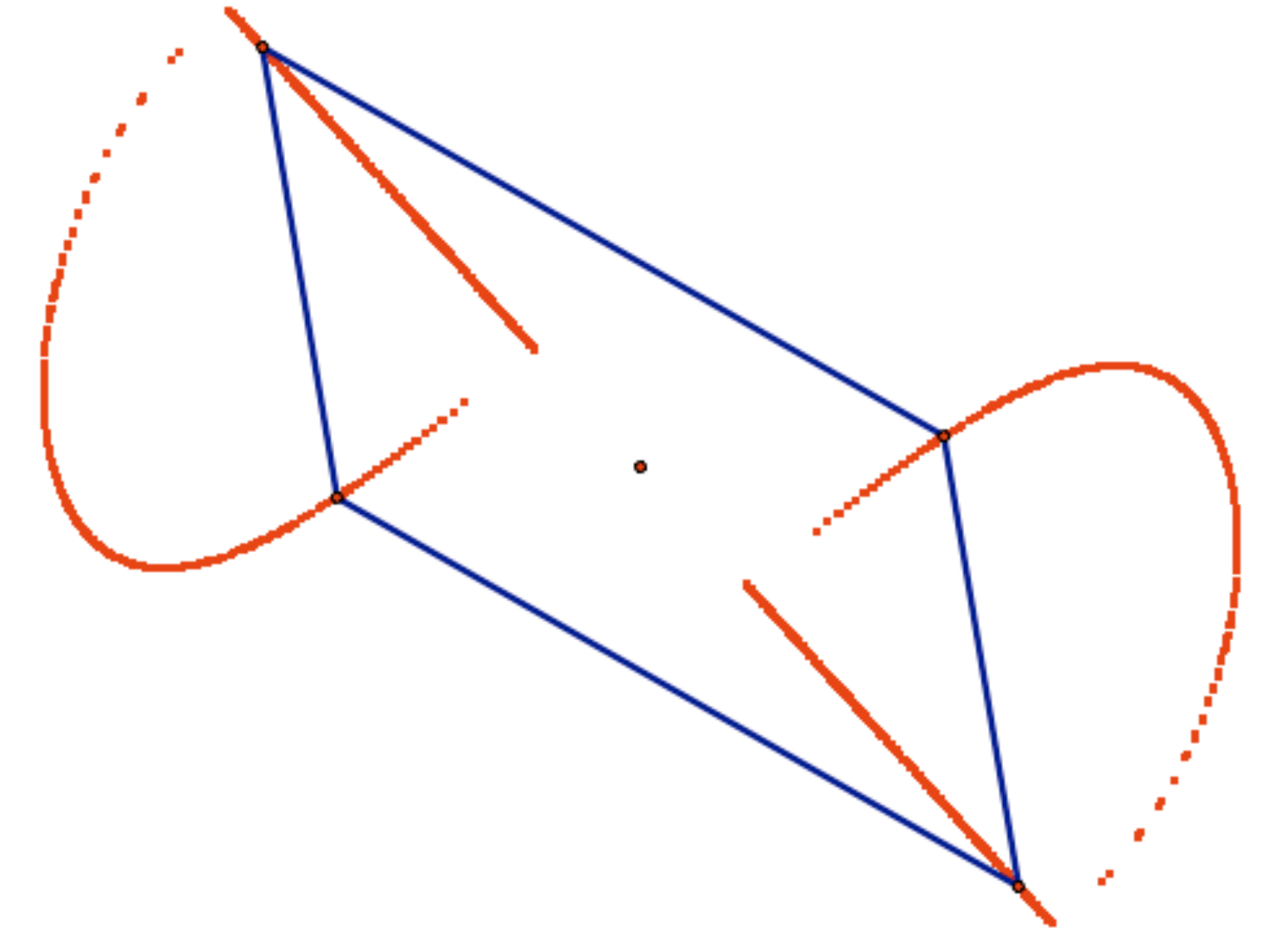}}
 \end{center}
  \caption{A $\mathcal{C}_2$-symmetric infinitesimal flex of the framework from Example 3.2.1 (a) and the path taken by the joints of the framework under the corresponding symmetry-preserving finite flex (b). }
    \label{fig:mechanismc2}
\end{figure}
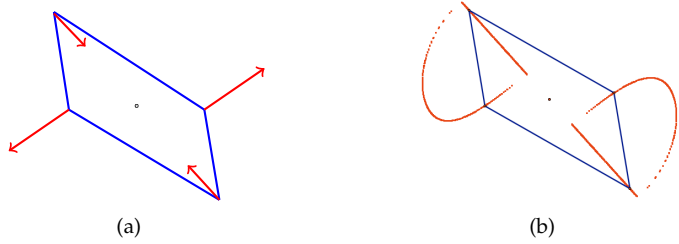

\noindent {\bf Example 3.2.2}
Let's apply Theorem \ref{symMaxthm} to the framework $(G,p)$ we considered in Example 3.2.1 (see also Figures \ref{quadc2pic} (a) and \ref{fig:mechanismc2}). We clearly have $dv_0=2\cdot 2=4$ and $e_0=2$. Further, we have $triv_{\mathcal{C}_2}=1$, since the only infinitesimal rigid motions that are $\mathcal{C}_2$-symmetric are the ones that correspond to rotations about the origin (see \citet{BS4, BS2} for details). Thus, we have
\begin{displaymath}
e_0=2<3=dv_0 - triv_{\mathcal{C}_2}\textrm{.}
\end{displaymath}
So, by Theorem \ref{symMaxthm}, we may conclude that any realization of $G$ which is `generic' modulo the half-turn symmetry has a symmetry-preserving finite flex (Figure~\ref{fig:mechanismc2}).
\eop

Note that the standard (non-symmetric) Maxwell count also detects a finite flex for the framework in Example 3.2, since for the graph $G$, we have $e=4<5=2v-3$. However, as shown in \citet{BSWWorbit}, there exist a range of interesting and famous examples in $3$-space, including the Bricard octahedra \citep{bricard, Stachel} and the flexible cross-polytopes \citep{Stachel4D}, which can be shown to be flexible via the symmetric count in Theorem \ref{symMaxthm}, but not with the standard non-symmetric Maxwell-Laman type counts for rigidity.

The following table shows the symmetric Maxwell type counts for a selection of point groups $\mathcal{S}$ in $3$-space. For simplicity at this stage, we assume that no joint and no bar is fixed by a non-trivial element in $\mathcal{S}$, so that all vertex orbits and edge orbits under the action of $\mathcal{S}$ have the same size $k_{\mathcal{S}}$. (Recall that a joint $p_i$ is fixed by $s\in \cal S$ if $s(p_i)=p_i$; a bar $\{p_i,p_j\}$ is fixed by $s\in \cal S$ if either $s(p_i)=p_i$ and $s(p_j)=p_j$ or $s(p_i)=p_j$ and $s(p_j)=p_i$). So, in particular, both the number of joints, $v$, and the number of bars, $e$ are divisible by $k_{\mathcal{S}}$. A necessary condition for rigidity in $3$-space is $e\geq 3v-6$ (recall Theorem \ref{maxtheorem}). For each group $\mathcal{S}$ in Table \ref{table:pointgroups}, $e$ is chosen to be the smallest number which satisfies $e\geq 3v-6$ and is divisible by $k_{\mathcal{S}}$; that is, $e$ is chosen to be the least number of edges for the framework to be rigid without symmetry and to be compatible with the symmetry constraints given by $\mathcal{S}$.
The integer $f_{\mathcal{S}}$ in the final column indicates an $f_{\mathcal{S}}$-dimensional space of $\mathcal{S}$-symmetric infinitesimal flexes if $f_{\mathcal{S}}>0$, and a $(-f_{\mathcal{S}})$-dimensional space of $\mathcal{S}$-symmetric self-stresses if $f_{\mathcal{S}}<0$. (See Section \ref{subsec:Stresses} for more on self-stresses.)

\begin{table}[h!]
 \centering
 \caption{Impact of some $3$-space point groups on counts for rigidity. }
\begin{tabularx}{1.\textwidth}{X}
\[ \begin{array}{|c|c|c|c|c|c|c|} \hline
  \text{\grey{$\mathcal{S}$}}  &   \grey{k_{\cal S}}  & \grey{triv_{{\cal S}}}& \grey{e}&  \text{\grey {$e_0$ }}&  \text{\grey {$3v_0 -triv_{{\cal S}}$}}& \text{\grey{$f_{\mathcal{S}}$}} \\
 \hline
\mathcal{C}_{1} &  1 & 6 &3v-6&3v_0-6 & 3v_0-6&0\\
\hline
\mathcal{C}_{i} & 2 & 3&3v-6& 3v_{0}-3&3v_{0}-3&0\\
\hline
\mathcal{C}_{2} & 2 &2 &3v-6& 3v_{0}-3&3v_{0}-2&1\\
\hline
 \mathcal{C}_{s} & 2 & 3&3v-6& 3v_{0}-3& 3v_{0}-3& 0 \\
\hline
\mathcal{C}_{2h} & 4 & 1 &3v-4& 3v_{0}-1& 3v_{0}-1&0\\
\hline
\mathcal{D}_{2h} & 8 & 0 &3v& 3v_{0}&3v_{0}-0&0\\
\hline
\mathcal{C}_{4}  & 4 & 2&3v-4&3v_{0}-1& 3v_{0}-2&(-1)\\
\hline
\mathcal{C}_{3} & 3 & 2&3v-6&3v_{0}-2 &3v_{0}-2&0\\
\hline
\end{array}  \]
\end{tabularx}
\label{table:pointgroups}
\end{table}%
The final column of Table \ref{table:pointgroups} indicates that at `generic' configurations, the frameworks with $\mathcal C_2$ symmetry always have a finite flex, while those with $\mathcal C_4$ symmetry are always stressed.

\section{Periodic frameworks }
\label{sec:PeriodicFrameworks}
\subsection{From infinite periodic frameworks to a finite orbit graph}
Let $\wt G$ be a simple infinite graph with finite degree at each vertex. Let $\wt \p$ be a placement of the vertices $V(\wt G)$ in $\mathbb R^3$, such that the resulting framework $(\wt G, \wt \p)$ is invariant with respect to three linearly independent translations $\bt_1, \bt_2, \bt_3 \in \mathbb R^3$. We assume without loss of generality that $\bt_1$ lies on the $x$-axis, and $\bt_2$ lies in the $xy$-plane. Let $L$ be the matrix whose rows are these translations:
$$L = \left[\begin{array}{c}\bt_1 \\\bt_2 \\\bt_3\end{array}\right] = \left[\begin{array}{ccc}t_{11} & 0 & 0 \\t_{21} & t_{22} & 0 \\t_{31} & t_{32} & t_{33}\end{array}\right].$$
We call the pair $\pfw$ a {\it periodic framework}. This definition can be adapted to describe periodic objects in $d$ dimensions, but here our focus will be on crystal-like structures in two and three dimensions. In two dimensions,
$$L = \left[\begin{array}{c}\bt_1 \\\bt_2 \end{array}\right] = \left[\begin{array}{cc}t_{11} & 0  \\t_{21} & t_{22}  \end{array}\right],$$
and the other definitions are similarly adapted.

\begin{figure}[h!]
\begin{center}
\includegraphics[width=3in]{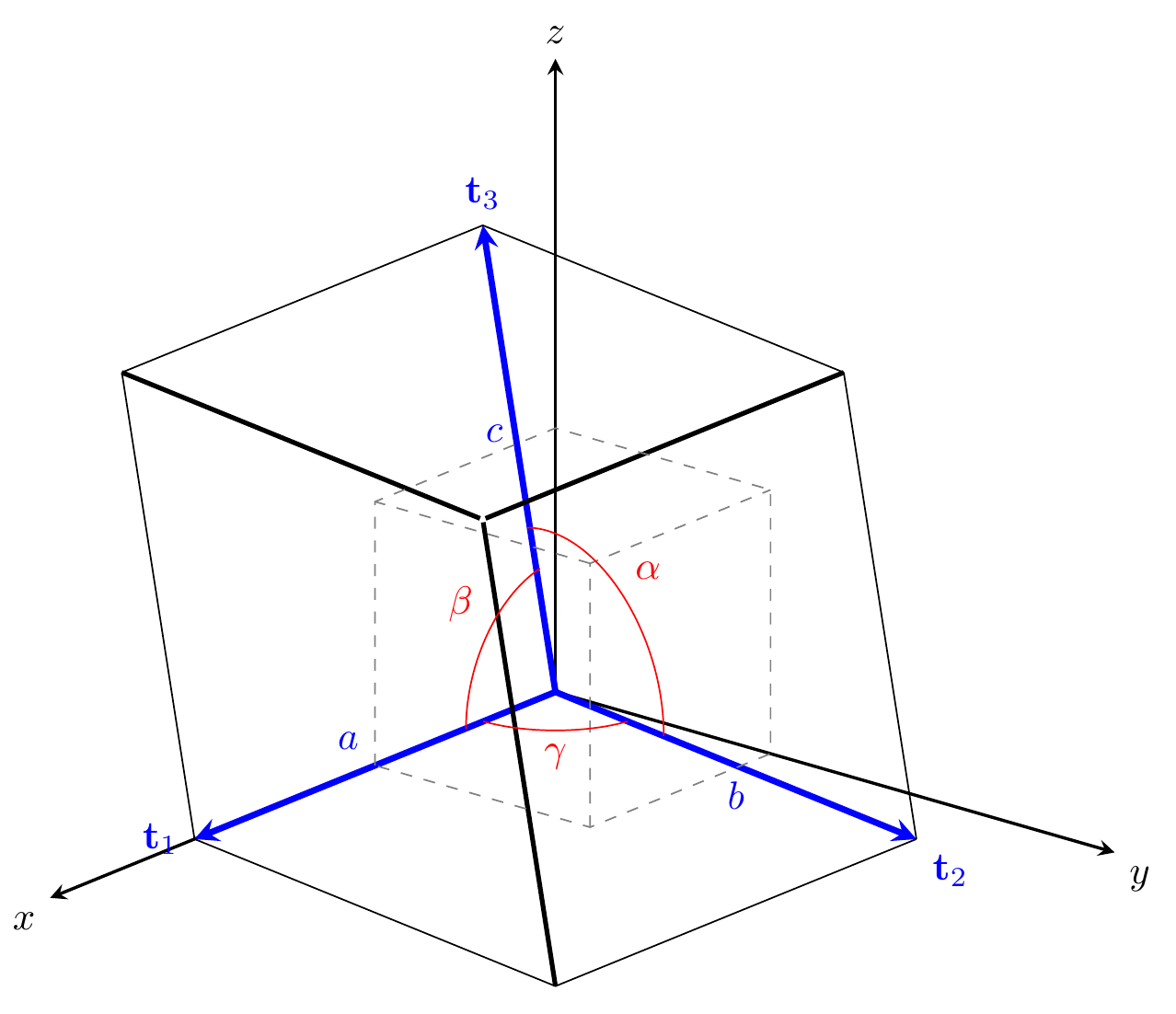}
\caption{\label{fig:unitCell}The unit cell defined by the translations $\bt_1, \bt_2, \bt_3$. The unit cube (all side lengths $1$) is shown in the dashed lines.}
\end{center}
\end{figure}

The three translations $\bt_1, \bt_2, \bt_3$ define a parallelepiped called the {\it unit cell}. This can be equivalently described by three lengths $a, b, c$, and three angles $\alpha, \beta, \gamma$, as illustrated in Figure \ref{fig:unitCell}. These coordinates are standard in crystallography texts \citep{Wiki}.  We may obtain one representation from the other by a change of coordinates.

A copy of the unit cell centered at the origin is given by
$$U = \left\{a_1\bt_1 + a_2\bt_2 + a_3 \bt_3 \ | \ -\frac{1}{2} \leq a_i < \frac{1}{2}\right\},$$
with its {\it boundary} defined to be the boundary of the closed parallelepiped:
$$\overline U = \left\{a_1\bt_1 + a_2\bt_2 + a_3 \bt_3  \ | \   |a_i| \leq \frac{1}{2}\right\}.$$
We assume without loss of generality that no vertex of $V(\wt G)$ lies on the boundary of the unit cell (we may simply translate the framework $\pfw$ until no vertex lies on the boundary).

The translations $\bt_1, \bt_2, \bt_3$ generate a crystalline lattice $\bt_1\mathbb Z \times \bt_2 \mathbb Z \times \bt_3 \mathbb Z$. This lattice partitions $\mathbb R^3$ into copies of the unit cell, each containing exactly one lattice point. Let each cell be centered at that lattice point, and let the cells be indexed by $\mathbb Z^3$ according to the lattice point they contain.

\begin{figure}[h!]
\begin{center}
\subfigure[$\pfw$]{\label{fig:periodicFramework}\includegraphics[width=1.5in]{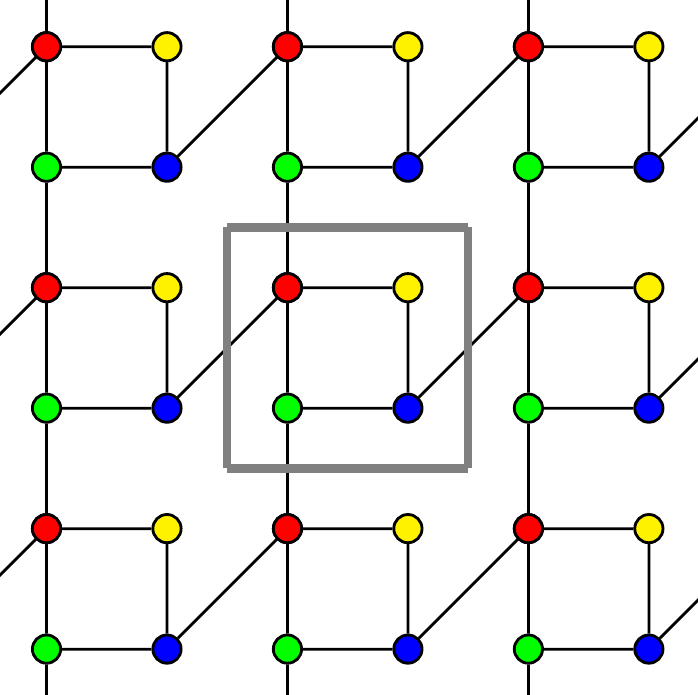}}\hspace{.5in}
\subfigure[$\pog$]{\label{fig:finiteGraph}\begin{tikzpicture}[auto, node distance=2cm]
\tikzstyle{vertex1}=[circle, draw, fill=red, inner sep=1pt, minimum width=3pt, font=\footnotesize];
\tikzstyle{vertex2}=[circle, draw, fill=yellow, inner sep=1pt, minimum width=3pt, font=\footnotesize];
\tikzstyle{vertex3}=[circle, draw, fill=blue, inner sep=1pt, minimum width=3pt, font=\footnotesize];
\tikzstyle{vertex4}=[circle, draw, fill=green, inner sep=1pt, minimum width=3pt, font=\footnotesize];
\tikzstyle{gain} = [fill=white, inner sep = 0pt,  font=\footnotesize, anchor=center];

\node[vertex1] (1) {$1$};
\node[vertex2] (2) [below right of=1] {$2$};
\node[vertex3] (3) [below left of=2] {$3$};
\node[vertex4] (4) [below left of=1] {$4$};

\path[thick] (1) edge (2)
(2) edge (3)
(3) edge   (1)
(1) edge [bend right]  (4)
(4) edge [bend right]   (1)
(4) edge  (3);

\pgfsetarrowsend{stealth}[ shorten >=2pt]
\path
(3) edge node[gain] {$(1,0)$} (1)
(1) edge [bend right] node[gain] {$(0,1)$} (4);

\pgfsetarrowsend{}

\path (1) edge  (2)
(2) edge  (3)
(4) edge [bend right]   (1)
(4) edge  (3);

\end{tikzpicture}}
\caption{A 2-dimensional periodic framework (a) with the $(0,0)$-cell indicated. The periodic orbit graph is shown in (b).  \label{fig:gainGraph}  }
\end{center}
\end{figure}
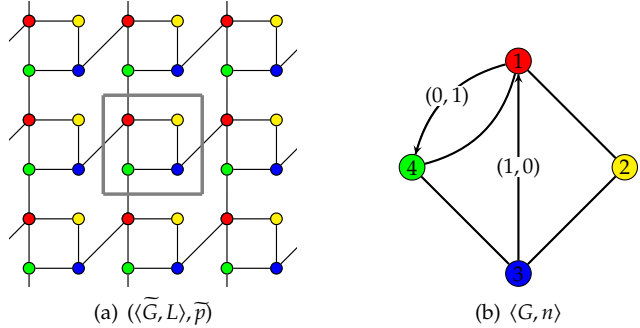

We now use this partition to define a finite labeled graph $G$ which represents the periodic framework, and can be used to study its rigidity. The edges of this graph are labeled invertibly by elements of the group $\mathbb Z^3$ in a way that captures the periodic structure of $\pfw$. Figure \ref{fig:gainGraph} illustrates this process for the analogous 2-dimensional case, in which edges are labeled by elements of the group $\mathbb Z^2$.

\begin{defin}\label{def:periodicOrbitGraph}
The \emph{periodic orbit graph} $\pog$ of a periodic framework $\pfw$ is a labeled multigraph whose vertex set $V(G)=\{1,\ldots, v\}$ consists of the vertices appearing in the $(0,0,0)$-cell of $\pfw$.  The edge set $E(G)$ and the labelling  $n: E(G) \rightarrow \mathbb Z^3$ on these edges are defined as follows:
\begin{description}
\item[{\bf Case 1:}] An edge $\wt e \in E(\wt G)$ whose length is completely contained within the $(0,0,0)$-cell must connect two distinct vertices of $V(\wt G)$ which also lie in the $(0,0,0)$-cell. For every such edge $\wt e$, define $e \in E(G)$ to be the edge connecting the corresponding vertices of $V(G)$. Assign this edge an arbitrary direction, and label it with $(0,0,0)$. For visual simplicity in our diagrams, edges labeled by the zero (identity) element appear as unlabeled, undirected edges.

\item[{\bf Case 2:}] Let $\wt e \in E(\wt G)$ be an edge that crosses the boundary of the $(0,0,0)$-cell. In particular, suppose  $\wt e$ connects the vertex $v$ in the $(0,0,0)$-cell with the vertex $w$ in the $(n_1, n_2, n_3)$-cell. Then define $e \in E(G)$ to be a directed edge of $G$ that originates in the vertex $v$, terminates at the vertex $w$. Assign this directed edge the label  $(n_1, n_2, n_3) \in \mathbb Z^3$.
\end{description}
\end{defin}
The periodic orbit graph $\pog$ contains a set of representatives of the vertex and edge orbits of $\pfw$ under the action of $\mathbb Z^3$. Let $\p$ be the restriction of $\wt \p$ to the vertices of the $(0,0,0)$-cell,  $\p: V \rightarrow U.$  We call the pair $\pofw$ the {\it periodic orbit framework} of $\pfw$, and the labeled multigraph $\pog$ will be called the {\it periodic orbit graph.}

Let the edges $E(G)$ be ordered. An edge $e_k$ of $\pog $ is denoted by $\{i,j; n_k\}$, where $i, j \in V(G)$ and $n_k \in \mathbb Z^3$. This edge corresponds to an equivalence class of bars in the periodic framework $\pfw$, which contains the bar $(\p_i,\p_j+n_k T)$, and all of its translates by $\sum a_i \bt_i$, for integers $a_i$. This edge can be equivalently represented by $\{j, i; -n_k\}.$

It has been shown \citep{Theran, Ross} that every ($d$-dimensional) periodic graph admits such a representation, and is invariant under the choice of unit cell of a specified size. It is also known that every finite directed multigraph whose edges are labeled with elements of the group $\mathbb Z^d$ can be realized as a $d$-dimensional periodic graph \citep{ Ross}. In general, directed multigraphs whose edges are labeled by elements of a group, with the reverse direction implicitly labeled by the inverse group element,  are known as {\it gain graphs} \citep{GrossTucker}.

\subsection{Periodic rigidity and infinitesimal rigidity}
\label{subsection:periodicrigidity}
We may define notions of rigidity and infinitesimal rigidity for the periodic framework $\pfw$ simply by a direct application of the ideas from Section \ref{sec:Preliminaries}. However, in this paper we are interested exclusively in motions and infinitesimal motions of $\pfw$ which preserve the periodicity of the framework. We say that $\pfw$ is {\it periodic (infinitesimally) rigid} if the only periodic (infinitesimal) motions of $\pfw$ are trivial.  By defining rigidity and infinitesimal rigidity for a periodic orbit framework $\pofw$, we are able to identify precisely these characteristics of the corresponding periodic framework, $\pfw$. Further details can be found in \citet{BorceaStreinuI, Theran, Ross}.

Let $\pofw$ be a periodic orbit framework, and let $e_k = \{i, j; n_k\}$ be an edge of $\pog$, with $n_k = (n_{k1}, n_{k2}, n_{k3})$. The {\it edge length} of $e_k$ is given by the Euclidean length of the vector $\p_i - (\p_j + n_kL)$.
\begin{eqnarray}
\|e_k\|^2 &  =  & \|\p_i - (\p_j + n_kL) \|^2 = \sum_{\ell=1}^3(\p_{i\ell} - (\p_{j\ell} + (n_{k}L)_{\ell}))^2 \nonumber \\
& = & (\p_{i1} - (\p_{j1} + n_{k1}t_{11} + n_{k2}t_{21} + n_{k3}t_{31}))^2 +  \nonumber \\
& & (\p_{i2} - (\p_{j2} + n_{k2}t_{22} +  n_{k3}t_{32}))^2 + (\p_{i3} - (\p_{j3} + n_{k3}t_{33}))^2
\label{eqn:edgeLength}
\end{eqnarray}
Because there are a finite number of edges of $\pog$, we may use the above periodic edge length to define \emph{rigidity} and \emph{flexibility} of the orbit framework $\pofw$. It is analogous to the definition presented in Section \ref{sec:Preliminaries} and we omit it here.

Letting the positions of the vertices $\p_i = (\p_{i1}, \p_{i2}, \p_{i3})$ and the generators of the lattice $t_{11}, \dots, t_{33}$ vary with time, infinitesimal motions of $\pofw$ can be found by differentiating (\ref{eqn:edgeLength}), with the assumption that $\|e_k\| = C$, a constant. We obtain
$$(\p_i -  (\p_j+n_kL)) \cdot (d\p_i/{dt} - d\p_j/{dt}) + \mathcal L(i,j; n_k)\cdot (dt_{11}/{dt}, \dots, dt_{33}/{dt}) = 0,$$
where $\mathcal L(i,j; n_k)$ is the $6$-tuple of coefficients of $(dt_{11}/{dt}, \dots, dt_{33}/{dt})$. For example, the coefficient corresponding to $t_{32}$ is $n_{k3}(\p_{i2} - (\p_{j2} + n_{k2}t_{22} +  n_{k3}t_{32}))$.

More generally, an {\it  infinitesimal motion} $(\bu, w)$ of a periodic orbit framework $\pofw$ in $\mathbb R^3$  is a pair of functions
$$\bu : V (G) \rightarrow \mathbb R^3,\  \textrm{and} \  w: L \rightarrow \mathbb R^6$$
such that
\begin{equation}
(\p_i -  (\p_j+n_kL)) \cdot (\bu_i - \bu_j) + \mathcal L(i,j; n_k)\cdot w = 0 \ \textrm{for all} \ \{i, j; n_k \} \in E(\pog).
\label{eqn:infPerMotion}
\end{equation}
An infinitesimal motion $(\bu, w)$ of $\pofw$ is called  {\it trivial}  if $w$ is the zero map, and there exists $\bu_0 \in \mathbb R^3$ such that $\bu(i) = \bu_0$ for all $i \in V(G)$.  This simple form follows from the fact that the only isometries of the whole space that preserve the periodic structure are translations.  If an infinitesimal motion is not trivial, then it is called an {\it infinitesimal flex}.   The periodic orbit framework $\pofw$ is {\it infinitesimally rigid} if every infinitesimal motion of $\pofw$ is a trivial one. Otherwise the framework is {\it infinitesimally flexible}.

Let $(\bu, w)$ be an infinitesimal motion of $\pofw$, with $\bu \neq 0$.  When $w:L \rightarrow \mathbb R^6$ is the zero-map, it indicates that the lattice vectors are fixed by the infinitesimal motion. We call such a motion $(\bu, 0)$ a {\it lattice-fixing infinitesimal motion}. If $w \neq 0$, then the infinitesimal motion $(\bu, w)$ is called {\it lattice flexing}. In this paper we regard the lattice-fixing motions as a specialization of the lattice-flexing motions, and hence we will assume all motions are lattice-flexing, unless otherwise noted.

The fixed lattice variation is interesting in its own right, and in two-dimensions admits a concise combinatorial characterization \citep{Ross}. It may also be of interest to consider partial flexing of the lattice. For example, we may ask that the three translation vectors are only allowed to scale, but the angles between them remain fixed at $90^{\circ}$. This would correspond to the translation matrix
$$L =  \left[\begin{array}{ccc}t_{11} & 0 & 0 \\0 & t_{22} & 0 \\0 & 0 & t_{33}\end{array}\right].$$
Such variations will be treated briefly at the end of this section.

\begin{remark} Any infinitesimal motion of a periodic orbit framework $\pofw$ can be extended to a periodic infinitesimal motion of the periodic framework $\pfw$. However, there will be some infinitesimal motions of $\pfw$ that do not preserve the periodic structure,  and therefore do not specialize to infinitesimal motions of $\pofw$. In particular, any infinitesimal motion of $\pfw$ that breaks the periodicity of $\pfw$ will not appear as a motion of $\pofw$. In other words, a periodic framework $\pfw$ may be infinitesimally flexible, yet infinitesimally periodic rigid in our analysis. \end{remark}

It has been demonstrated \citet{BorceaStreinuI, Theran, Ross} that the vector space of periodic infinitesimal motions of a periodic framework $\pfw$ is equivalent to the vector space of infinitesimal motions of the periodic orbit framework $\pofw$. To make this connection, we use a periodic orbit matrix, which we shall now define.

\subsection{Orbit rigidity matrices for periodic frameworks}
\label{subsec:PeriodicOrbitMatrices}
Let $\pofw$ be a three-dimensional periodic orbit framework with  $v=|V(G)|$ and $e = |E(G)|$. The rigidity matrix $\mathbf{R}\pofw$ for the periodic orbit framework is the $e \times (3v + 6)$ matrix with the row corresponding to the edge $\{i,j; n_k\}$ given by
\[  \renewcommand{\arraystretch}{0.8}
     \bordermatrix{  &&  i &  &  j &  & L     \cr
\{i,j; n_k\} &    0 \cdots 0 & \big(\p_i-(\p_j+n_kL)\big) & 0  \cdots  0 & \big(\p_j-(\p_i - n_kL)\big) & 0  \cdots  0 & \mathcal L(i, j; n_k)
     \cr }.\]
The entry $\mathcal L(i, j; n_k)$ is a six-tuple representing the coefficients of $(dt_{11}/{dt}, \dots, dt_{33}/{dt})$. 
Loops may appear in the periodic orbit framework. The row of $\mathbf{R}\pofw$ corresponding to the loop edge $\{i, i; n_{\ell}\}$ is
\[  \renewcommand{\arraystretch}{0.8}
     \bordermatrix{  & &  i &  & L
     \cr
 &    0 \cdots 0 &  0  \cdots  0 & 0  \cdots  0 & \mathcal L(i, i; n_{\ell})
     \cr }.\]
 Note that the matrix $\mathbf{R}\pofw$ is identical to the `augmented compatibility matrix' used by \citet{hutgue}.

\medskip

\noindent {\bf Example 4.3.1}
We consider a two-dimensional periodic orbit framework with two vertices and five edges.
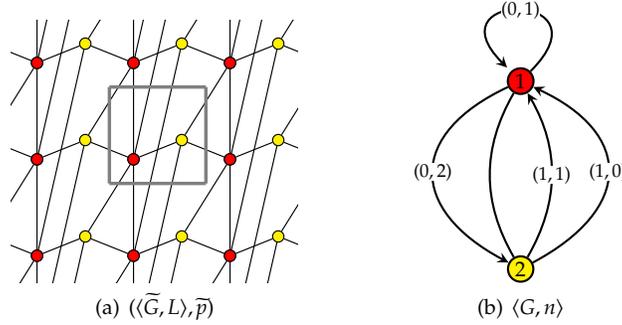
\begin{figure}[h!]
\begin{center}
\subfigure[$\pfw$]{\label{fig:2vertexPeriodic}\includegraphics[width=1.5in]{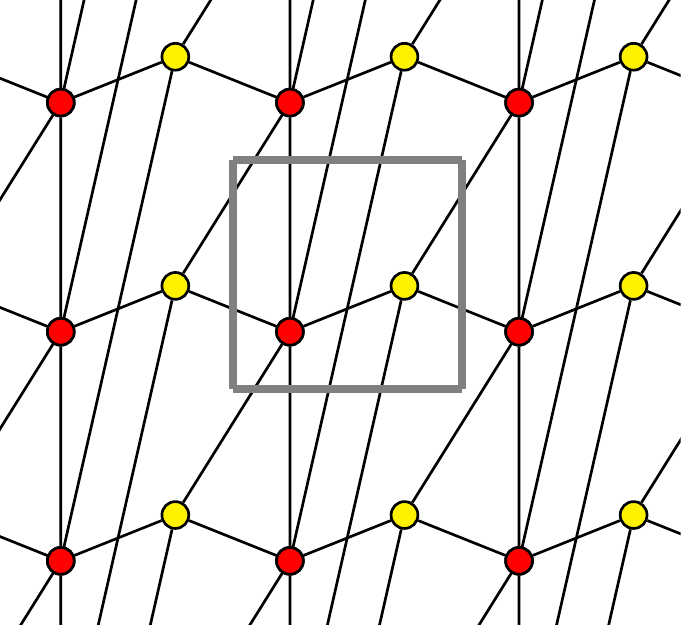}} %
\hspace{.5in}
\subfigure[$\pog$]{\label{fig:2vertexOrbit}\begin{tikzpicture}[auto, node distance=2cm, thick]
\tikzstyle{vertex1}=[circle, draw, fill=red, inner sep=1pt, minimum width=3pt, font=\footnotesize];
\tikzstyle{vertex2}=[circle, draw, fill=yellow, inner sep=1pt, minimum width=3pt, font=\footnotesize];
\tikzstyle{vertex3}=[circle, draw, fill=blue, inner sep=1pt, minimum width=3pt, font=\footnotesize];
\tikzstyle{vertex4}=[circle, draw, fill=green, inner sep=1pt, minimum width=3pt, font=\footnotesize];
\tikzstyle{gain} = [fill=white, inner sep = 0pt,  font=\scriptsize, anchor=center];

\node[vertex1] (1) {$1$};
\node[vertex2] (2) at (0,-2.5) {$2$};

\draw (1) edge [bend right]  (2);

\pgfsetarrowsend{stealth}[ shorten >=2pt]
\path -((1) edge [->,>=stealth,shorten >=2pt, loop, looseness=15] node[gain, midway] {$(0,1)$} (1);
\draw (1) .. controls (-1.5,-.7) and (-1.5, -1.7)  .. node[gain] {$(0,2)$} (2);
\draw (2) edge [bend right] node[gain] {$(1,1)$} (1);
\draw (2) .. controls (1.5, -1.7) and (1.5,-.7)  .. node[gain] {$(1,0)$} (1);
\pgfsetarrowsend{}

\end{tikzpicture}}
\caption{A two vertex example in $\mathbb R^2$.}
\end{center}
\end{figure}
The rigidity matrix $\mathbf{R}\pofw$ for this framework has 5 rows and 7 columns. Below it is broken into two sections: the four columns corresponding to variables $p_1 = (a, b), p_2 = (c, d)$, and the three columns corresponding to the three non-zero variables in $L$: $t_{11}, t_{21}, t_{22}$. Note that the edge $\{1,2; (0,0)\}$ has the entry $(0,0,0)$ in the $ (t_{11}, t_{21}, t_{22})$ columns, but is non-zero elsewhere. In contrast, the loop edge $\{1,1;(1,0)\}$ has zero entries everywhere except in the columns corresponding to $L$.
\[\mathbf{R}\pofw = \bordermatrix{& p_1 = (a, b) & p_2=(c,d) & (t_{11}, t_{21}, t_{22}) \cr
\{1,2; (0,0)\} & \big(p_1-p_2\big) & \big(p_2-p_1) & (0,0,0) \big) \cr
\{1,2; (0,2)\} & \big(p_1-[p_2+(0,2)L]\big) &\big( p_2-[p_1-(0,2)L]\big) &  (\ast, \ast, \ast) \cr
\{2,1; (1,0)\} & \big(p_1-[p_2-(1,0)L]\big) &\big( p_2-[p_1+(1,0)L]\big) & (\ast, \ast, \ast) \cr
\{2,1; (1,1)\} & \big(p_1-[p_2-(1,1)L]\big) & \big(p_2-[p_1+(1,1)L]\big) & (\ast, \ast, \ast) \cr
\{1,1;(0,1)\} & 0 & 0 & (0,0,1)}\]
where the three columns corresponding to $(t_{11}, t_{21}, t_{22})$ (the non-zero entries of $L$) are given by:
\[\bordermatrix{
& t_{11} & t_{21} & t_{22}\cr
& 0 & 0 & 0 \cr
& 0 & 0 & -2(b-(d+2t_{22}))\cr
& -(c - (a +t_{11} )) & 0 & 0\cr
& -(c - (a +t_{11} + t_{21})) & -(c - (a +t_{11} + t_{21})) & -(d-(b + t_{22}))\cr
& 0 & 0 & 1
}\]
The two trivial motions of $\pofw$ are represented by the column vectors $(1,0,1,0,0,0,0)^T$ and $(0,1,0,1,0,0,0)^T$, which are always solutions to the linear system $\mathbf{R}\pofw \cdot (u,w) = 0$.  These are translations of the whole structure.
If we are only interested in the lattice-fixing infinitesimal motions of $\pofw$, we may omit the columns corresponding to $(t_{11}, t_{21}, t_{22})$ and the translations still appear in the modified matrix.
\eop

Returning to three dimensions, we may associate an infinitesimal motion $(\bu, w)$ of the periodic orbit framework $\pofw$ with a column vector $(\bu, w)^{T}$ in $\mathbb R^{3v+6}$. The equations in (\ref{eqn:infPerMotion}) may then be written as the solutions to the linear system described by $\mathbf{R}\pofw \cdot (\bu, w)^T = 0$. There will always be three trivial solutions corresponding to the three trivial motions, and hence the maximum rank of $\mathbf{R}\pofw$ is $3v+3$.

If we wish to consider only the lattice-fixing infinitesimal motions of $\pofw$, then we may omit the final six columns corresponding to $\mathcal L(i,j; n_k)$. Our rigidity matrix is then of dimension $e \times 3v$, with maximum rank $3v -3$. In general, the maximal rank of a fixed-lattice rigidity matrix is $dv - d$, and the maximal rank of a flexible-lattice rigidity matrix is $dv - d + {d+1 \choose 2}$ \citep{BorceaStreinuI}.

As in the symmetric setting, there also exists a modified notion of generic for periodic frameworks (for further details see \citet{Ross}). For our purposes, it will be important to know only that generic rigidity of the periodic orbit framework $\pofw$ depends on the underlying periodic orbit graph $\pog$. Furthermore, for generic frameworks, infinitesimal rigidity and rigidity are equivalent \citep{Theran, Ross}.

As shown in \citet{BorceaStreinuI, Theran, Ross}, the vector space of periodic infinitesimal motions of a periodic framework  $\pfw$ is equivalent to the vector space of infinitesimal motions of the periodic orbit framework $\pofw$. The following theorem states that the vector space of periodic infinitesimal motions of a periodic framework corresponds to the kernel of the periodic rigidity matrix.

\begin{theorem} \citep{BorceaStreinuI, Ross} Let $\pofw$ be a periodic orbit framework corresponding to the periodic framework $\pfw$. The kernel of the corresponding periodic rigidity matrix $\mathbf{R}\pofw$ is isomorphic to the space of periodic infinitesimal motions of the associated periodic framework $\pfw$.
\label{thm:periodicKernel}
\end{theorem}

\begin{cor}\citep{BorceaStreinuI}
The periodic framework $\pfw$ is infinitesimally periodic rigid in $\mathbb R^d$ if and only if the rank of the rigidity matrix for the corresponding orbit framework $\mathbf{R}\pofw$ is $dv - d +  {d+1 \choose 2}$.
\label{cor:periodicRank}
\end{cor}

Returning to Example 4.1, the rank of the matrix corresponding to generic positions of the vertices is 5, which is maximal on two vertices, and hence $\pofw$ is infinitesimally rigid. If we are only interested in the rigidity of the framework on a fixed lattice,  the rank of the lattice-fixing portion of $\mathbf{R}\pofw$ (the first four columns) is 2, which again is maximal. Note that this means that three of the five edges of our example are redundant on a fixed lattice.

From Theorem \ref{thm:periodicKernel} and a periodic version of Theorem \ref{thm:flexes} we obtain a periodic Maxwell type counting rule for detecting finite periodic flexes:

\begin{theorem} Let $\pfw$ be a periodic framework in dimension $d$ with a corresponding orbit framework $\pofw$, where $v = |V(G)|$ and $e=|E(G)|$. If $$e < dv - d + {d+1 \choose 2}\ =\  dv + {d \choose 2},$$
then $\pofw$ has an infinitesimal flex, which corresponds to a periodic infinitesimal flex of $\pfw$.

Furthermore, for generic positions of the vertices of $G$ relative to the generating lattice $L$, $\pofw$ has a finite flex, which corresponds to a periodic finite flex of $\pfw$.
\label{thm:periodicMaxwell}
\end{theorem}

Theorem \ref{thm:periodicMaxwell} can be adapted for the fixed lattice with the count $e < dv - d$, or for any other variation of the flexible lattice. For dimensions $2$ and $3$, Table \ref{table{table:latticeParameters}} shows the number of lattice parameters corresponding to each of the lattice variants in the following list:
\begin{enumerate}
\item [(i)] {\em fully flexible lattice}:  all variations of the lattice shape are permitted;
\item [(ii)] {\em distortional change}: keep the volume fixed but allow the shape of the lattice to change;
\item [(iii)] {\em scaling change}: keep the angles fixed but allow the scale of the translations to change independently;
\item [(iv)] {\em hydrostatic change}: keep  the shape of the lattice unchanged but scale to change the volume;
\item [(v)] {\em fixed lattice}: allow no change in the lattice.
\end{enumerate}

\begin{table}[h!]
 \centering
 \caption{Number of parameters corresponding to types of lattice deformations with no added symmetry, in two and three dimensions.}
\begin{tabularx}{1.\textwidth}{X}
\[ \begin{array}
{|c|c|c|c|c|c|c|c|} \hline
 \grey{Lattice Def}&   \grey{2-D}   & \grey{3-D}\\ \hline
\text{flexible}& 3 & 6 \\
\hline
\text{distortional}& 2 & 5\\
\hline
\text{scaling}&  2&3\\
\hline
\text{hydrostatic}& 1&1 \\
\hline
\text{fixed} & 0& 0 \\
\hline
\end{array}  \]
\end{tabularx}
\label{table{table:latticeParameters}}
\end{table}

Why might we study one of these variants?  For crystals, we might focus on short time-scale vibrations, during which the large motions of distant atoms needed for a flexible lattice could not happen. In this case we effectively study a fixed lattice with local variation and all velocities small.  Or we might study slow responses to general pressure, given by a fully flexible lattice.  In between, we could study responses to  pressures and constraints of various types, with various boundary conditions, which correspond to various intermediate situations.

Another setting which produces periodic structures is simulation of large sphere packings by simulations with a modest number of spheres, and a periodic bounding box to give a better approximation than a fixed boundary.   Here, which case applies will depend on the variation of the periodic bounding box which the simulation chooses to permit.

\section{Periodic frameworks with symmetry}
\label{sec:PeriodicSymmetry}
In the previous two sections we have built up the orbit matrix  for finite frameworks under point groups ${\cal S}$ in $\R^{d}$ and for periodic frameworks with groups $\Z^{d}$.   Counting the rows and columns of these orbit matrices, the Maxwell type counts of rows vs columns minus trivial motions give necessary conditions for a framework to be generically rigid (minimally rigid).   Recall that for $\Z^{d}$ we had several variants, ranging from the fully flexible lattice with ${d+1\choose 2}$ columns added for the lattice variables, to the fixed lattice with no columns added for lattice deformations.

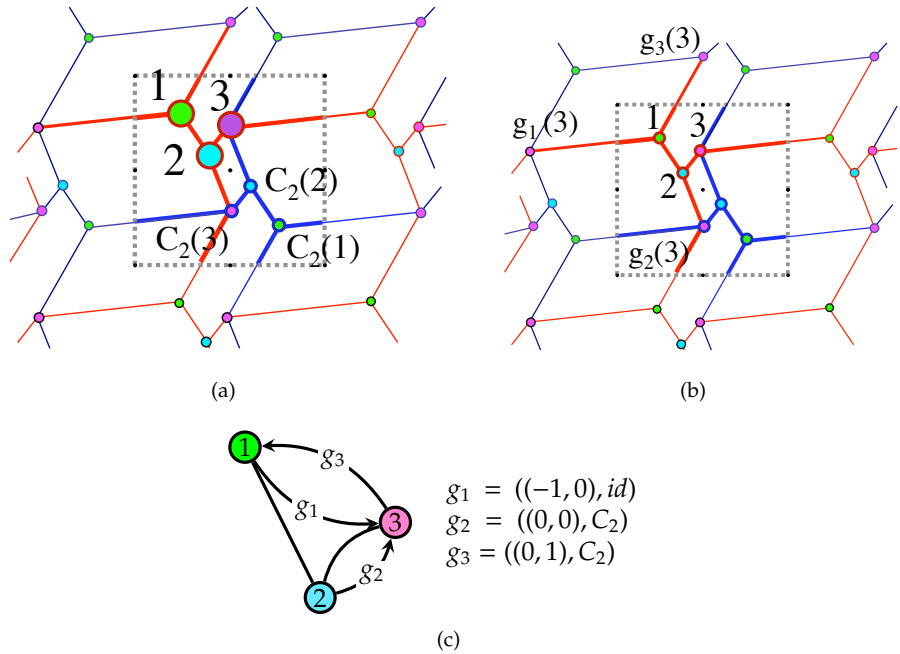
\begin{figure}[ht]
  \begin{center}
    \subfigure[] {\includegraphics [width=.5\textwidth]{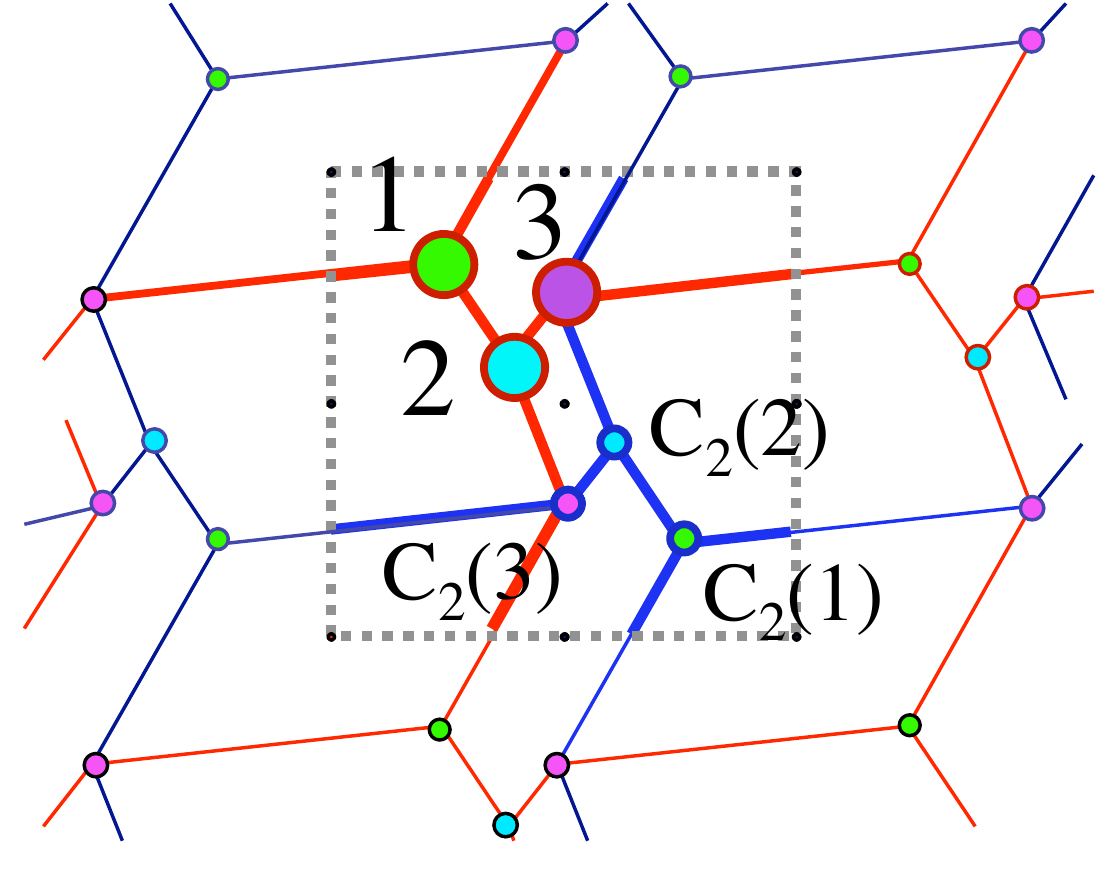}} \quad
   \subfigure[] { \includegraphics[width=.45\textwidth]{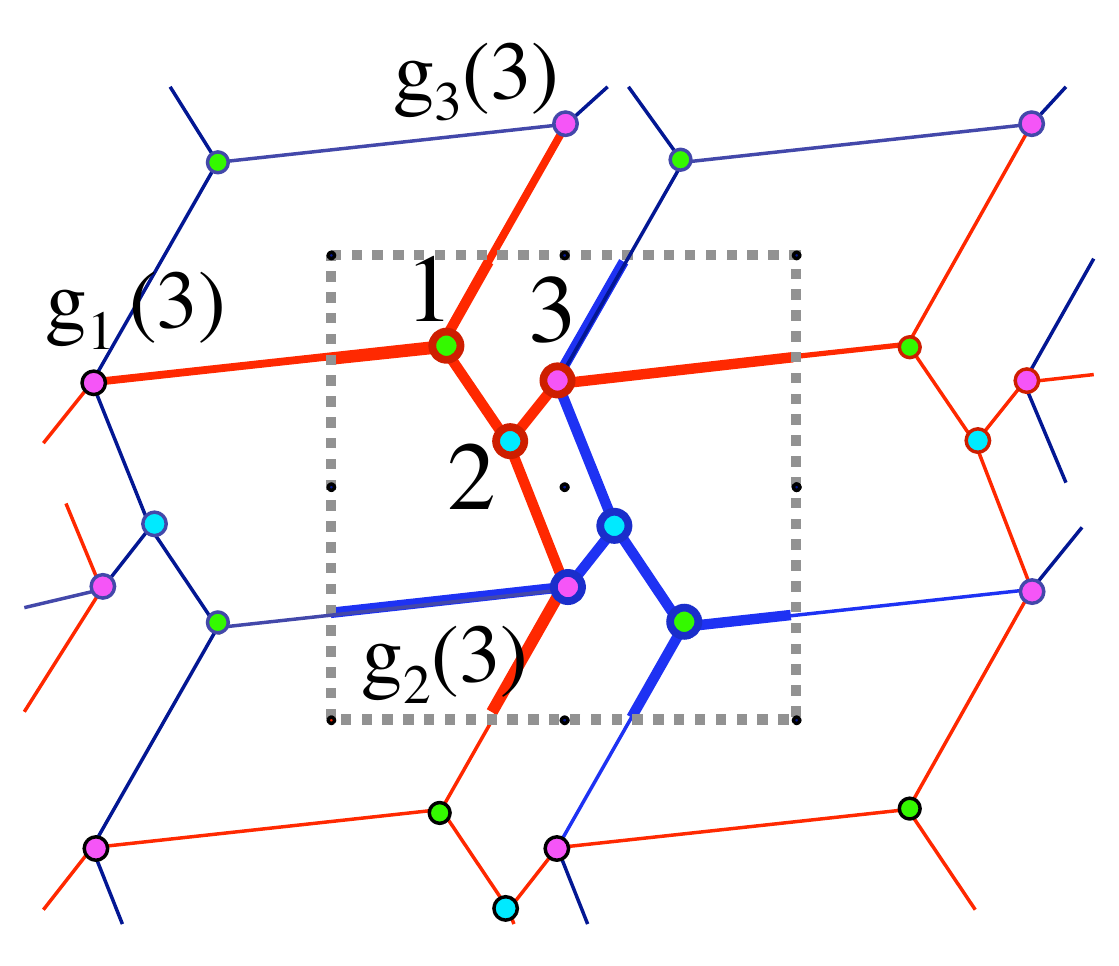}}
   \subfigure[] { \begin{tikzpicture}[auto, node distance=2cm, very thick]
\tikzstyle{vertex1}=[circle, draw, fill=green, inner sep=1pt, minimum width=3pt];
\tikzstyle{vertex2}=[circle, draw, fill=sky, inner sep=1pt, minimum width=3pt];
\tikzstyle{vertex3}=[circle, draw, fill=pink, inner sep=1pt, minimum width=3pt];
\tikzstyle{vertex4}=[circle, draw, fill=green, inner sep=1pt, minimum width=3pt];
\tikzstyle{gain} = [fill=white, text=black, inner sep =1pt,  anchor=center];
\node[vertex1] (1) {$1$};
\node[vertex2] (2) at (1, -2) {$2$};
\node[vertex3] (3) at (2, -1) {$3$};
\path (1) edge (2)
(2) edge [bend left] (3);
\pgfsetarrowsend{stealth}[ shorten >=2pt]
\path (1) edge [bend right] node[gain] {$g_1$} (3)
(2) edge [bend right] node[gain] {$g_2$} (3)
(3) edge [bend right] node[gain] {$g_3$} (1) ;
\pgfsetarrowsend{}
\node[text width = 1.2in] at (4.2, -1) {$g_1 =  ( (-1,0), id)$ $g_2 = ( (0,0), C_2)$ $g_3 =( (0,1), C_2)$};
\end{tikzpicture}}
     \end{center}
    \caption{A plane framework with $\Z^{2} \rtimes \mathcal{C}_{2}$ symmetry can be labeled with the elements of the group (a), or in short hand with `gains' (b) as in the gain graph (c)  \label{fig:gainsample}}
    \end{figure}

We now turn our attention to periodic frameworks with added symmetry. These frameworks have orbit graphs whose edges are labeled by elements of groups of the form $\Z^{d}\rtimes {\cal S}$. An example of such a framework with its corresponding orbit graph is shown in Figure \ref{fig:gainsample}.
\begin{defin}\label{def:symmetricPeriodicOrbitGraph}
Let $\pfw$ be a periodic framework with symmetry group $\Z^{d}\rtimes {\cal S}$. The \emph{symmetric periodic orbit graph} $\spog$ corresponding to this framework is the labelled multigraph with one representative for each equivalence class of edges and vertices under the action of $\Z^{d}\rtimes {\cal S}$. The labelling of the edges $g:E(G) \rightarrow \Z^{d}\rtimes {\cal S}$ is determined in the manner described in Definitions \ref{defsymorbgr} and \ref{def:periodicOrbitGraph}.
\end{defin}

\begin{figure}[h!]
\begin{center}
\includegraphics[width=3.5in]{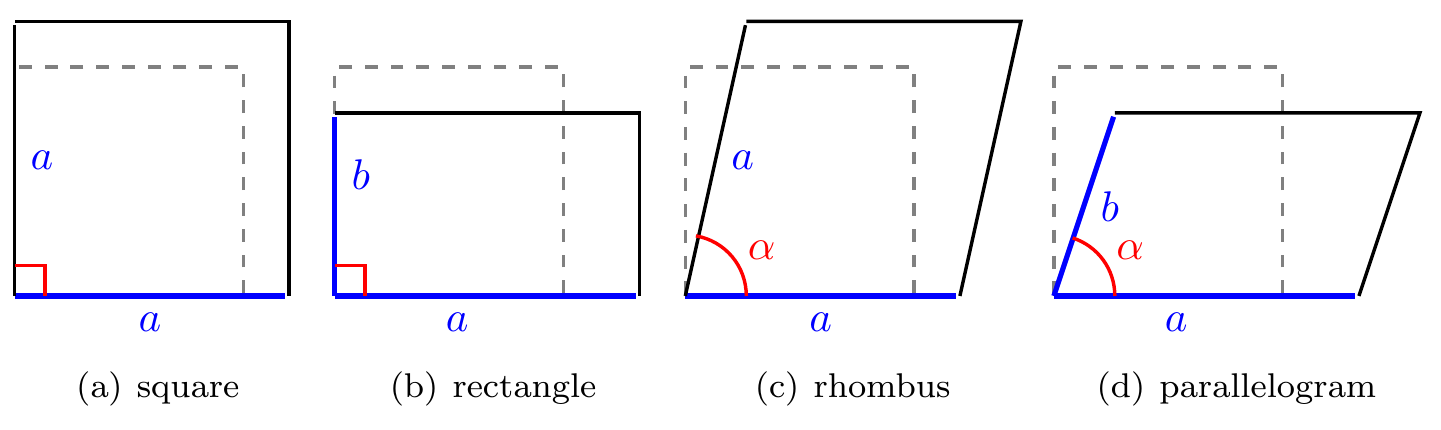}
\caption{The four planar crystal systems. The number of lattice parameters are (a) 1, (b) 2, (c) 2, (d) 3. \label{fig:2DcrystalSystems}}
\end{center}
\end{figure}

\begin{figure}[h!]
\begin{center}
\includegraphics[width=5in]{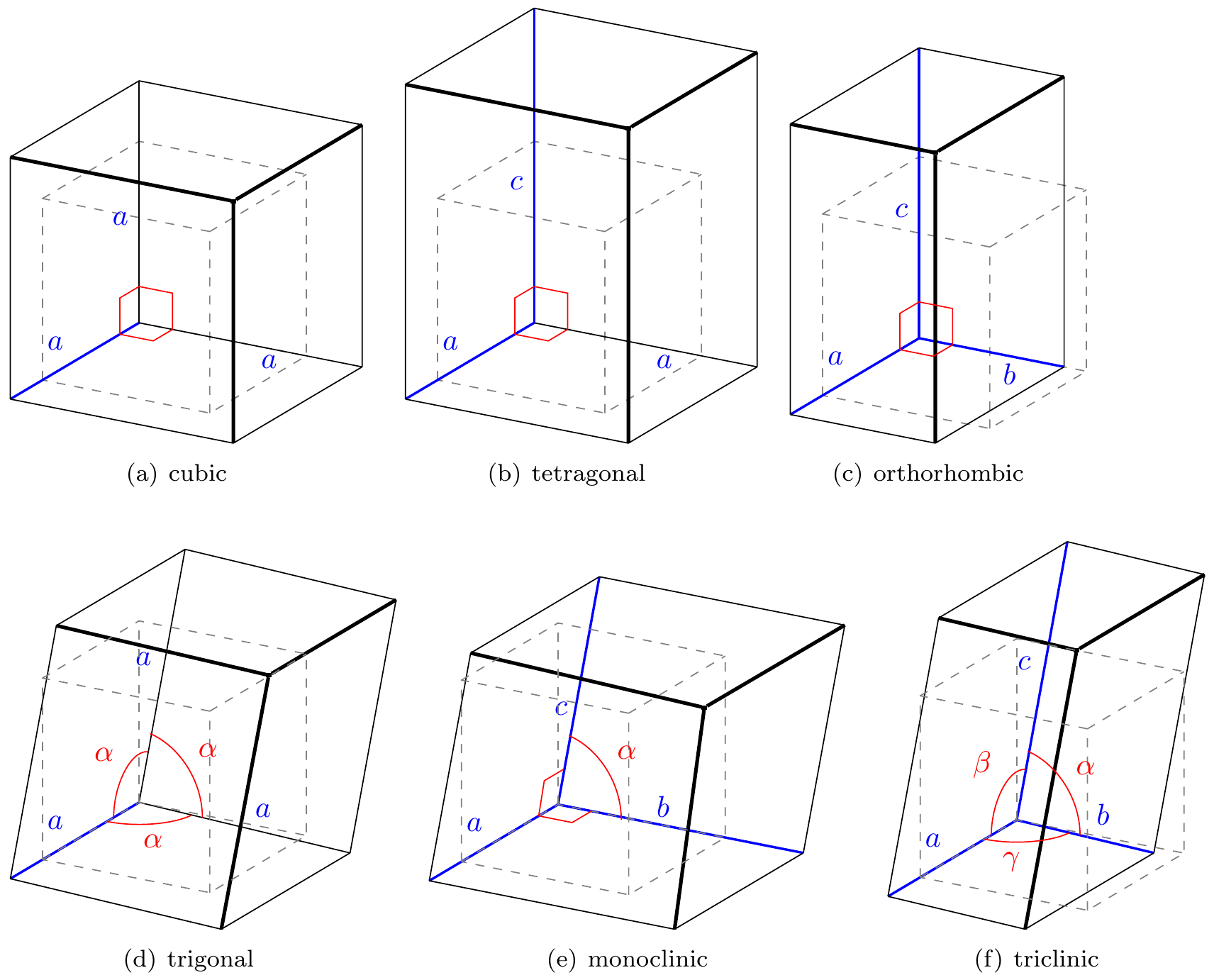}
\caption{The six crystal systems addressed in this work. The number of lattice parameters are (a) 1, (b) 2, (c) 3, (d) 2, (e) 4, (f) 6. \label{fig:3DcrystalSystems}}
\end{center}
\end{figure}

The symmetry group $\Z^{d}\rtimes {\cal S}$ of a symmetric periodic framework will determine its {\it crystal system}, which is a characterization of the parameters which determine the unit cell. This is usually defined by the number and arrangement of lengths and angles determining the cell, and these parameters represent the variations of lattice shapes which preserve the given symmetry \citep{Wiki, Tables}.  In the plane we will consider four different crystal systems, shown in Figure \ref{fig:2DcrystalSystems}, and in space we consider the six crystal systems shown in Figure \ref{fig:3DcrystalSystems}.

The crystal system of a framework specifies the maximum number of parameters that determine the lattice, and therefore the number of lattice columns of our orbit rigidity matrix. We may further reduce the number of lattice columns by changing the type of lattice that we are considering: (flexible, distortional, scaling, hydrostatic or fixed), although it should be noted that the lattice system will partially determine these choices. For instance, for a two-dimensional framework with a rhombus unit cell, scaling and hydrostatic will be identical.

\begin{remark} If we transform the lattice to the unit cube, by an affine transformation which preserves the symmetry, then the lattice parameters represent the number of non-zero partial derivatives of the length of bars, for variables from the lattice pattern. \end{remark}

The rigidity of the symmetric periodic framework can be studied using the orbit matrix. Letting $v_0$ and $e_0$ represent the number of vertices and edges in the orbit graph $G$, the orbit matrix has dimension $e_0 \times (dv_0 + \ell_{\mathcal S})$, where $\ell_{\mathcal S}$ describes the number of columns corresponding to the lattice parameters. This number will vary depending on a) the crystal system corresponding to the symmetry group $\mathcal S$, and b) the type of lattice variation we are considering (i.e. flexible, distortional, scaling, hydrostatic, or fixed).

As with the orbit matrices for ${\cal S}$, the second number which is important is the number of trivial motions which preserve all the group operations.  Since we are working with periodic frameworks, we are looking for what translations also preserve the symmetries in ${\cal S}$ within the orbit matrix.  We denote by $t_{{\cal S}}$ the dimension of the space of points which is fixed by all elements of the group ${\cal S}$.   This space is also called the symmetry element of the group.  For our calculations here, this  can only be one point, a line, a plane, or all of $3$-space.

It is the combination of these two numbers $\ell_{\cal S}$ and $t_{{\cal S}}$, plus the number of orbits of edges and vertices (corresponding to the number of edges $e_0$ and vertices $v_0$ of the orbit graph $\spog$), which generates the predictions of the number of non-trivial motions (if any) which occur when the vertices of the framework are in generic position. As before, we assume that no edges or vertices are fixed by the action of the group $\Z^{d}\rtimes {\cal S}$.

In Section \ref{sec:PlanePeriodicSymmetry} we describe symmetric periodic frameworks with two samples of ${\cal S}$  in the plane in detail, followed by tables covering all the groups within our analysis.   In Section \ref{sec:SpacePeriodicSymmetry} we outline two samples of ${{\cal S}}$ in 3-space, again followed by tables for the relevant groups.

In the tables we  give three distinct notations for each group: the Schoenflies notation  used by chemists \citep{bishop, Hall}, the Hermann-Mauguin notation  used internationally by crystallographers \citep{Tables}, and the orbifold notation used by mathematicians \citep{Conway}.

Recall that the Schoenflies notation was briefly introduced in Section \ref{subsec:SymmetricOrbitMatrices1}.
 
In the Hermann-Mauguin notation, $n$-fold rotational symmetries are denoted by $n$, and these axis numbers are written down in decreasing order of $n$. If an $n$-fold rotational axis is contained in a mirror, then an $m$ is written after the corresponding axis number. If there exists a mirror which is perpendicular to an $n$-fold rotational axis, then this is denoted by placing the symbol $/m$ after the corresponding axis number. Finally, the notation $\overline{n}$ for an axis number $n$ indicates that an inversion, followed by the $n$-fold rotation, is a symmetry of the structure.

In the orbifold notation for wallpaper groups in the plane, an $n$-fold rotational symmetry is denoted by $n$, and a mirror symmetry is denoted by an asterisk, $*$. A number before an asterisk indicates a center of pure rotation, whereas a number after an asterisk indicates  a center of rotation with mirrors through it. If there are only translational symmetries present, then this is denoted by the symbol $\circ$. The orbifold notation for groups in $3$-space is similar; note, however, that for a group in $3$-space, an $x$ indicates the presence of an inversion symmetry, whereas in $2$-space, $x$ indicates a glide reflection.
  
For an interdisciplinary audience, the simultaneous use of these three notations seems to be an appropriate presentation, and such comparative columns can be found in multiple sources, including the Wikipedia pages for crystal systems \citep{Wiki}.

Note that the assumption that the group has the form $\Z^{d}\rtimes {\cal S}$ means that we do not cover all the plane wallpaper groups, or the full set of space groups.   Specifically, we will not include groups which have glide reflections as generators of the group, or those which have $6$-fold rotations.   These other groups will also have orbit matrices, but require an alternative analysis for comparisons and counts.  We return to this issue in Section \ref{subsec:OtherGroups}.

\section{$2$-D periodic frameworks with symmetry: $\Z^{2}\rtimes \cal S$}
\label{sec:PlanePeriodicSymmetry}
\subsection{$\Z^{2} \rtimes \mathcal{C}_{2}$ - half-turn symmetry in the plane lattice}
\label{subsec:PlanePeriodicHalf}
Half-turn symmetry in the plane is equivalent to inversion in the point axis. This symmetry fits an arbitrary parallelogram for the lattice (Figure \ref{fig:2DcrystalSystems}(d)), and $\ell_{\mathcal{C}_{2}}=3$.
We will consider periodic plane frameworks with symmetry $\Z^{2} \rtimes \mathcal{C}_{2}$ for two variations of the lattice: (1) a fully flexible lattice; (2) a fixed lattice.

\medskip

\noindent {\bf Example 6.1.1: fully flexible lattice $\Z^{2} \rtimes \mathcal{C}_{2}$.}
The original (non-symmetric) necessary count for a periodic framework on the fully flexible lattice to be minimally rigid is
$e=2v+1$ (recall Theorem \ref{thm:periodicMaxwell}).  To permit half-turn symmetry, with no vertex or edge fixed by the half-turn, we will need to start with the modified count  $2e_{0}=2(2v_{0})+2$, where $v_0$ and $e_0$ are the numbers of vertices and edges of the orbit graph, respectively. Dividing by $2$, this gives  $e_{0}=2v_{0}+1$.

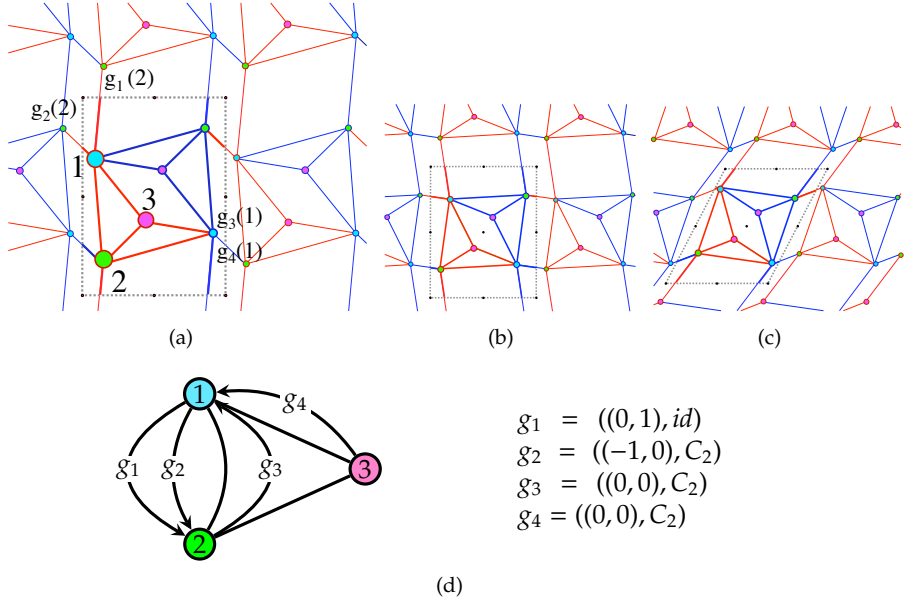
\begin{figure}[ht]
  \begin{center}
   \subfigure[] {\includegraphics [width=.40\textwidth]{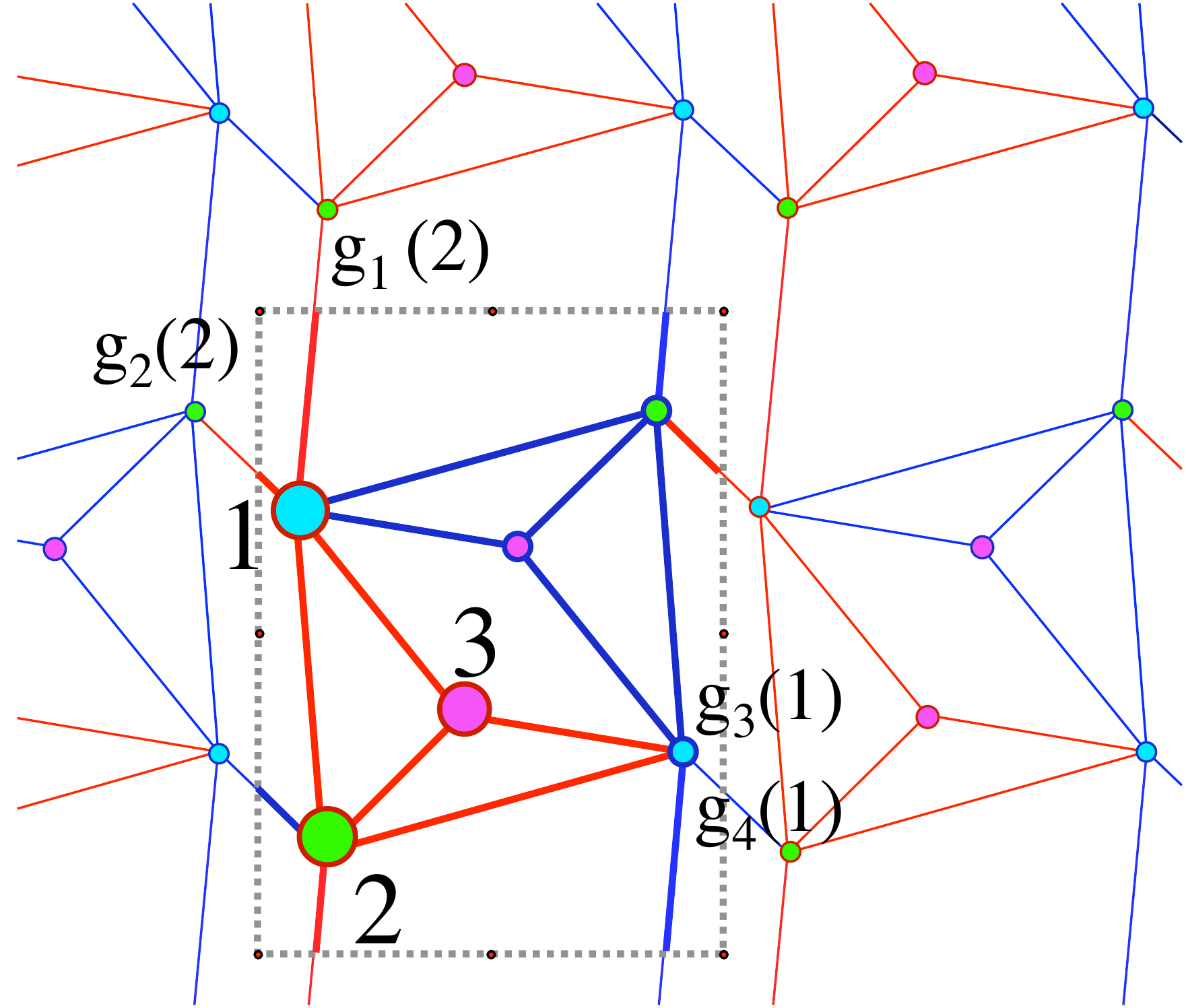}}
   \subfigure[] { \includegraphics[width=.28\textwidth]{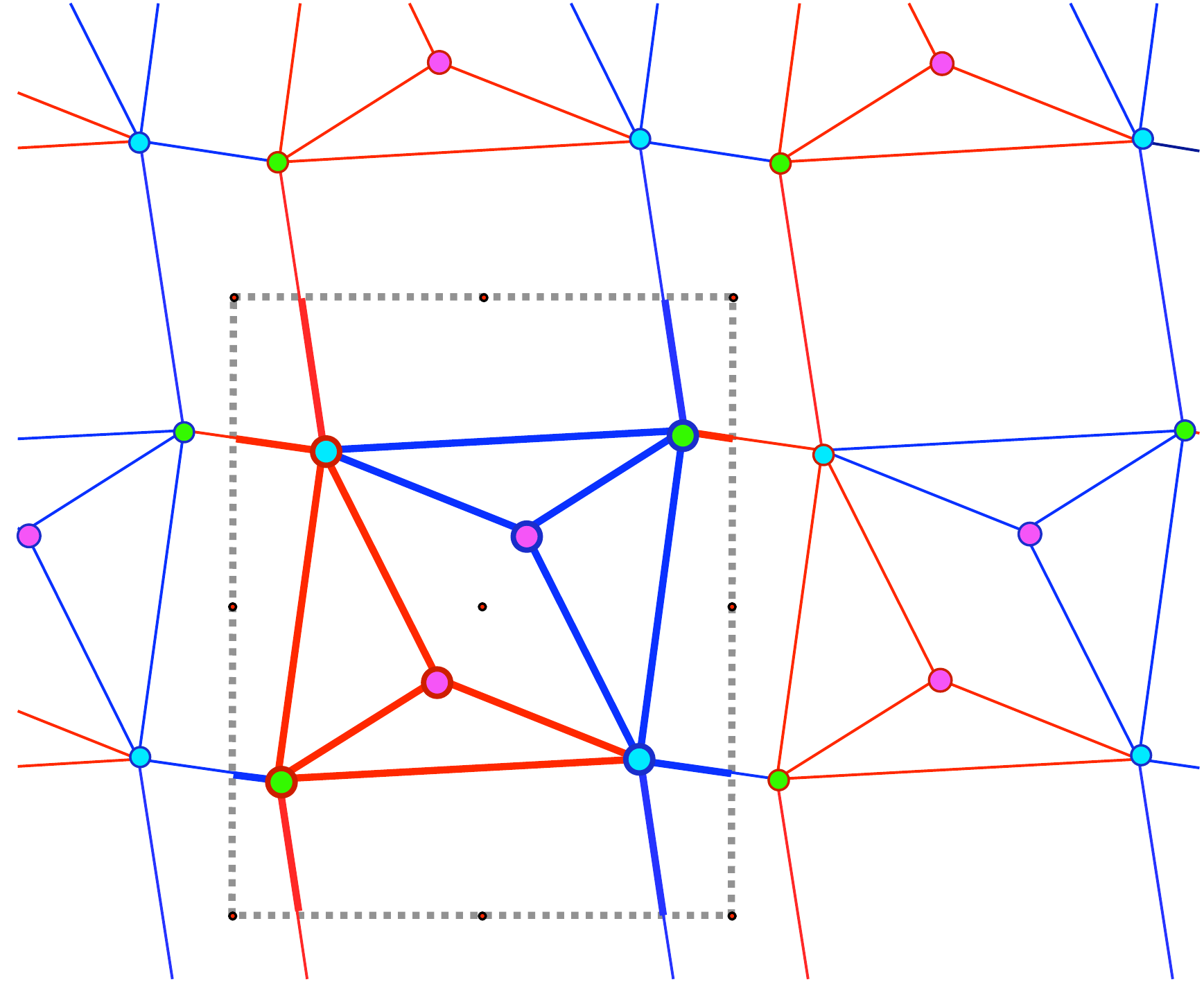}}
    \subfigure[] { \includegraphics[width=.28\textwidth]{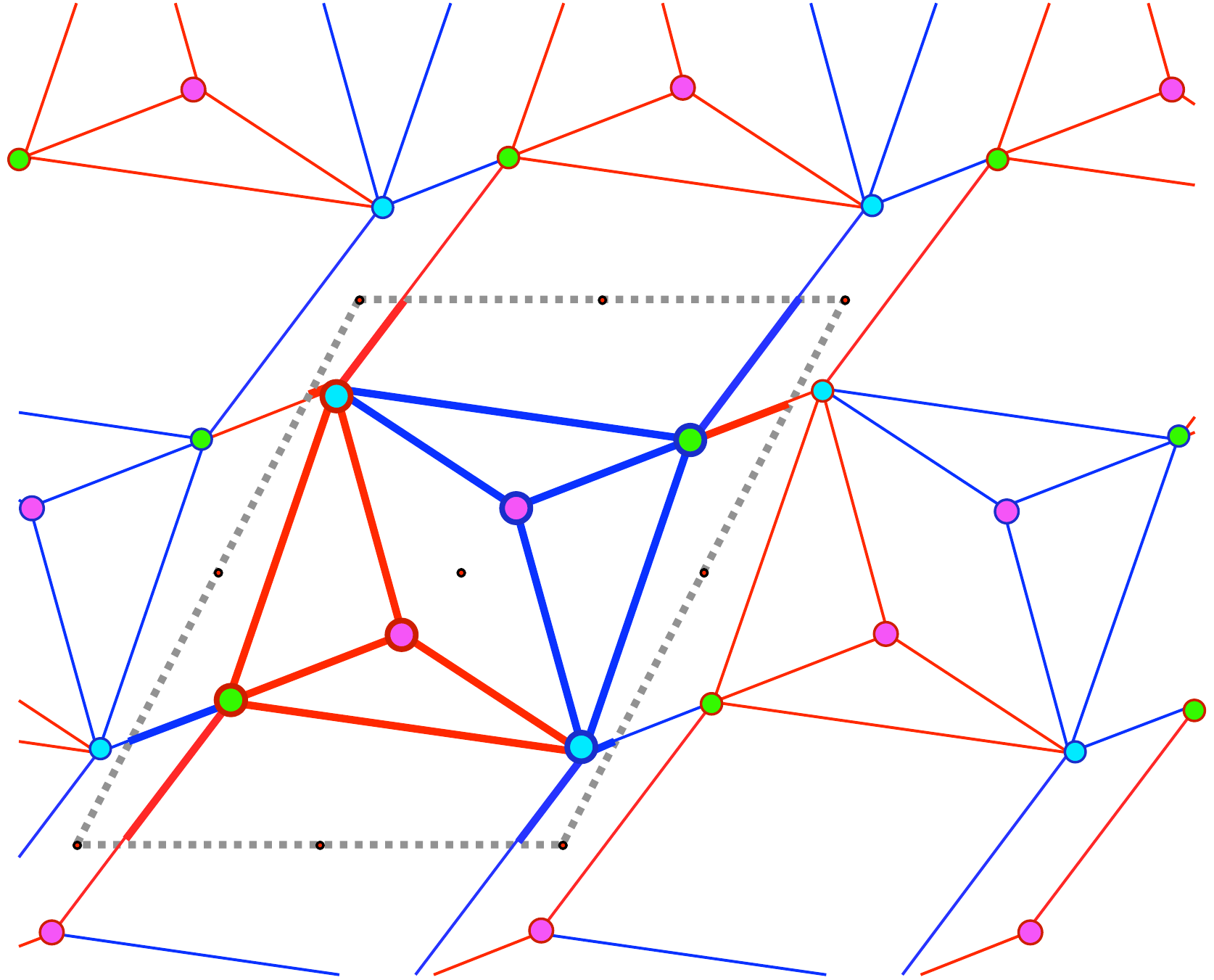}}
     \subfigure[] {\begin{tikzpicture}[auto, node distance=2cm, very thick]
\tikzstyle{vertex2}=[circle, draw, fill=green, inner sep=1pt, minimum width=3pt]; \tikzstyle{vertex1}=[circle, draw, fill=sky, inner sep=1pt, minimum width=3pt]; \tikzstyle{vertex3}=[circle, draw, fill=pink, inner sep=1pt, minimum width=3pt]; \tikzstyle{vertex4}=[circle, draw, fill=green, inner sep=1pt, minimum width=3pt]; \tikzstyle{gain} = [fill=white, text=black, inner sep =1pt,  anchor=center];

\node[vertex1] (1) {$1$};
\node[vertex2] (2) at (0, -2) {$2$};
\node[vertex3] (3) at (2.2, -1) {$3$};

\draw (2) -- (3) -- (1);

\pgfsetarrowsend{stealth}[ shorten >=2pt]

\draw (1) edge [bend right] node[gain] {$g_2$}
(2);
\draw (1) .. controls (-1.2,-.7) and (-1.2, -1.3)  .. node[gain] {$g_1$} (2);

\draw (2) .. controls  (1.2, -1.3) and (1.2,-.7)  .. node[gain] {$g_3$} (1);
\draw (3) edge [bend right] node[gain] {$g_4$} (1);
\pgfsetarrowsend{}
\draw (1) edge [bend left] 
(2);

\node[text width = 1.4in] at (6, -1)
{$g_1 = ((0,1), id)$
$g_2 = ((-1,0), C_2)$
$g_3 = ((0,0), C_2)$
$g_4 = ((0,0), C_2)$};
\end{tikzpicture}
}
   \end{center}
    \caption{A generically rigid graph on a fully flexible lattice, realized with 2-fold symmetry  has several non-trivial flexes changing the lattice. Its periodic symmetric orbit graph is pictured in (d).
 }
    \label{fig:2FoldFlexF}
    \end{figure}

Under the half-turn symmetry with a fully flexible  lattice, the orbit matrix has $2$ columns under each orbit of vertices, plus $\ell_{\mathcal{C}_{2}}=3$ columns for the three parameters for the lattice deformations.
 Further, we clearly have $t_{\mathcal{C}_{2}}=0$ since there are  no infinitesimal  trivial motions which preserve the half-turn symmetry along with the periodic lattice.   This creates the necessary symmetric Maxwell condition
$$e_{0}\ \geq \ 2v_{0}+3-0$$
for periodic rigidity.
However, as shown above, for a graph that was previously minimally rigid without the symmetry, we have $e_{0}=2v_{0}+1\ < \   2v_{0}+3$.
This gap predicts that a graph which counted to be minimally rigid without symmetry, realized generically with half-turn symmetry on a fully flexible lattice, now has two degrees of (finite) flexibility.
As an example, consider the snapshots of three configurations with the same edge lengths but changing angles and lengths of the unit cell in Figure~\ref{fig:2FoldFlexF}.  Together these snap shots confirm the predicted two degrees of freedom.

The orbit matrix corresponding to the framework pictured in Figure \ref{fig:2FoldFlexF} has the following form:

        \begin{displaymath}
  \renewcommand{\arraystretch}{0.8}
     \bordermatrix{  &p_1 &  p_2 &  p_3 & (t_{11}, t_{21}, t_{22})
     \cr
    \{1,2\} &  \big(p_1-p_2\big)  & \big(p_2-p_1\big) & 0 &(0,0,0)
     \cr
    \{2,3\} &  0& \big(p_2-p_3\big)  & \big(p_3-p_2\big) &(0,0,0)
     \cr
     \{3,1\} &   \big(p_1-p_3\big)&0  & \big(p_3-p_1\big) &(0,0,0)
     \cr
    \{1,2; g_1\} &  \big(p_1-g_1(p_2)\big)&  \big(p_2-g_1^{-1}(p_1)\big)  &0 &(\ast, \ast, \ast)
     \cr
   \{1,2; g_2\} & \big(p_1-g_2(p_2)\big)&  \big(p_2-g_2^{-1}(p_1)\big)  &0 &(\ast, \ast, \ast)
     \cr
    \{2,1; g_3\} & \big(p_1 -g_3^{-1}(p_2) \big)  & \big(p_2 - g_3(p_1)\big) & 0 &(\ast, \ast, \ast)
   \cr
    \{3,1; g_4\} & \big(p_1 -g_4^{-1}(p_3) \big)  &0 &  \big(p_3 - g_4(p_1)\big) &(\ast, \ast, \ast)
   \cr
     } \textrm{.}
    \end{displaymath}
     \eop

\medskip

\noindent {\bf Example 6.1.2: fixed lattice $\Z^{2} \rtimes \mathcal{C}_{2}$.} The original (non-symmetric) necessary count for any minimally rigid periodic framework on the fixed lattice is  $e=2v-2$.
   With added $\mathcal{C}_{2}$ symmetry, we have $e=2e_{0}$ and $v=2v_{0}$ (all orbits have $k_{\mathcal{C}_2}=2$), so a minimally rigid orbit graph, realized with  $\mathcal{C}_{2}$ symmetry, will have $2e_{0} =2(2v_{0})-2$, or
$e_{0}=2v_{0}-1$.

\begin{figure}[ht]
  \begin{center}
  \subfigure[] {\includegraphics [width=.445\textwidth]{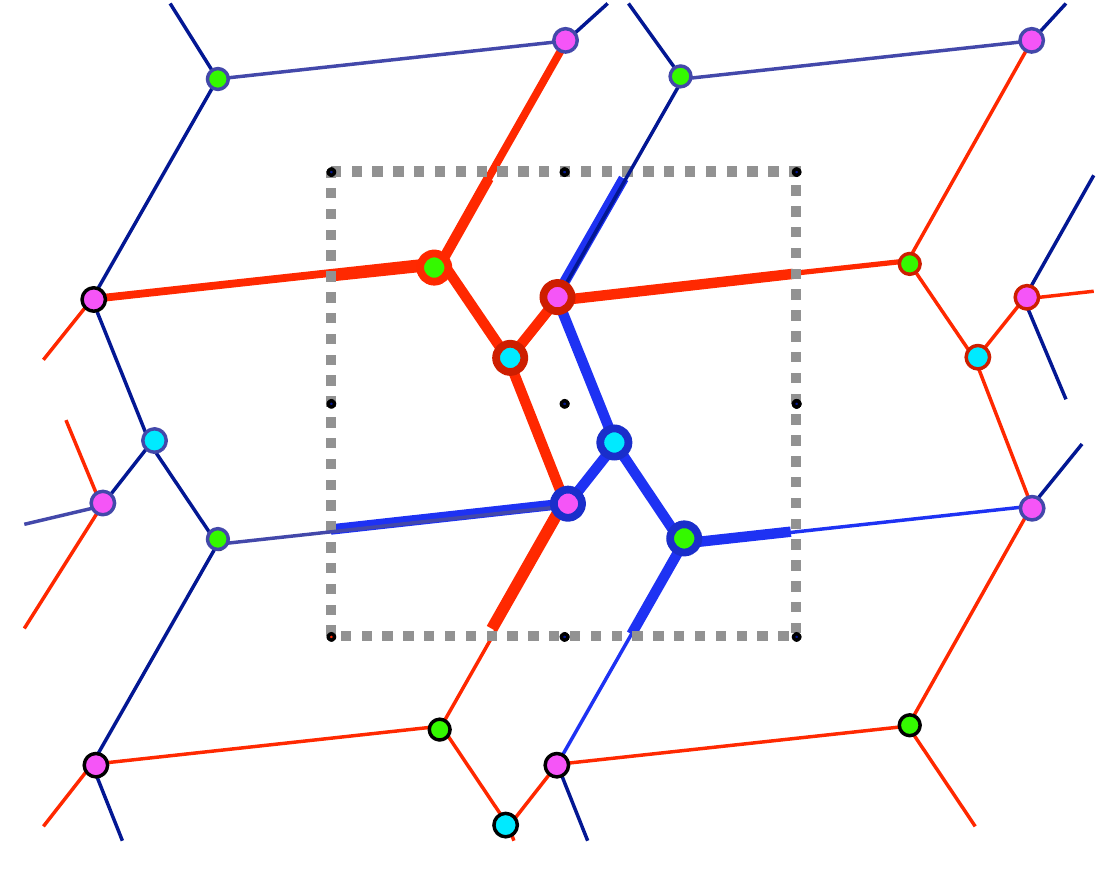}} \quad
   \subfigure[] { \includegraphics[width=.48\textwidth]{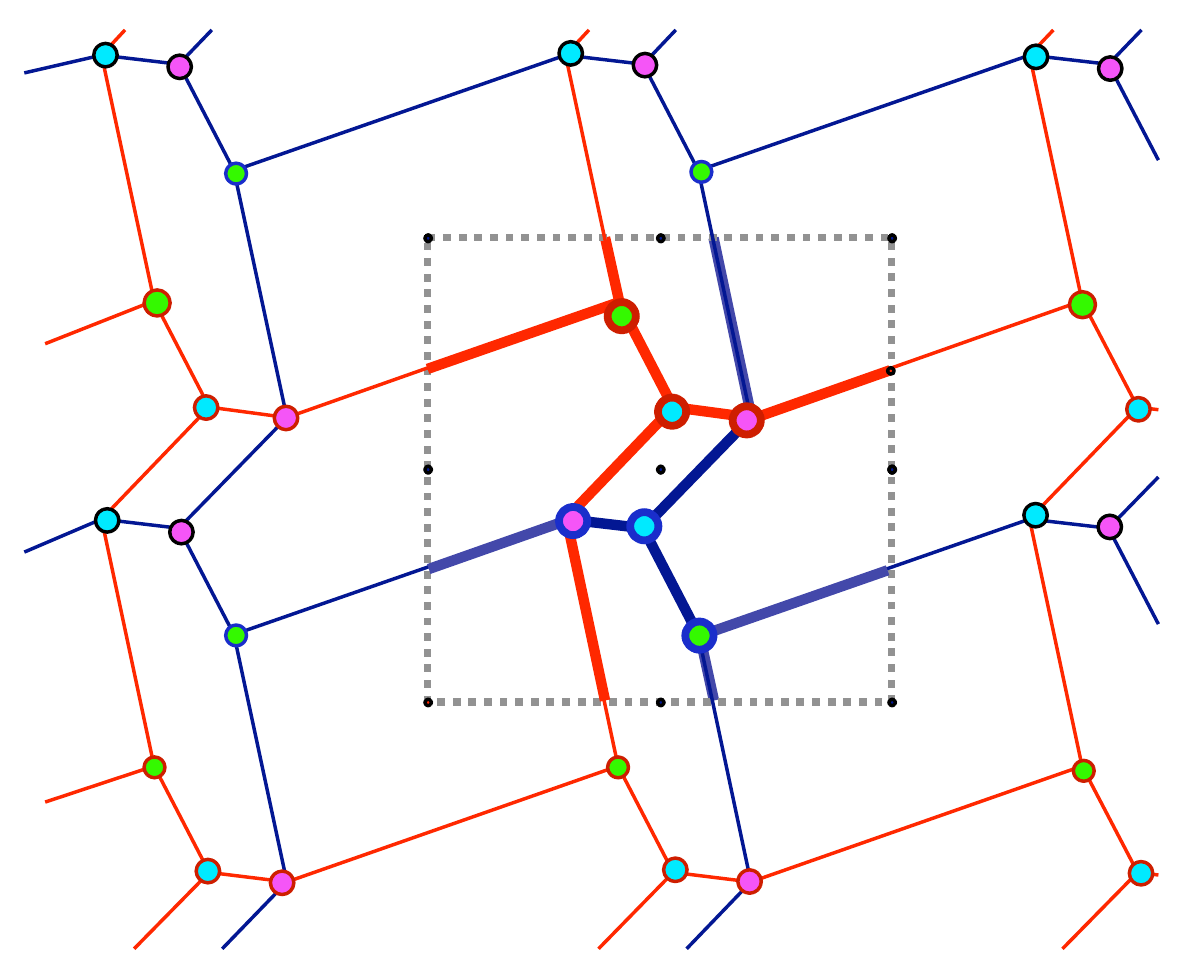}}
     \end{center}
    \caption{A plane framework with $\Z^{2} \rtimes \mathcal{C}_{2}$ symmetry has a non-trivial flex on the fixed lattice.}
    \label{fig:2FoldFixed}
    \end{figure}

As the example below illustrates, with a fixed lattice, the orbit matrix has $2$ columns under each orbit of vertices.  Further, there are no translations which preserve the half-turn symmetry along with the periodic lattice, and hence we have $t_{\mathcal{C}_2}=0$.   This creates the necessary symmetric Maxwell condition
      $$e_0 \geq 2v_{0}$$
for periodic rigidity.
However, as shown above, if the graph was chosen to be minimally rigid without the symmetry, we have $e_0 = 2v_{0}-1$.
The gap $e_0 = 2v_{0}-1 < 2v_{0}$ shows that with the added half-turn symmetry, a minimally rigid graph will become flexible within the fixed lattice.  Figure~\ref{fig:2FoldFixed} shows the sample graph already presented in Figure~\ref{fig:gainsample} with $v_{0}=3$, $e_{0}= 5$ and $e_{0}= 5< 6= 2v_{0}$, with two realizations with the same edge lengths - illustrating snapshots of a non-trivial motion, as predicted.

Here is the orbit matrix of the framework depicted in Figure \ref{fig:2FoldFixed} on the fixed lattice, with joints $p_1, p_2, p_3$:
        \begin{displaymath}
  \renewcommand{\arraystretch}{0.8}
     \bordermatrix{  &p_1 &  p_2 &  p_3 &
     \cr
    \{1,2\} &  \big(p_1-p_2\big)  & \big(p_2-p_1\big) & 0 &
     \cr
    \{2,3\} &  0& \big(p_2-p_3\big)  & \big(p_3-p_2\big) &
     \cr
    \{1,3; g_1\} &  \big(p_1-g_1(p_3)\big)& 0  & \big(p_3-g_1^{-1}(p_1)\big) &
     \cr
   \{2,3; g_2\} & 0&  \big(p_2-g_2(p_3)\big)& \big(p_3-g_2^{-1}(p_2)\big) &
     \cr
    \{1,3; g_3\} & \big(p_1 -g_3(p_2) \big)  & \big(p_2 - g_3^{-1}(p_1)\big) & 0
   \cr
     } \textrm{.}
    \end{displaymath} \eop


In Table \ref{table{HalfTurnPlane}} we summarize the $(\Z^{2} \rtimes \mathcal{C}_{2})$-symmetric Maxwell type counts for each of the lattice variants. For simplicity at this stage, we again assume that no joint and no bar is fixed by the half-turn, so that all vertex orbits and edge orbits of the periodic orbit graph under the action of the group have the same size $k_{\mathcal{C}_2}=2$.  For each type of lattice deformation, we always assume that $e$ is chosen to be the least number of edges for the framework to be rigid without symmetry and to be compatible with the symmetry constraints given by $\Z^{2} \rtimes \mathcal{C}_{2}$. The number $f_{\mathcal{C}_2}$ in the final column denotes the dimension of the space of $(\Z^{2} \rtimes \mathcal{C}_{2})$-symmetric infinitesimal flexes in each case. For `generic' configurations, these extend to finite symmetry-preserving flexes.

\begin{table}[h!]
 \centering
 \caption{Plane lattice deformations with $\mathcal{C}_{2}$ symmetry.}
\begin{tabularx}{1.\textwidth}{X}
\[ \begin{array}
{|c|c|c|c|c|c|c|c|} \hline
 \grey{Lattice Def}&   \grey{\cal S}   & \grey{k_{\cal S}}  & \grey{t_{{\cal S}}}& \grey{\ell_{\cal S}}
 & \text{\grey {rows }}& \text{\grey {columns -$t_{\cal S}$}}& \grey{f_{\cal S}} \\ \hline
\text{flexible}& \mathcal{C}_2 & 2&  0&3&e_{0}=2v_{0}+1  &  2v_{0}+3&2\\
\hline
\text{distortional}& \mathcal{C}_2 & 2&  0&2&e_{0}=2v_{0}  &  2v_{0}+2&2\\
\hline
\text{hydrostatic}& \mathcal{C}_2 & 2&  0&1&e_{0}=2v_{0}-1  & 2v_{0}+1&2\\
\hline
\text{fixed}& \mathcal{C}_2 & 2& 0&0&e_{0}=2v_{0}-1 &  2v_{0}&1\\
\hline
\end{array}  \]
\end{tabularx}

\label{table{HalfTurnPlane}}
\end{table}

\subsection{ $\Z^{2}\rtimes \mathcal{C}_{s}$ - mirror symmetry in the plane lattice}
\label{subsec:PlanePeriodicMirror}

The mirror parallel to a side of the lattice restricts the possible lattices to rectangles.  
This mirror symmetry is preserved by translation along the line of the mirror, so $t_{\mathcal{C}_{s}}=1$.

We will consider periodic plane frameworks with symmetry $\Z^{2} \rtimes \mathcal{C}_{s}$ again in two layers: (1) a fully flexible lattice; (2) a fixed lattice.

\medskip

\noindent {\bf Example 6.2.1: fully flexible lattice $\Z^{2} \rtimes \mathcal{C}_{s}$.}    The original (non-symmetric) necessary count for any minimally rigid periodic framework on the fully flexible lattice is $e=2v+1$ (recall Theorem \ref{thm:periodicMaxwell}).   To permit mirror symmetry, with no vertex or edge fixed by the mirror, we will need to start with the shifted count $2e_{0}=2(2v_{0})+2$ or equivalently $e_{0}=2v_{0}+1$.

\begin{figure}[ht]
  \begin{center}
    \subfigure[]  { \includegraphics[width=.65\textwidth]{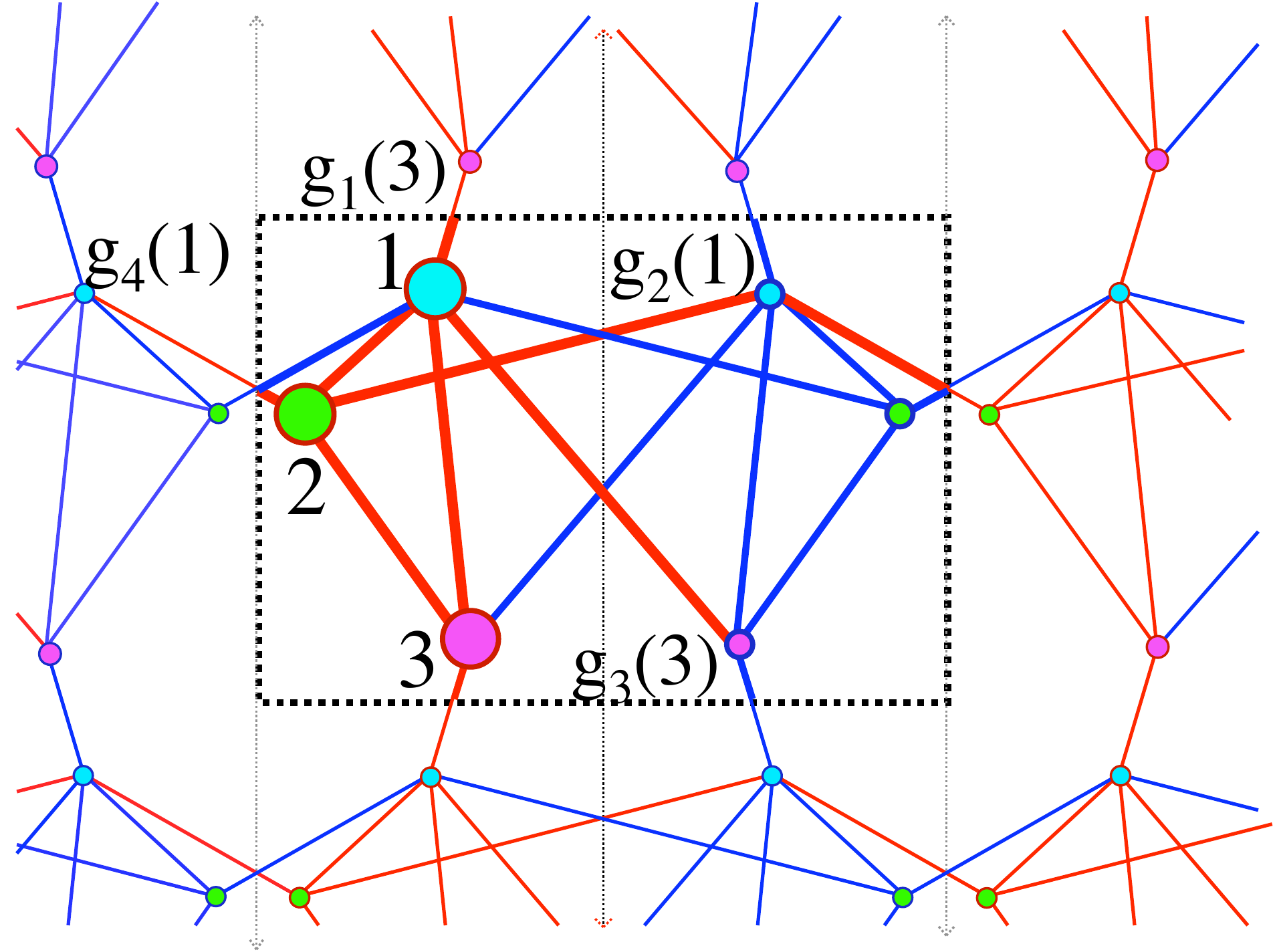}}%
  \subfigure[] {\begin{tikzpicture}[auto, node distance=2cm, very thick]
\tikzstyle{vertex1}=[circle, draw, fill=pink, inner sep=1pt, minimum width=3pt];
\tikzstyle{vertex2}=[circle, draw, fill=yellow, inner sep=1pt, minimum width=3pt];
\tikzstyle{vertex3}=[circle, draw, fill=sky, inner sep=1pt, minimum width=3pt];
\tikzstyle{vertex4}=[circle, draw, fill=green, inner sep=1pt, minimum width=3pt];
\tikzstyle{gain} = [fill=white, inner sep =1pt,  anchor=center];
\node[vertex3] (1) {$1$};
\node[vertex4] (2) at (-1.3, -1.2) {$2$};
\node[vertex1] (3) at (.4, -3) {$3$};
\path (1) edge (2);
\path(2) edge  (3);
\path(1) edge  (3);
\pgfsetarrowsend{stealth}[ shorten >=2pt]
\path  (1) edge [bend left] node[gain] {$g_3$} (3);
\path(2) edge [bend left] node[gain] {$g_2$} (1);
\path (2) edge [bend right] node[gain] {$g_4$} (1);
\draw (1) .. controls (1.6,-.8) and (1.8, -1.8)  .. node[gain] {$g_1$} (3);

\pgfsetarrowsend{}
\node[text width = 1in] at (0.5, -4.5) {$g_1 = ((0,1), id)$ $g_2 = ((0,0), \sigma)$ $g_3 = ((0,0), \sigma)$ $g_4 = ((-1,0), \sigma)$};
\end{tikzpicture}}
     \end{center}
    \caption{The mirrors (vertical lines in (a)) fit only with the two scalings and this framework prevents those scalings. The orbit graph corresponding to this framework is shown in (b).}
    \label{fig:PlaneMirrorF}
         \end{figure}
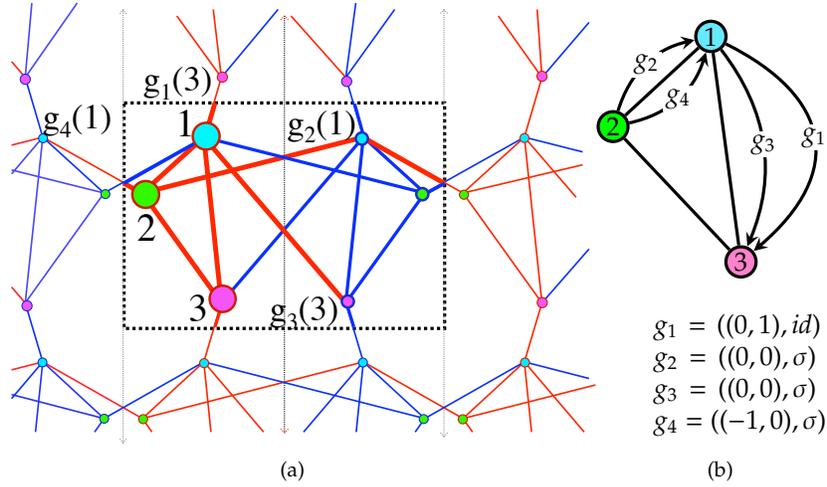

Under the mirror symmetry with a flexible lattice which preserves the symmetry, the orbit matrix has $2$ columns under each orbit of vertices, plus $\ell_{\mathcal{C}_{s}}=2$ lattice scaling columns for the mirror preserving flexes of the lattice.  Since $t_{\mathcal{C}_{s}}=1$, we have the necessary symmetric Maxwell condition
$$e \geq  2v_{0}+2-1 = 2v_{0}+1$$
for periodic rigidity.
This inequality, together with the previous condition for minimal rigidity without the mirror symmetry, suggests that there is no added flexibility from this mirror symmetry.  The example in Figure~\ref{fig:PlaneMirrorF} illustrates such a situation with $v_{0} =3$, $e_{0} =7$, and  $e_{0} = 7=2v_{0} +1$.  It is indeed rigid on a flexible lattice,  up to vertical translation along the mirror line.
 \eop

\medskip

\noindent {\bf Example 6.2.2: fixed lattice $\Z^{2} \rtimes \mathcal{C}_{s}$.} The original necessary count for a periodic framework on the fixed lattice to be minimally rigid is $e=2v-2$. With added mirror symmetry, we have $e=2e_{0}$ and $v=2v_{0}$, so a minimally rigid orbit graph, realized with  $\mathcal{C}_{s}$ symmetry, will have $2e_0 =2(2v_{0})-2$, or $e_{0}=2v_{0}-1$.

Under the mirror with a fixed lattice, the orbit matrix has $2$ columns under each orbit of vertices. Moreover, we have  $t_{\mathcal{C}_{s}}=1$ since the translation along the axis preserves the mirror symmetry along with the periodic lattice.   This creates the necessary symmetric Maxwell condition
      $$e_{0} \geq 2v_{0}-1$$
for periodic rigidity.
Together with the count for minimal rigidity without symmetry, this suggests that there is no added flexibility from this symmetry.
 \eop

It turns out that for mirror symmetry, all of the variants of lattice deformations produce no added motions.

\subsection{Table of groups for the fully flexible lattice in $2$-dimensions}
\label{subsec:PlaneTableFlex}

Examples 6.1.1 and 6.1.2 indicate a process that can be applied to other plane symmetries which preserve the lattice. Each row in Table \ref{table:PlaneFlex} presents the calculation for a given plane wall-paper group which is presented as
$\Z^{2}\rtimes \cal S$.  Recall that we are not including the plane wall-paper groups that have core glide reflections or $3$-fold and $6$-fold rotations, since they require some significant modifications of the simple pattern presented here (see also Section \ref{subsec:OtherGroups}). Thus, we do not have 17 lines in the table.

In each row of Table \ref{table:PlaneFlex}, the calculation has several parts - each producing an integer:

\begin{enumerate}
\item the number of edge orbits, $e_0$, so that   $k_{{\cal S}}e_{0}
\geq  2(k_{{\cal S}} v_{0}) +1$, which guarantees that we have at least the number of edges needed for a rigid periodic framework without symmetry.  This means we need to add a modified constant $\lceil{\frac{1} {k_{\cal S}}}\rceil$.  For Table \ref{table:PlaneFlex}, this value is always $1$, and the number of rows is always $e_{0}=2v_{0}+1$.

\item $t_{\cal S}$ which is the dimension of the space of translations contained in the symmetry element of ${\cal S}$.  This will be $2$ for the identity group, $1$ for a single mirror, and $0$ otherwise.
 \item $\ell_{\cal S}$ which is the dimension of the space of lattice deformations which preserve all the symmetries in ${\cal S}$ or equivalently, the number of independent parameters in the lattice system (edge lengths and angles).

\item the comparison of these numbers as the number of rows $e_{0}$ compared to the number of columns minus $t_{\cal S}$:
$2v_{0}+\ell_{\cal S} - t_{\cal S}$.
\item the difference $f_{\cal S}=2v_{0}+\ell_{\cal S} - t_{\cal S} - (2v_{0}+1)=\ell_{\cal S} - t_{\cal S} -1$ which is the dimension of the guaranteed extra non-trivial motions of the symmetric framework, over the rigidity which the original count without symmetry promised.
\end{enumerate}


\begin{table}[h!]
 \centering
 \caption{The added flexibility induced by basic symmetries on a fully flexible 2-D lattice for $\Z^{2}\rtimes{\cal S}$.}

\begin{tabularx}{1.\textwidth}{X}
\[ \begin{array}{|c|c|c|c|c|c|c|c|c|c|c|} \hline
\grey{Lat} &  \grey{\text Sch_{\cal S}}&   \grey{H{\text-}M_{\cal S}}  & \grey{\text orb_{\cal S}}& \grey{k_{\cal S}}  & \grey{t_{\cal S}}&\grey{ \ell_{\cal S}}&\grey{\text rows}
& \grey{\text {columns} -t_{\cal S}}&\grey{f_{\cal S }} \\
 \hline
\text{par} & \mathcal{C}_{1} & 1& \circ & 1 & 2&3&2v+1&   2v+3-2 &0\\
\hline
'' & \mathcal{C}_{2} &2 & 2222& 2 & 0&3&2v_{0}+1& 2v_{0}+3-0&2\\
\hline
'' & \mathcal{C}_{2v} &2m&2*22& 4 &0&2&2v_{0}+1& 2v_{0}+2-0&1\\
 \hline
\text{rect} &  \mathcal{C}_{s} & m&**& 2 & 1&2&2v_{0}+1& 2v_{0}+2-1&0\\
\hline
'' & \mathcal{C}_{2v} &2/m& *2222& 4 & 0&2&2v_{0}+1& 2v_{0}+2-0&1\\
\hline
 \text{square}&\mathcal{C}_{4} &4&442& 4 & 0&1&2v_{0}+1& 2v_{0}+1-0&0\\
\hline
'' &\mathcal{C}_{4v}&4m &*442 & 4 & 0& 1&2v_{0}+1& 2v_{0}+1-0&0\\
\hline
\end{array}  \]
\end{tabularx}

\label{table:PlaneFlex}
\end{table}
\subsection{Table of groups for the fixed lattice in $2$-dimensions}
\label{subsec:PlaneTableFixed}

As we mentioned earlier, it can be of interest to consider periodic frameworks where the lattice is fixed.
In the following table, each row will present the corresponding calculation for a given plane wall-paper group which is presented as
$\Z^{2}\rtimes \cal S$.  As above, this analysis does not include rows for the hexagonal tiling or groups which include glide reflections.

As in the previous section, in each row of Table \ref{table:PlaneFixed}, the calculation has several parts, each producing an integer.

\begin{table}[h!]
 \centering
 \caption{The added flexibility  induced by basic symmetries on a fixed 2-D lattice for $\Z^{2}\rtimes{\cal S}$.}

\begin{tabularx}{1.\textwidth}{X}
\[ \begin{array}{|c|c|c|c|c|c|c|c|c|} \hline
 \grey{Lat} &   \grey{\text Sch_{\cal S}} & \grey{H{\text-}M_{\cal S}} & \grey{\text orb_{\cal S}}& \grey{k_{\cal S}}  &\grey{t_{\cal S}}&
 \grey{\text {rows}}& \grey{ \text {columns} -{t_{\cal S}}}&\grey{f_{\cal S }}\\ \hline
\text{par} & \mathcal{C}_{1} & 1 & \circ & 1 & 2&2v-2 & 2v-2 &0\\
\hline
'' & \mathcal{C}_{2} & 2 & 2222& 2 & 0&2v_{0}-1& 2v_{0}&1\\
\hline
''  & \mathcal{C}_{2v} & 2m & 2*22& 4 &0&2v_{0}& 2v_{0}-0&0\\
\hline
\text{rect} & \mathcal{C}_{s} & m &**& 2 & 1&2v_{0}-1& 2v_{0}-1&0\\
\hline
\end{array}  \]
\end{tabularx}
\label{table:PlaneFixed}
\end{table}
Note that the number of rows (edge orbits), $e_0$,  is now such that   $k_{{\cal S}}e_{0}
\geq  2(k_{{\cal S}} v_{0}) -2$. This guarantees that we have at least the number of edges needed for a rigid periodic framework on the fixed lattice without symmetry.    This means we need to subtract a modified constant $c=\lfloor{\frac{2} {k_{\cal S}}}\rfloor$.  For Table \ref{table:PlaneFixed}, $c$ is $2$, $1$, or $0$.

Since for a fixed lattice, we clearly have $\ell_{\cal S}=0$ for each group $\mathcal{S}$, the corresponding column is omitted in Table \ref{table:PlaneFixed}.

Analogously to Table \ref{table:PlaneFlex}, the final column of Table \ref{table:PlaneFixed} shows the difference $f_{\cal S}=2v_{0}- t_{\cal S} - (2v_{0}-c)=c- t_{\cal S}$ which is the dimension of the guaranteed extra non-trivial motions of the symmetric framework, over the rigidity which the original count without symmetry promised.

Note that Table \ref{table:PlaneFixed} does not include all the point groups from Table \ref{table:PlaneFlex}. The groups we omitted only produce $0$'s in the last column.

\section{$3$-D periodic frameworks with symmetry: $\Z^{3}\rtimes \cal S$}
\label{sec:SpacePeriodicSymmetry}

We now apply the basic patterns of the previous sections to investigate the types of counts which arise for periodic structures with added symmetry in $3$-space.   As happened in the plane, these symmetries can have three impacts:
\begin{enumerate}
\item[(a)] the symmetry can restrict the possible shapes of the lattice cell or  equivalently,  the symmetry constraints leave a specific subset of $\ell_{\cal S}$ flexes of the lattice structure which preserve the desired symmetry.
\item[(b)]  the symmetry can block some, or all, of the translations of the lattice structure, altering the  basic count of $t_{\cal S}$;
\item[(c)] the symmetry determines the order of the group, that is, the size $k_{\cal S}$ of the orbits.
\end{enumerate}

\subsection{$\Z^{3}\rtimes \mathcal{C}_{i}$ - inversive symmetry in space}
\label{subsec:3DPeriodicInversive}
Consider the inversive symmetry in $3$-space with the center of symmetry at the origin.  This operation (which in the Schoenflies notation is called $i$) takes a joint $p$ to a joint $-p$.   In many tables of crystal symmetry, this symmetry operation is called central symmetry, and the crystals are called centrosymmetric.  All shapes of lattices are possible, and these fit into the triclinic lattice system (three angle choices).   In the Schoenflies notation, if inversion is the only non-trivial symmetry operation, the group is written as $\mathcal{C}_{i}$. In the Hermann-Mauguin notation, it is written as $\bar 1$, and in the orbifold notation, it is written as $1x$.   

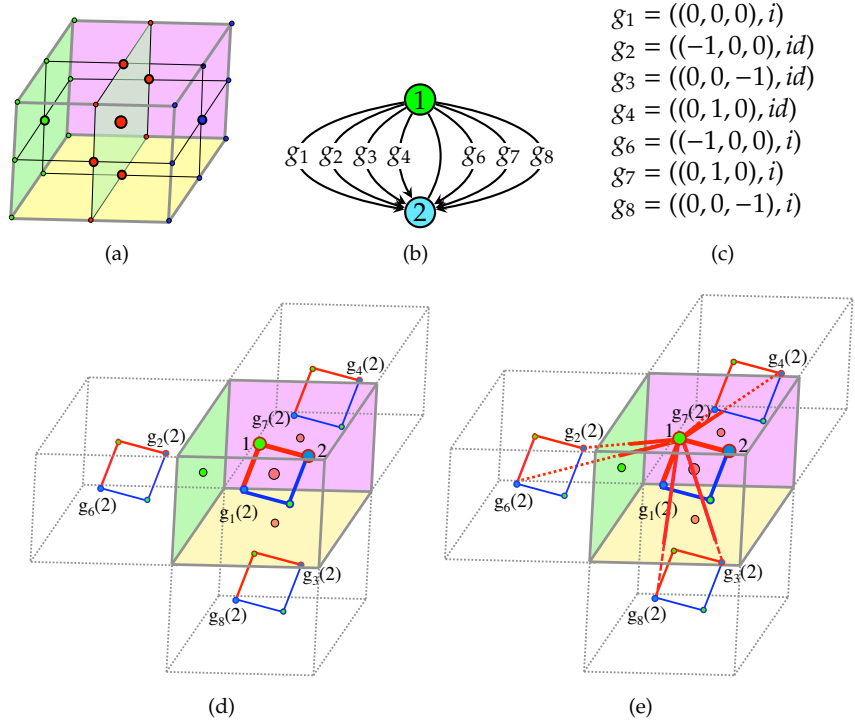
\begin{figure}[ht]
  \begin{center}
   \subfigure[] {\includegraphics [width=.25\textwidth]{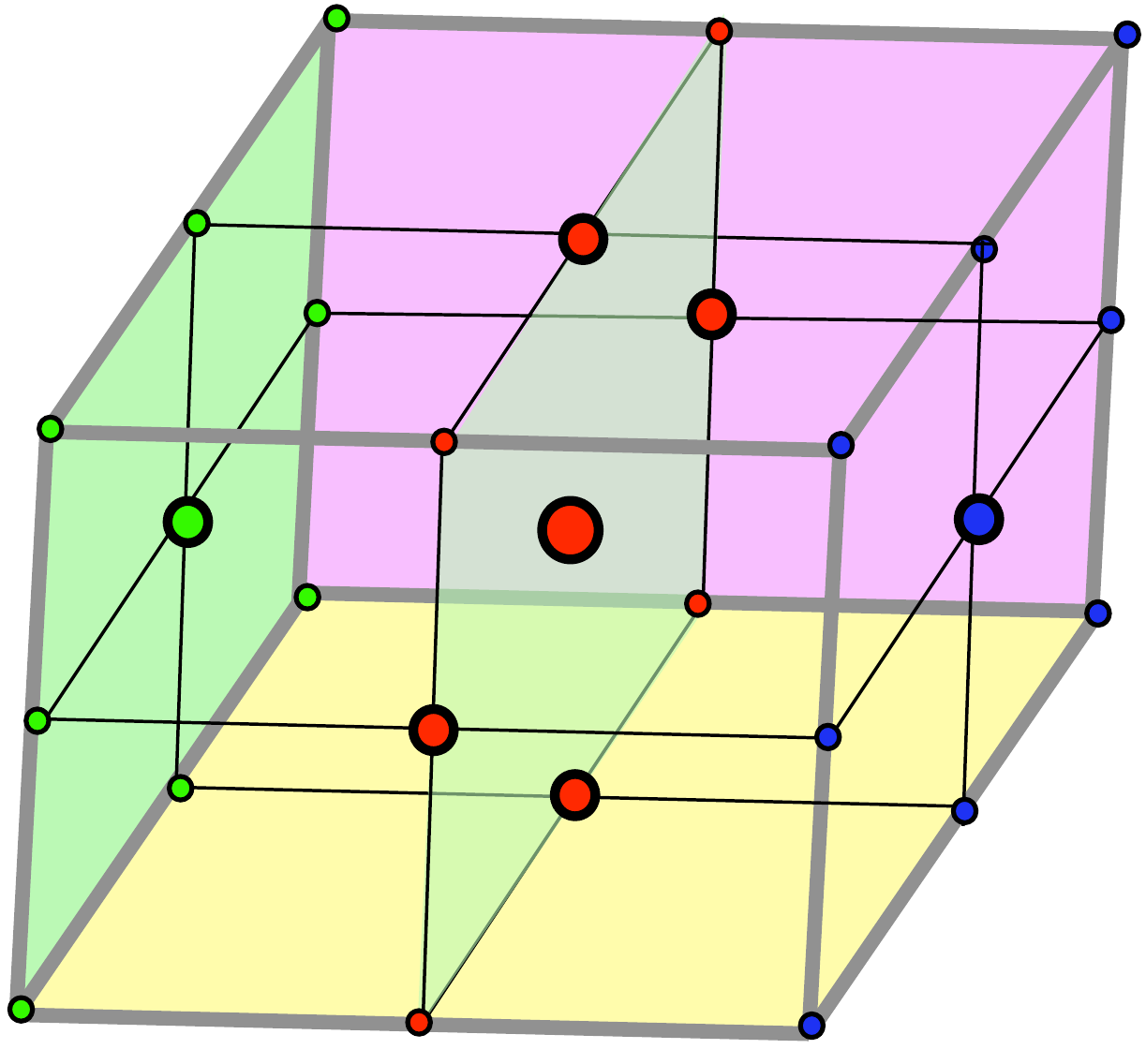}} \hspace{.1in}
       \subfigure[] {\begin{tikzpicture}[auto, node distance=2cm, thick, scale=.75]
\tikzstyle{vertex1}=[circle, draw, fill=green, inner sep=1pt, minimum width=3pt]; \tikzstyle{vertex2}=[circle, draw, fill=sky, inner sep=1pt, minimum width=3pt]; \tikzstyle{vertex3}=[circle, draw, fill=pink, inner sep=1pt, minimum width=3pt]; \tikzstyle{vertex4}=[circle, draw, fill=green, inner sep=1pt, minimum width=3pt]; \tikzstyle{gain} = [fill=white, text=black, inner sep =1pt,  anchor=center];
\node[vertex1] (1) {$1$}; \node[vertex2] (2) at (0, -2) {$2$};
\draw (1) edge [bend left] 
(2);
\pgfsetarrowsend{stealth}[ shorten >=2pt]
\draw (1) edge [bend right] node[gain] {$g_4$}
(2);
\draw (1) .. controls (-1.2,-.7) and (-1.2, -1.3)  .. node[gain] {$g_3$} (2);
\draw (1) .. controls (-2,-.6) and (-2, -1.4)  .. node[gain] {$g_2$} (2);
\draw (1) .. controls (-2.8,-.5) and (-2.8, -1.5)  .. node[gain] {$g_1$} (2);
\draw (1) .. controls (1.2,-.7) and (1.2, -1.3)  .. node[gain] {$g_6$} (2);
\draw (1) .. controls (2,-.6) and (2, -1.4)  .. node[gain] {$g_7$} (2);
\draw (1) .. controls (2.8,-.5) and (2.8, -1.5)  .. node[gain] {$g_8$} (2);
\pgfsetarrowsend{}

\end{tikzpicture}} 
\subfigure[] {\begin{tikzpicture}[auto, node distance=2cm, thick, scale=.75]
\tikzstyle{vertex1}=[circle, draw, fill=green, inner sep=1pt, minimum width=3pt]; \tikzstyle{vertex2}=[circle, draw, fill=sky, inner sep=1pt, minimum width=3pt]; \tikzstyle{vertex3}=[circle, draw, fill=pink, inner sep=1pt, minimum width=3pt]; \tikzstyle{vertex4}=[circle, draw, fill=green, inner sep=1pt, minimum width=3pt]; \tikzstyle{gain} = [fill=white, text=black, inner sep =1pt,  anchor=center];
\node at (-1.5, 0) {};
\node[text width = 1.3in] at (1,0)
{$g_1 = ((0,0,0), i)$\\
$g_2 = ((-1,0,0), id)$\\
$g_3 = ((0,0,-1), id)$\\
$g_4 = ((0,1,0), id)$\\
$g_6= ( (-1,0,0), i)$\\
$g_7 = ( (0,1,0), i)$\\
$g_8 = ((0,0,-1), i)$
};
\end{tikzpicture}}\hspace{.1in}
    \subfigure[] { \includegraphics[width=.44\textwidth]{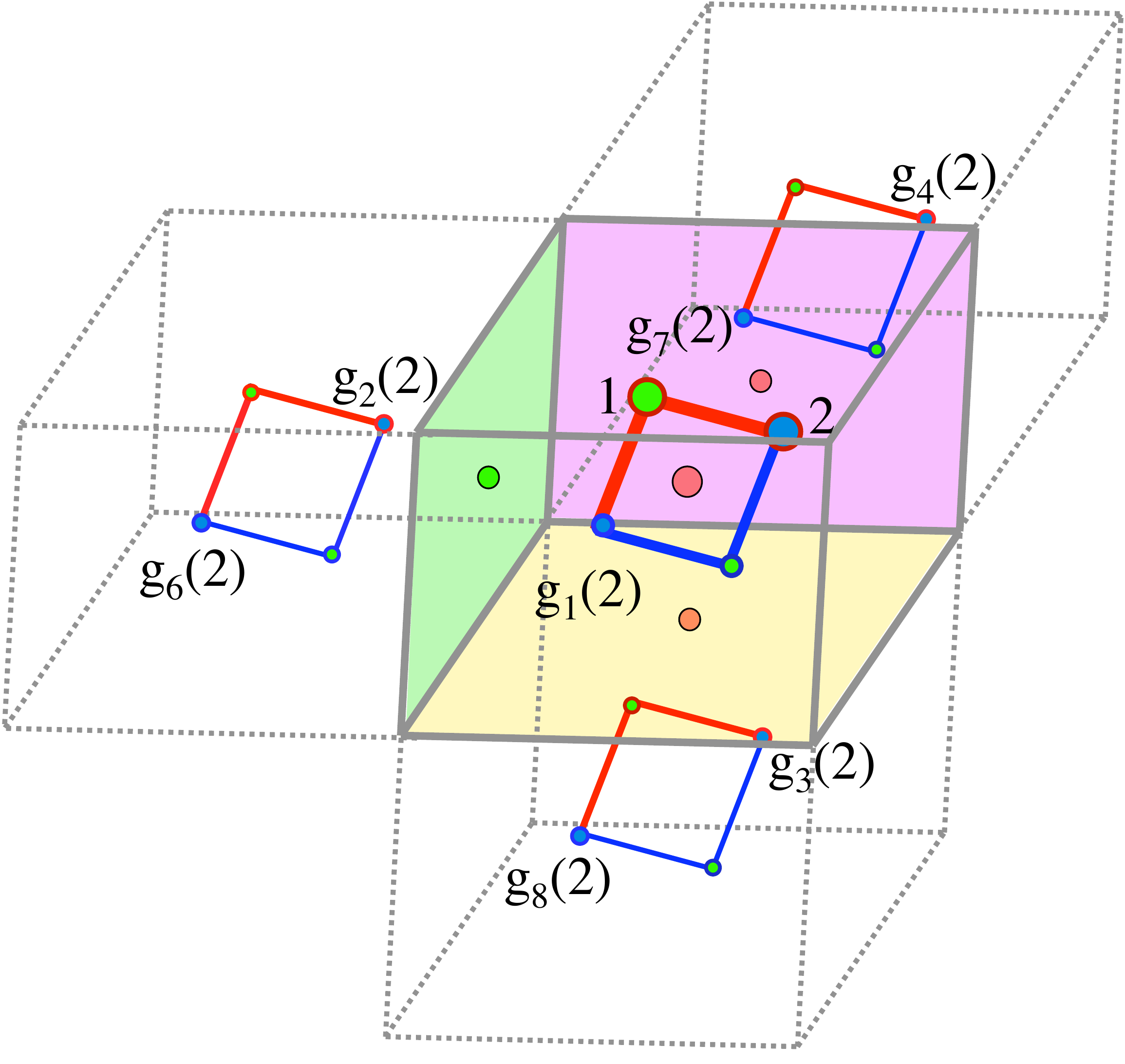}}
     \subfigure[] { \includegraphics[width=.45\textwidth]{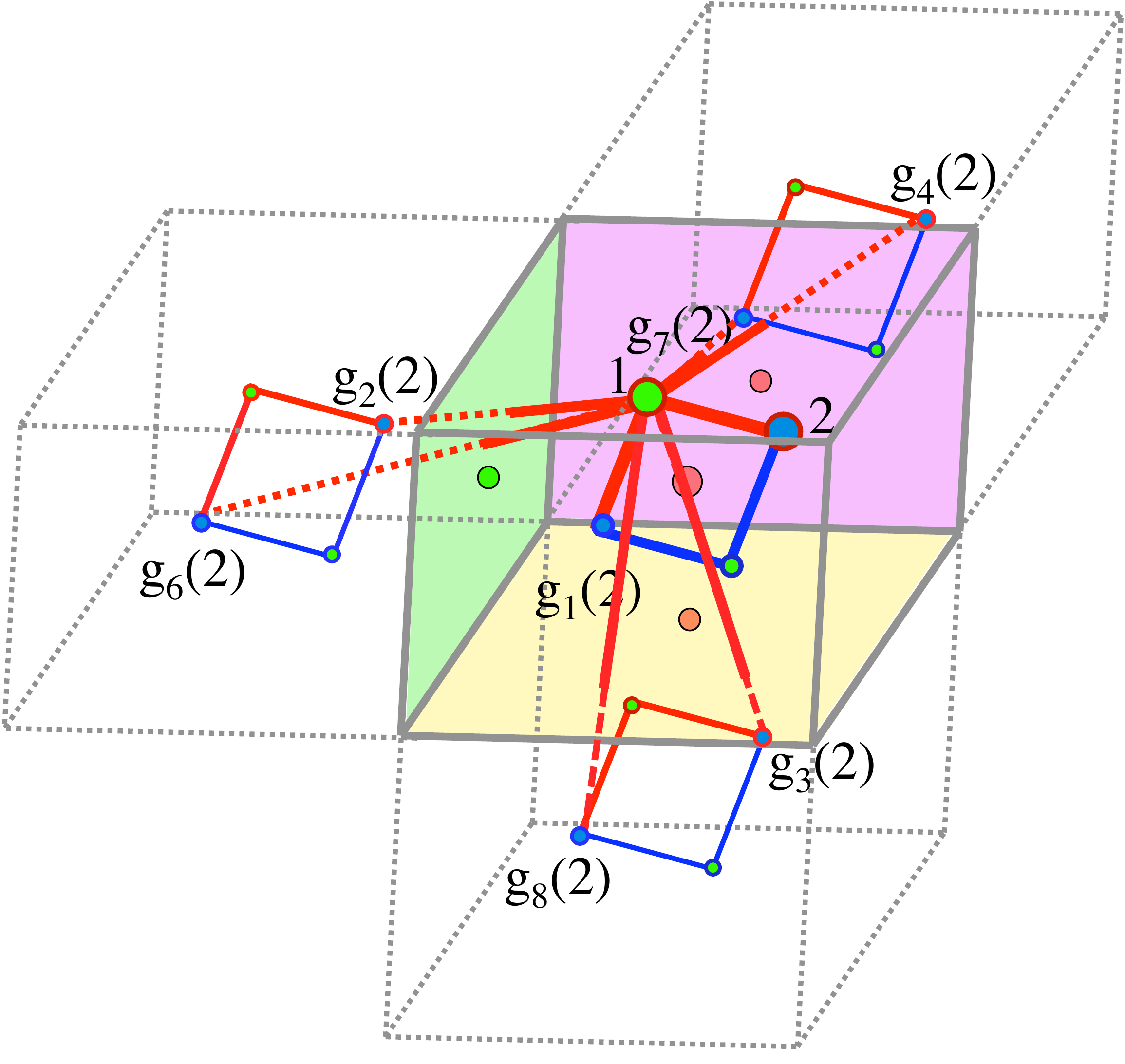}}
  \caption{In $3$-D, one center of inversion repeats with half the period (a). An orbit framework with $2$ orbits of vertices is shown in (b), with the group elements associated with the directed edges listed in (c). Parts (d) and (e) illustrate building up the corresponding symmetric-periodic  framework, moving from $2$ to $8$ orbits of edges (d). }
    \label{fig:3DInversion}
    \end{center}
    \end{figure}

As in the plane, if we have a center of inversion $c$, and a translation vector $t$ then there is another inversion centered at $c+{1\over 2}t$.
So, given the lattice of translations $\Z^{3}$ and one center of inversion at the origin, there is a full lattice of inversions, with translations ${1\over 2}\Z^{3}$, and the group of operations on the framework is written $\Z^{3}\rtimes \mathcal{C}_{i}$ (see also Figure~\ref{fig:3DInversion}(a)).

 \medskip

\noindent {\bf Example 7.1.1: fully flexible lattice $\Z^{3} \rtimes \mathcal{C}_{i}$.}
The necessary count for any minimally rigid non-symmetric periodic framework on the fully flexible lattice is $e=3v+3$ (recall Theorem \ref{thm:periodicMaxwell}). To permit inversion symmetry ($k_{\mathcal{C}_{i}}=2$) we will need to start with the shifted count $2e_{0} = 3(2v_{0})+4$ or equivalently $e_{0} = 3v_{0}+2$.

Since the full flexibility of the lattice fits with the inversive symmetry, we still have $\ell_{\mathcal{C}_{i}} = 6$. Further,  when we move to the symmetric periodic orbit matrix under inversive symmetry, all of the infinitesimal translations disappear from the kernel, so that $t_{\mathcal{C}_{i}}=0$. This gives rise to the symmetric Maxwell condition
$$e_0 \geq  3v_{0}+6 $$
for periodic rigidity.
The gap $e_{0} = 3v_{0}+2 \  < \ 3v_{0}+6$ implies that a graph which counted to be minimally rigid without symmetry, realized generically with inversive symmetry on a fully flexible lattice now has a space of (finite) flexes of dimension $4$.
\eop

\medskip

\noindent {\bf Example 7.1.2: fixed lattice $\Z^{3} \rtimes \mathcal{C}_{i}$.}
The necessary count for any minimally rigid non-symmetric periodic framework on the fixed lattice to be minimally rigid is
$e= 3v-3$. To permit inversive symmetry, we will need to start with the shifted count $2e_{0} = 3(2v_{0})-2$ or equivalently $e_{0} \ =\ 3v_{0}-1$.

Since we again have $t_{\mathcal{C}_{i}}=0$, we obtain the necessary symmetric Maxwell condition
$$e_0 \geq  3v_{0} $$
for periodic rigidity.
The gap $e_{0}= 3v_{0}-1\ <\ 3v_{0}$ predicts a non-trivial finite flex in generic realizations with inversive symmetry on the fixed lattice.
 \eop

As a summary, here is the impact of inversive symmetry for each of the variants of lattice flexibility introduced in Section \ref{subsec:PeriodicOrbitMatrices}:
\begin{table}[h!]
 \centering
 \caption{$3$-D lattice deformations with $\mathcal{C}_{i}$ symmetry.}
\begin{tabularx}{1.\textwidth}{X}
\[ \begin{array}
{|c|c|c|c|c|c|c|c|} \hline
 \grey{Lattice Def}&   \grey{\cal S}   & \grey{k_{\cal S}}  & \grey{t_{\cal S}}&\grey{\ell_{\cal S}}
 & \text{\grey {rows }}& \text{\grey {columns -$t_{\cal S}$}}& \grey{f_{\cal S }} \\ \hline
\text{flexible}& \mathcal{C}_{i} & 2&  0&6&e_{0}=3v_{0}+2  &  3v_{0}+6&4\\
\hline
\text{distortional}& \mathcal{C}_{i} & 2&  0&5&e_{0}=3v_{0}+1  &  3v_{0}+5&4\\
\hline
\text{scaling}& \mathcal{C}_{i} &  2& 0&3&e_{0}=3v_{0}   &  3v_{0}+3&3\\
\hline
\text{hydrostatic}& \mathcal{C}_{i} & 2&  0&1&e_{0}=3v_{0}-1  &  3v_{0}+1&2\\
\hline
\text{fixed}& \mathcal{C}_{i} & 2& 0&0&e_{0}=3v_{0}-1 &  3v_{0}&1\\
\hline
\end{array}  \]
\end{tabularx}

\label{table{inversive}}
\end{table}

\subsection{$\Z^{3}\rtimes \mathcal{C}_{2}$ and $\Z^{3}\rtimes \mathcal{C}_{s}$ - half-turn and mirror symmetry in space}
\label{subsec:3DPeriodicHalf}

Assume we have a 2-fold rotational axis along the $z$ direction.  This places the pattern into the monoclinic crystal system:  one face of the lattice is a parallelogram (perpendicular to the axis) and two faces are parallel to the axis and perpendicular to the parallelogram face.
  For this type of lattice, there are $4$ lattice parameters: the scale of each of the generating translations, and the one angle between the two generating translations of the parallelogram.


\medskip

\noindent {\bf Example 7.2.1: fully flexible lattice $\Z^{3} \rtimes \mathcal{C}_{2}$.}
With a fully flexible lattice, the necessary minimal number of edges for a periodic framework to be rigid and to be compatible with half-turn symmetry is $2e_{0} = 3(2v_{0}) +4$, or $e_{0} = 3v_{0} +2$.
In the orbit matrix, there are four columns corresponding to the lattice deformations, so the necessary symmetric Maxwell type count for periodic rigidity is
$$  e_{0} \geq 3v_{0} +4-1 = 3v_{0} +3. $$
Since we started with $e_{0} = 3v_{0} +2\ <\ 3v_{0} +3  $, we predict a non-trivial symmetry preserving finite flex for generic realizations with half-turn symmetry on the flexible lattice.  \eop

\medskip

\noindent {\bf Example 7.2.2: fixed lattice $\Z^{3} \rtimes \mathcal{C}_{2}$.}
With a fixed  lattice, the necessary minimal number of edges for a periodic framework to be rigid and to be compatible with half-turn symmetry is $2e_{0} \geq 3(2v_{0}) -2$, or $e_{0} \geq 3v_{0} -1$.  The necessary symmetric Maxwell type count for periodic rigidity on the fixed lattice is
$$  e_{0} \geq 3v_{0} -1. $$
Thus, we do not detect any added motions in this case.  \eop

In Table \ref{table{C2}} we present the $(\Z^{3} \rtimes\mathcal{C}_{2})$-symmetric Maxwell type counts for each type of lattice deformation.

\begin{table}[h!]
 \centering
 \caption{$3$-$D$ lattice deformations with $\mathcal{C}_{2}$ symmetry.
}
\begin{tabularx}{1.\textwidth}{X}
\[ \begin{array}
{|c|c|c|c|c|c|c|c|} \hline
 \grey{Lattice Def}&   \grey{\cal S}   & \grey{k_{\cal S}}  & \grey{t_{\cal S}}& \grey{\ell_{\cal S}}
 & \text{\grey {rows }}& \text{\grey {columns -$t_{\cal S}$}}& \grey{f_{\cal S }}
 \\ \hline
\text{flexible}& \mathcal{C}_{2} & 2&  1&4&e_{0}=3v_{0}+1  &  3v_{0}+4-1&2\\
\hline
\text{distortional}& \mathcal{C}_{2} & 2&  1&3&e_{0}=3v_{0}  &  3v_{0}+3-1&2\\
\hline
\text{scaling}& \mathcal{C}_{2} &  2& 1&3&e_{0}=3v_{0}   &  3v_{0}+3-1&2\\
\hline
\text{hydrostatic}& \mathcal{C}_{2} & 2& 1&1&e_{0}=3v_{0}-1  &  3v_{0}+1-1&1\\
\hline
\text{fixed}& \mathcal{C}_{2} & 2& 1&0&e_{0}=3v_{0}-1 &  3v_{0}-1&0\\
\hline
\end{array}  \]
\end{tabularx}
\label{table{C2}}
\end{table}

\medskip


Consider a periodic  framework in space with mirror symmetry. For this new group, there are only two key calculations to be done:
\begin{enumerate}
\item $t_{\mathcal{C}_{s}}=2$, since the two translations on directions within the mirror will (instantaneously) preserve the mirror.
\item $\ell_{\mathcal{C}_{s}}=4$.  Although there initially appear to be two alignments for the mirror: (i) parallel to two translation axes and perpendicular to another or (ii)  containing an axis of translation,  these turn out to be two variations of the same larger space tiling, and crystallographers only consider the first version.
 In this case we have an orthorhombic lattice system, and we have four parameters, $\ell_{\mathcal{C}_{s}}=4$.
 \end{enumerate}

\medskip

\noindent {\bf Example 7.2.3: fully flexible lattice $\Z^{3} \rtimes \mathcal{C}_{s}$.}
As before, we start with the following initial count without symmetry: $2e_{0} = 3(2v_{0}) +4$, or $e_{0} = 3v_{0} +2$.
From the periodic symmetric orbit matrix we obtain the following necessary symmetric Maxwell type count for periodic rigidity: $$e_{0} \geq 3v_{0}+4-2 = 3v_{0} +2.$$
This suggests that there is no additional flexibility in the structure when mirror symmetry is added.  \eop

It turns out that for mirror symmetry, all of the variants of lattice deformations produce no added motions.

\subsection{Table of groups for the fully flexible lattice in $3$-dimensions}
\label{subsec:3DPeriodicTableFlex}

Following the process illustrated in the previous examples, we can track the necessary increases in flexibility which follow from minimal generically rigid periodic frameworks for various symmetry groups $\Z^{3}\rtimes \mathcal{S}$ in $3$-space.
As before, this does not include rows for the groups  with $6$-fold rotational symmetry, or any patterns with glide reflections. They will require some significant modifications of the simple pattern presented here.

Analogous to the tables in Sections \ref{subsec:PlaneTableFlex} and \ref{subsec:PlaneTableFixed}, in each row of Table \ref{table:3DFlex}, the calculation has several parts - each producing an integer.

The number of rows (edge orbits), $e_0$,  is such that   $k_{{\cal S}}e_{0}\geq  3(k_{{\cal S}} v_{0}) +3$, which guarantees that we have at least the number of edges needed for a rigid periodic framework without symmetry.    This means we need to add a modified constant $c=\lceil{\frac{3} {k_{\cal S}}}\rceil$.  For Table \ref{table:3DFlex}, $c=3$ for  $k_{{\cal S}}=1$, $c=2$ for  $k_{{\cal S}}=2$, and $c=1$ for all bigger orbit sizes.

As usual, $t_{\cal S}$ is the dimension of the space of translations contained in the symmetry element of ${\cal S}$.  This will be $t_{\cal S}=3$ for the identity group, $t_{\cal S}=2$ for a single mirror,  $t_{\cal S}=1$ for a single rotation (with or without a mirror along the axis), and $t_{\cal S}=0$ if only a point is fixed.

In Table \ref{table:3DFlex} we compare the number of rows, $e_0$, with the number of columns minus $t_{\cal S}$, $3v_{0}+\ell_{\cal S} - t_{\cal S}$; the difference $f_{\cal S}=3v_{0}+\ell_{\cal S} - t_{\cal S} - (3v_{0}+c)=\ell_{\cal S} - t_{\cal S} -c$  is the dimension of the guaranteed extra non-trivial motions of the symmetric framework over the rigidity which the original count without symmetry promised.


\begin{table}[h!]
 \centering
 \caption{The added flexibility induced by within basic symmetries on a fully flexible $3$-D lattice for $\Z^{3}\rtimes{\cal S}$.}
\begin{tabularx}{1.\textwidth}{X}
\[ \begin{array}{|c|c|c|c|c|c|c|c|c|c|} \hline
 \grey{Lat.\ System }&   \grey{Sch_{\cal S}} & \grey{H{\text-}M_{\cal S}} & \grey{orb_{\cal S}}&  \grey{k_{\cal S}}  & \grey{t_{\cal S}}&  {\grey{\ell_{{\cal S}}}} & \text{\grey {rows}}& \text{\grey {columns -$t_{\cal S}$}}& \grey{f_{\cal S }} \\ \hline
\text{triclinic}& \mathcal{C}_{1}  & 1& 11& 1 & 3&6&3v+3 & 3v+6-3&0\\
\hline
\text{''}&\mathcal{C}_{i}  & \bar{1}& 1x& 2 & 0&6&3v_{0}+2 & 3v_{0}+6-0&4\\
\hline
\text{monoclinic}& \mathcal{C}_{2} & 2& 22& 2 & 1&4&3v_{0}+2 & 3v_{0}+4-1&1\\
\hline
\text{''}& \mathcal{C}_{s} &m& 1* & 2& 2&4&3v_{0}+2 & 3v_{0}+4-2&0\\
\hline
\text{''}& \mathcal{C}_{2h} & 2/m &2*& 4 & 0&4&3v_{0}+1 & 3v_{0}+4-1&2\\
\hline
\text{orthorhom}& \mathcal{C}_{2v} &222  &*22& 4 & 1&3&3v_{0}+1 & 3v_{0}+3-1&1\\
\hline
\text{''}& \mathcal{D}_{2} &mm2  &222& 4 & 0&3&3v_{0}+1 & 3v_{0}+3-0&2\\
\hline
\text{''}& \mathcal{D}_{2h} &mmm  &*222& 8 & 0&3&3v_{0}+1 & 3v_{0}+3-0&2\\
\hline
\text{tetragonal}& \mathcal{C}_{4} & 4 & 44 & 4 & 1&2&3v_{0}+1 & 3v_{0}+2-1&0\\
\hline
\text{''}& \mathcal{S}_{4} & \bar{2}& 2x &4 & 0&2&3v_{0}+1 & 3v_{0}+2-0&1\\
\hline
\text{''}& \mathcal{C}_{4h} & 4/m &4* & 8 & 0&2&3v_{0}+1 &  3v_{0}+2-0&1\\
\hline
\text{''}& \mathcal{C}_{4v} &4mm  &*44 & 8 & 1&2&3v_{0}+1 &  3v_{0}+2-1&0\\
\hline
\text{''}& \mathcal{D}_{2d} &\bar{4}2m  &2*2& 8 & 0&2&3v_{0}+1 & 3v_{0}+2-0&1\\
\hline
\text{''}& \mathcal{D}_{4} & 422 &422& 8 & 0&2&3v_{0}+1  &  3v_{0}+2-0&1\\
\hline
\text{''}& \mathcal{D}_{4h} &4/mmm  &*422& 16 & 0&2&3v_{0}+1 &  3v_{0}+2-0&1\\
\hline
\text{trigonal}&\mathcal{C}_{3} &   3 &33& 3 & 1&2&3v_{0}+1 & 3v_{0}+2-1&0\\
\hline
\text{''}& \mathcal{S}_{6} &\bar3   &3x& 6 & 0 &2&3v_{0}+1 & 3v_{0}+2-0&1 \\
\hline
\text{''}& \mathcal{D}_{3} &32   &322& 6 & 0 &2&3v_{0}+1 & 3v_{0}+2-0&1 \\
\hline
\text{''}& \mathcal{C}_{3v} &3m   &*33& 6 & 1&2&3v_{0}+1 & 3v_{0}+2-1&0\\
\hline
\text{''}& \mathcal{D}_{3d} &\bar3m   &2*3& 12 & 0&2&3v_{0}+1 & 3v_{0}+2-0&1\\
\hline
\text{cubic}& \mathcal{T} & 23  &332& 12 & 0&1&3v_{0}+1 & 3v_{0}+1-0&0\\
\hline
\text{''}& \mathcal{T}_{h} & m\bar3   &3*2& 24 & 0&1&3v_{0}+1 & 3v_{0}+1-0&0\\
\hline
\text{''}& \mathcal{T}_{d} &\bar43m   &*332& 24 & 0&1&3v_{0}+1 & 3v_{0}+1-0&0\\
\hline
\text{''}& \mathcal{O} &432   &432& 24 & 0&1&3v_{0}+1 & 3v_{0}+1-0&0\\
\hline
\text{''}& \mathcal{O}_{h} &m\bar3 m   &*432& 48 & 0&1&3v_{0}+1 & 3v_{0}+1-0&0\\
\hline
\end{array}  \]
\end{tabularx}
\label{table:3DFlex}
\end{table}

\subsection{Table of groups for the fixed lattice in $3$-dimensions}
\label{subsec:3DPeriodicTableFixed}

In Table \ref{table:3DFixed}, we track the necessary increases in flexibility which follow from minimal generically rigid periodic frameworks on a fixed lattice for various symmetry groups in $3$-space. This analysis is analogous to the one in the previous section.  We simply remove the column for $\ell_{\cal S}$ which is always $0$, and work with the modified counts.

The entries $(-1)$ in Table \ref{table:3DFixed} indicate that, for this group, the symmetry guarantees that there is a symmetric self-stress in the symmetric framework (see also Section \ref{subsec:Stresses}).  Because the patterns of $0$ and occasional $(-1)$ become clear quickly, we do not fill in all rows of the matrix.


\begin{table}[h!]
 \centering
 \caption{The added flexibility induced by symmetries on a fixed 3-D lattice for $\Z^{3}\rtimes{\cal S}$.}
\begin{tabularx}{1.\textwidth}{X}
\[ \begin{array}{|c|c|c|c|c|c|c|c|c|c|c|} \hline
 \grey{Lat.\ System }&   \grey{Sch_{\cal S}} & \grey{H{\text-}M_{\cal S}} & \grey{orb_{\cal S}}&  \grey{k_{\cal S}}  & \grey{t_{{\cal S}}}
 &  \text{\grey {rows}}&  \text{\grey {columns -$t_{{\cal S}}$}}& \grey{f_{\cal S }}
 \\ \hline
\text{triclinic}& \mathcal{C}_{1} &1 & 11& 1 & 3&3v-3 & 3v-3&0\\
\hline
\text{''} & \mathcal{C}_{i} & \bar{1} & 1x& 2 & 0&3v_{0}-1 & 3v_{0}-0&1\\
\hline
\text{monoclinic}&\mathcal{C}_{2} & 2 & 22& 2 &1 &3v_{0}-1 & 3v_{0}-1&0\\
\hline
\text{''}& \mathcal{C}_{s} & m &1*& 2 & 2&3v_{0}-1 & 3v_{0}-2& (-1)\\
\hline
\text{''}& \mathcal{C}_{2h} & 2/m &2*& 4 & 0&3v_{0} & 3v_{0}-0&0\\
\hline
\text{orthorhomb.}& \mathcal{C}_{2v} & 222 &*22& 4 & 1&3v_{0} & 3v_{0}-1&(-1)\\
\hline
\text{''}& \mathcal{D}_{2} & mm2 &222& 4 & 0 &3v_{0} & 3v_{0}-0&0\\
\hline
\text{''}& \mathcal{D}_{2h} & mmm &*222& 8 & 0 &3v_{0} & 3v_{0}-0&0\\
\hline
\text{tetragonal}& \mathcal{C}_{4} & 4&44 & 4 & 1&3v_{0} & 3v_{0}-1&(-1)\\
\hline
\text{''}& \mathcal{S}_{4} & \bar{2} &2x & 4 & 0&3v_{0} & 3v_{0}-0&0\\
\hline
\end{array}  \]
\end{tabularx}
\label{table:3DFixed}
\end{table}


\section{Extensions and further observations}
\label{sec:Further}

\subsection{Adjusting for fixed joints and bars} \label{subsec:loops}

To keep the analysis simpler, up to this point, we have not considered joints or bars which are fixed by a non-trivial symmetry operation in $\cal S$.  This was for simplicity of our counts (keeping a fixed $k_{\cal{S}}$ for all vertex and edge orbits), but incorporating these modifications are not a barrier to the final analysis.
The methods presented here also apply to symmetric periodic frameworks with fixed structural components, as has been shown in the previous work on finite symmetric frameworks \citep{BSWWorbit}. We give a brief indication of how this can work, in two cases.

\medskip

We first consider the case where the given periodic framework has a bar (but no joint) which is fixed by a non-trivial symmetry operation in the group $\cal S$. In this case, the comparison of counts needs to be adjusted, since there are now edge orbits of different sizes - the size of an edge orbit corresponding to a bar which is fixed by a non-trivial symmetry operation will no longer be equal to the order of the group $\cal S$. This is illustrated in the following example.

\medskip

\noindent{\bf Example 8.1.1}  The necessary count for any minimally rigid non-symmetric periodic framework on the fully flexible lattice is $e=3v+3$. If a framework with $\Z^{3}\rtimes \mathcal{C}_{i}$ symmetry has exactly one bar  which is fixed by the inversion, then we have $e= 2e_{0}-1$. Thus, in this case, we do not need to shift the count $e=3v+3$ to allow for inversive symmetry. We have $ 2e_{0}-1 = 3(2v_{0})+3$, and hence $e_0=3v_{0}+2 $. The symmetric Maxwell condition for periodic rigidity is $e_0 \geq 3v_{0} +6$. So we detect $4$ added degrees of flexibility. \eop

The method for adjusting all of the table entries, for any number of fixed bars, should now be accessible to the reader.

\medskip

Similarly, for a periodic framework with symmetry $\Z^{d}\rtimes \cal S$, there can also be \emph{joints} in the unit cell which are fixed by some non-trivial symmetry operations in the group $\cal S$. For these, the methods applied in \citet{BSWWorbit} also immediately transfer. Note, however, that if there exist \emph{joints} that are fixed by non-trivial symmetry operations in $\cal S$, we may not only have vertex orbits of different sizes, but the sets of columns corresponding to the vertices in the orbit matrix may now also be of varying size.

Depending on what subspace the joint is now restricted to,  the number of corresponding columns in the orbit matrix is reduced.  If a joint of a $3$-dimensional structure, for example, is restricted to a mirror plane, the number of columns will be reduced to $2$.  If it is restricted to a line (e.g. the intersection of two mirrors, or an axis of rotation) there will only be $1$ column, and if it is restricted to a point (e.g. the center of inversion, or the intersection of a mirror and an axis) the number of columns will be $0$.

We give two samples to illustrate this. In both cases we assume that there are no bars which are fixed by a non-trivial symmetry operation in $\cal S$.

\medskip

\noindent{\bf Example 8.1.2} The necessary count for any minimally rigid non-symmetric periodic framework on the fully flexible lattice is $e=3v+3$. If a framework with $\Z^{3}\rtimes \mathcal{C}_{s}$ symmetry has exactly one joint which is fixed by the reflection in $\mathcal{C}_s$, then we have $v= 2v_{0}-1$. Thus, in this case, we do not need to shift the count $e=3v+3$ to allow for mirror symmetry. We have
$e= 2e_{0}= 3(2v_{0}-1)+3=6v_0$, and hence $e_0=3v_{0}$. The symmetric Maxwell condition for periodic rigidity is $$e_0 \geq 3(v_{0}-1)+2+4-2=3v_0+1.$$ So we detect one added degree of flexibility. \eop

\medskip

\noindent{\bf Example 8.1.3} As noted above, the necessary count for any minimally rigid non-symmetric periodic framework on the fully flexible lattice is $e=3v+3$. Recall from Example 7.2.3 that we need to shift this count to $e=3v+4$ if we want to permit $\Z^{3}\rtimes \mathcal{C}_{s}$ symmetry, where no joint and no bar is fixed by the reflection in $\mathcal{C}_s$. No added motion was detected in this case. We now start with a framework whose underlying periodic graph satisfies the count $e=3v+4$, and we further assume that exactly two of the joints are fixed by the reflection in $\mathcal{C}_s$ (and hence lie on two distinct points on the corresponding mirror plane). Then we have $v=2v_0-2$. Thus, we have
$e= 2e_{0}= 3(2v_{0}-2)+4=6v_0-2$, and hence $e_0=3v_{0}-1$. The symmetric Maxwell condition for periodic rigidity is $$e_0 \geq 3(v_{0}-2)+2(2)+4-2=3v_0.$$ So we now detect one added degree of flexibility, although the framework is overbraced, generically. \eop

\subsection{Other groups}
\label{subsec:OtherGroups}

Given the wide array of lattice systems and space groups which occur in crystals, we have only analyzed some of the possibilities.  The groups we considered were restricted to the form $\Z^{d}\rtimes \cal S$.

It is a small change to include groups with glide reflections as generators which are not semi-direct products.  Such a group does include $\Z^{d}$ as a subgroup, so it admits a representation as a periodic orbit graph.  The glide reflection acts within the periodic orbit graph,  and the same counts and methods presented here can be applied (see also Example 8.1.3). Note, however, that the group acting on the periodic orbit framework is not a point group (i.e., a group of isometries leaving a point fixed),  so it is technically distinct from the groups we have considered so far.

\medskip

\noindent{\bf Example 8.2.1}
We start with the necessary count for minimal rigidity of a non-symmetric periodic framework on the fully flexible lattice: $e=3v +3$. This needs to be shifted to $2e_{0}=3(2v _{0})+ 4$ or $e_0=3v_0+2$ to permit glide reflectional symmetry. As in the case of $\Z^{3}\rtimes \mathcal{C}_s$, there are $4$ parameters for the lattice deformations in the orbit matrix, and there are still two translations which preserve the glide reflection, along with the periodic lattice. This leads to the following symmetric Maxwell condition for periodic rigidity:  $e_{0} \geq 3v_{0}+ 4-2=3v_0+2$. So, as in the case of $\Z^{3}\rtimes \mathcal{C}_s$, we do not detect any flexes with these counts.

Note that placing a joint on the mirror plane of the glide reflection does not mean that this joint is fixed by the glide reflection, but it does restrict the number of degrees of freedom for this joint - that is, the number of columns corresponding to this joint in the orbit matrix - to 2.   Thus, shifting two joints onto the glide plane will generate an added degree of flexibility, as it did for a single mirror (recall Example 8.1.3).
\eop

It is also simple to include a single screw symmetry in 3-space, such as a $2$-fold screw.  This will have the same impact as the $2$-fold rotation, generating an additional flexibility in the flexible lattice.   However, we have not completed the analysis for more general mixtures of glide reflections, and screw symmetries, with the various compatible lattices and groups.  We anticipate that it will be possible to extend all of the results to the space groups which share the underling lattice systems we explored here.

We have also not analyzed the symmetries which include $6$-fold rotations in the plane or in $3$-space - those related to the triangular lattice.  These triangle groups require some additional care with the gain graphs and the orbit counts when they are placed into the unit cells as we factor out the  translations. Again, we anticipate these lines of the tables will fall into place once this analysis is completed.  In short, we anticipate that it will be possible to give analogous orbit matrices and counts for all of the plane wallpaper groups and all of the 230 space groups.

However, the class already covered includes a number of key examples, and these example do illustrate that symmetry can add flexibility to structures which otherwise count to be generically rigid.


\subsection{Higher dimensions}
\label{subsec:HigherDim}

The basic results for rigidity matrices for periodic structures and for orbit rigidity matrices for symmetric structures all extend to arbitrary dimensions \citep{BorceaStreinuI,
BorceaStreinuII,
BSWWorbit}.

All of the techniques for combined analyses also extend to arbitrary dimensions $d$ for  groups of the form  $\Z^{d}\rtimes \cal S $.  As an example, we summarize the extension for the groups $\Z^{d}\rtimes \mathcal{C}_{i} $ and $\Z^{d}\rtimes \mathcal{C}_{s} $.

For (non-symmetric) generic rigidity on a  fully flexible lattice in $d$-space, a necessary condition is
$ e \geq dv + {d\choose 2}  $ (recall Theorem \ref{thm:periodicMaxwell}). In both of the following examples, we will assume that we have the minimal number of edges for the graph to be rigid on a fully flexible lattice and to be compatible with the given symmetry, with no joint and no bar fixed by a non-trivial symmetry operation.


\medskip

\noindent{\bf Example 8.3.1} Consider a periodic framework with inversive symmetry in dimension $d$ and the associated orbit matrix.  We have
$ 2e_{0} = d(2v_{0}) + {d\choose 2} $ if ${d\choose 2} $ is even, and $ 2e_{0} = d(2v_{0}) + {d\choose 2} +1$ otherwise. Further, we have the following column count for the orbit matrix:  $dv_{0} + {d+1\choose 2} $.  There are clearly no residual translations which preserve the inversive symmetry along with the periodic lattice. So, since
$$ e_0= dv_{0} + \left\lceil{{d\choose 2} / 2} \right\rceil < dv_{0} + {d+1\choose 2},$$
we detect $f_{\mathcal{C}_i}={d+1\choose 2}-\left\lceil {{d\choose 2}/ 2}\right\rceil $ degrees of flexibility. This is a number which is growing as $ \lfloor {d\choose 2}/2\rfloor +d$,  which is quadratic. \eop

\medskip

\noindent{\bf Example 8.3.2} For mirror symmetry, we also have a pattern across the dimensions.   In dimension $d$, there will be one lattice vector perpendicular to the mirror, and the rest will be parallel to the mirror.  This removes a set of $d-1$ angles as lattice variables, leaving $\ell_{\mathcal{C}_s}={d\choose 2}+1$.  There are also $d-1$ residual translations. So, since
$$e_{0} = dv_{0} + \left\lceil{{d\choose 2} / 2} \right\rceil \ < \ dv_{0} + {d\choose 2}+1-(d-1),$$
we detect $f_{\mathcal{C}_{i}}=  {d\choose 2} +1 -(d-1) -\left\lceil {{d\choose 2}/ 2}\right\rceil $ degrees of flexibility.  This simplifies as $f_{\mathcal{C}_{i}}=  \lfloor {d\choose 2}/2\rfloor +2 - d $.  For $d>3$ this is positive and growing quadratically. \eop

Similar methods can be used to analyze a number of other groups of the form  $\Z^{d}\rtimes \cal S $ in higher dimensions, provided we have worked out the flexibility of the lattices which support the symmetry, in that space.

\subsection{Topology of orbit frameworks}\label{subsec:Topology}
We may regard orbit frameworks as graphs embedded on appropriately chosen orbifolds (generalizations of manifolds, which locally resemble Euclidean space). This orbifold is defined by the original setting of the framework ($\mathbb R^2$ or $\mathbb R^3$), modulo the symmetry group. For example, periodic frameworks have symmetry group $\mathbb Z^d$, and may be viewed as graphs on the $d$-dimensional topological torus $\mathbb R^d / \mathbb Z^d$. We can think of this as the unit cell with pairs of opposite faces identified. Similarly, a plane framework with $\mathcal C_n$ symmetry can be regarded as a framework on a cone, with opening angle $2\pi/n $.  The orbit matrix now also provides conditions for rigidity and flexibility for the (orbit) frameworks on this surface.

For periodic frameworks with additional symmetry, the underlying orbifold may be more exotic. For example, periodic frameworks with mirror symmetry in the plane or space (given by the groups $\mathbb Z^2 \rtimes \mathcal C_{\cal S}$ and $\mathbb Z^3 \rtimes \mathcal C_{\cal S}$ respectively) correspond to frameworks on  $2$- or $3$-spheres $\mathbb S^2$ and $\mathbb S^3$, but with a flat metric.  Frameworks in $3$-space with inversive symmetry ($\mathbb Z^3 \rtimes \mathcal C_{i}$) have an orbifold with topology of $\mathbb P^3$, projective $3$-space.  Similar statements are possible for all frameworks which admit an orbit framework under the action of their symmetry group.  Again, the periodic symmetric orbit matrices represent the rigidity matrices for frameworks actually living in these more exotic spaces, with flat metrics.  The results here give some necessary conditions for rigidity on these orbifold surfaces.

As an additional topological and geometric layer, the more detailed studies of the orbifolds for space groups describe a number of these orbifolds as fiberfolds \citep{Conway}.  These are essentially `fibered prisms' over the plane orbifold.  For a group $\cal{S}$ in the plane, we have the stretched group $\bar{\cal{S}}$ such that $\bar{s}(x,y,z) = (s(x,y),z)$ for each group element.  The simplest extension of our counts occurs for the scaling lattice deformations.   In this setting, we have $k_{\bar{\cal{S}}}=k_{\cal{S}}$,
$t_{\bar{\cal{S}}}=t_{\cal{S}}+1$ and $\ell_{\bar{\cal{S}}}=\ell_{\cal{S}}+1$.  Comparing $e_{0}=2v_{0}$ for the plane without symmetry with the scaling counts $e_{0}=2v_{0}+\ell_{\cal{S}} -t_{\cal{S}}$, we have a flexibility of $\ell_{\cal{S}} -t_{\cal{S}}$.  Doing the same calculation in $3$-space, we find  $\ell_{\bar{\cal{S}}}  -t_{\bar{\cal{S}}}=\ell_{{\cal{S}}}  - t_{{\cal{S}}} $. We conclude that the flexibility induced by $ \cal{S}$ and by $\bar{\cal{S}}$ is the same.

In general, it is an open problem to investigate this connection across dimensions more thoroughly.  The  examples in \S\ref{subsec:HigherDim} give a hint of how some of this might work.
It is also an  interesting problem to  predict how flexibility of  a periodic symmetric graph in one dimension connects to flexibility for some associated periodic symmetric graph in the next dimension up.

\subsection{Statics and stresses}\label{subsec:Stresses}
In this paper, we have not focused on stresses - row dependencies of the (periodic) rigidity matrices, or symmetric stresses - row dependencies of the orbit matrices \citep{BSWWorbit, WW, W1}.  Imposing symmetry on a periodic structure, however, may not only increase the flexibility of the structure, but it can also increase the dimension of the space of self-stresses of the structure (recall Table \ref{table:3DFixed}, for example). There are essentially three possible sources for these additional self-stresses.

As a first source of stress, in a number of cases, we had to have additional edges to ensure the count of edges was divisible by $k_{\cal S}$.  This of course leads to additional row dependencies of the (periodic) rigidity matrix which can be detected by simply counting the number of vertices and edges of the underlying periodic orbit graph of the structure. These stresses may or may not be symmetric.

Secondly, where introducing symmetry takes a minimally rigid periodic framework to an infinitesimally flexible framework this is a guarantee that the symmetric periodic framework has a space of self-stresses of dimension equal to the dimension of the added infinitesimal flexes. Note that if we detect a \emph{symmetric} infinitesimal flex in the periodic framework from the orbit counts, then the framework must have an asymmetric self-stress. Similarly, if we detect a \emph{symmetric} self-stress in the periodic framework from the orbit counts, then the framework must have an asymmetric infinitesimal motion. This is because the orbit counts cannot detect the presence of `paired' symmetric infinitesimal flexes and self stresses.


The third possibility is that although our  orbit counts do not detect any self-stresses, the rows of the orbit matrix are not actually independent.
To determine this requires additional direct computation of the rank of the orbit matrix.

The classic plane example of Bottema's mechanism -  $K_{4,4}$ with symmetry group $\mathcal{C}_{2v}$ in the plane - is known to be flexible not from the simple orbit counts, but from these counts plus the added information that the orbit matrix has a row dependency \citep{BSWWorbit}.
One can make periodic structures from this example, and create added flexes due to the periodic symmetry, which are not predicted simply by the counts presented here for the rows and columns of the orbit matrix.

\subsection{Sufficient conditions for rigidity or independence with symmetry}
\label{subsec:Sufficient}
In the classical work for finite frameworks, there are results where combinatorial properties on a graph, and its subgraphs, guarantee that the rigidity matrix is full rank, or independent, at generic configurations.  The most famous example in the plane is Laman's theorem \citep{Laman, W1}.   In $3$-space, some partial sufficiency results come from graphs constructed by specific types of inductive constructions, or planar graphs which are $3$-connected \citep{W1}.

One can ask about combinatorial properties which guarantee that a periodic orbit matrix  has full rank.  This has been determined in the plane for both the flexible lattice  \citep{Theran}, and the fixed lattice \citep{Ross}.  In addition, for some plane symmetry groups, there are combinatorial characterizations of graphs which are minimal rigid and maximal independent frameworks at symmetric generic configurations, as well as conjectures for other cases \citep{BS3}.  It is then a natural question to be explored to find combinatorial characterizations of graphs which are minimal periodic rigid with a given symmetry, in the plane.  This may be accessible for at least some cases.

One can also ask for characterizations of special classes of symmetric periodic rigid graphs in $3$-D, for example, frameworks built by connecting symmetric convex polyhedra in a periodic form, or frameworks built by inductive constructions based on edge and vertex splitting.

In summary,  our initial exploration presented here raises a large number of interesting unsolved problems, and areas for further research.  We invite the reader to join in this exploration.

\section*{Acknowledgements}

The authors thank Simon Guest, Stephen Power, and Michael Thorpe for helpful conversations and correspondence.  Their assistance with context and current literature was invaluable.  We also thank Gerald Audette for assistance sorting out the vocabulary and notation in crystallography. Finally we are grateful to the anonymous referees for numerous useful comments and suggestions.

\bibliographystyle{humanbio}
\bibliography{SymmetryPeriodic_Rev}
\end{document}